\documentclass{amsart}
\usepackage{epsfig}
\usepackage{graphicx}
\usepackage{eepic}
\usepackage{amssymb,latexsym,ifthen}
\newtheorem{Def}{Definition}

\newtheorem{Thm}{Theorem}

\newtheorem{Prop}{Proposition}
\newtheorem{Claim}{Claim}
\newtheorem{lemma}{Lemma}

\newcommand{\rank}{\operatorname{rank}}
\newcommand{\s}{\scriptscriptstyle}
\newcommand{\p}{\scriptstyle}
\newtheorem{Cor}{Corollary}
\newtheorem{Conj}{Conjecture}

\newcommand{\diag}{\operatorname{diag}}
\newcommand{\doublesubscript}[3]{
\displaystyle\mathop{\displaystyle #1_{#2}}_{#3}}
\title{Grassmannians and Cluster Algebras}

\begin{document}
\maketitle

\section{Introduction}

\bigskip
\bigskip
\noindent
This paper follows the program of study initiated by 
S. Fomin and A. Zelevinsky in \cite{CA1}, and
demonstrates in Theorem \ref{ch2,thm.main} that the homogeneous
coordinate ring of the Grassmannian $\Bbb{G}(k,n)$ is a 
{\it cluster algebra of geometric type}.

\bigskip
\noindent
In its most naive incarnation, a cluster algebra is
a commutative ring generated inside an
ambient field by a family 
of distinguished generators called {\it cluster
variables} which are grouped into 
families or {\it clusters}. These generators are produced,
recursively, inside the ambient field by means
of a process called {\it mutation} from 
an initial fixed family of 
indeterminates. 
The ambient field in the present
case is simply the field of rational functions
$\Bbb{C}\Big( \Bbb{G}(k,n) \Big)$ on
the Grassmannian $\Bbb{G}(k,n)$ and the
initial cluster will be a special family
of Pl\"ucker coordinates. 

\bigskip
\noindent
The description of the coordinate ring
$\Bbb{C}\Big[ \Bbb{G}(k,n) \Big]$ as a cluster
algebra is motivated by the study of the
{\it dual canonical basis}, as put forth
by Lusztig in \cite{Lusztig}, for the coordinate
ring of $SL(n, \Bbb{C} )$ and
its Grassmannians. 
A fundamental characteristic of this basis, 
conjectured first by 
\cite{BZ} and then later proved by \cite{Reineke},
is that it enjoys the following multiplicative
property: if the product $f \cdot g$ of two basis
elements $f$ and $g$ is itself in the basis
then $f$ and $g$ {\it quasi-commute} inside the quantized
coordinate ring $\Bbb{C}_q \Big[ \Bbb{G}(k,n) \Big]$.
With this in mind,
\cite{Leclerc} combinatorially described all pairs of Pl\"ucker
coordinates in the Grassmannian 
$\Bbb{G}(k,n)$ which quasi-commute 
while \cite{Scott} proved that the set of monomials 
consisting of products of 
pairwise quasi-commuting Pl\"ucker coordinates 
is in fact a basis for the coordinate ring
$\Bbb{C}\Big[ \Bbb{G}(2,n) \Big]$.
In an attempt to explicitly compute the dual canonical bases
in general, \cite{CA1} introduced the formalism of 
cluster algebras and anticipates that a large portion
of the dual canonical basis for $\Bbb{C}\Big[ \Bbb{G}(k,n) 
\Big]$ will entail all monomials consisting of products
of cluster variables attached to a common cluster.
  
\bigskip
\noindent
Elements of the dual canonical basis are expected to be {\it positive} 
in the following sense. A point $p$ in the affine cone of the
real Grassmannian $\Bbb{G}(k,n)$
is called {\it positive} if all Pl\"ucker coordinates $\Delta^I(p)$
are positive; a regular function $f$ on the affine cone 
is positive if it takes positive values on positive
points. In section 6 the author proves that cluster variables 
are positive for any Grassmannian. This type of positivity 
is evidence for a {\it positive Laurent} decomposition for
cluster variables conjectured in \cite{CA1}.
 
\bigskip
\noindent
Theorem \ref{ch2,thm.main} of this paper generalizes Proposition 12.6 of
\cite{CA2} which establishes that
the homogeneous coordinate ring of the
Grassmannian $\Bbb{G}(2,n)$ of 2-planes in $n$-space is
a cluster algebra of geometric type. 
Since there is a natural open embedding of
$SL(k,\Bbb{C})$ into the {\it middle} 
Grassmannian $\Bbb{G}(k,2k)$, Theorem \ref{ch2,thm.main}
implies that the coordinate rings of $SL(n,\Bbb{C})$
and its maximal {\it double Bruhat cell}
$SL(n,\Bbb{C})^{w_0,w_0}$ are cluster algebras.
This coincides with the recent work
of \cite{CA3} which proves that the coordinate
ring of any {\it double Bruhat cell} \ $G^{u,v}$ of any 
semi-simple algebraic group is a cluster algebra. 
A feature of the
coordinate ring of any Grassmannian, implicit in the proof
of Theorem 3, is that as a cluster algebra it coincides
with its {\it upper algebra} - as defined by 
\cite{CA3}. This is noteworthy given that
it is not known whether the 
coordinate rings of all double Bruhat cells have
this property.

\bigskip
\noindent
In \cite{Shapiro} the coordinate ring
of a related variety is also studied from the
vantage point of cluster algebras. Specifically,
\cite{Shapiro} considers the open subvariety 
$\Bbb{G}_0(k,n)$ consisting of points in the
Grassmannian $\Bbb{G}(k,n)$ for which the 
Pl\"ucker coordinate $\Delta^{[1 \dots k]}$
is non-zero. A cluster algebra $\mathcal{A}_{\Bbb{G}(k,n)}$
is constructed which is shown to be a subalgebra
of the coordinate ring of $\Bbb{G}_0(k,n)$
and compatible with the Sklyanin Poisson 
structure inherited from $\Bbb{G}(k,n)$. 

\bigskip
\noindent
In proving that $\Bbb{C} \Big[ \Bbb{G}(k,n) \Big]$
is a cluster algebra the author makes use of 
a generalization of double wiring arrangements,
due to A. Postnikov \cite{Postnikov}, called {\it alternating
wiring arrangements}. A key ingredient in the 
proof establishes a bijective correspondence
between a certain class of alternating
wiring arrangements and those clusters in 
$\Bbb{C} \Big[ \Bbb{G}(k,n) \Big]$
consisting only of Pl\"ucker coordinates.
Under this correspondence, the mutation process
between clusters of Pl\"ucker coordinates 
manifests itself geometrically as a certain local 
transformation of the
associated alternating wiring arrangements.

\bigskip
\noindent
One the most striking features of the theory
of cluster algebras is the classification of
cluster algebras of {\it finite type}; i.e. those having
only finitely many cluster variables.
The classification, as put forth in \cite{CA2}, is
identical to the Cartan-Killing classification of
semisimple Lie algebras under which every
Dynkin graph is canonically associated to
a finite type cluster algebra. This pairing
is made more fantastic by a natural
bijective correspondence between cluster variables
and the set of {\it almost positive roots}
- i.e. positive and negative simple roots -
of the associated root system. 
Fomin and Zelevinsky have already shown in \cite{CA2}
that the cluster algebra $\Bbb{C} \Big[ \Bbb{G}(2,n+3) \Big]$
corresponds to the root system of type $A_n$. 

\bigskip
\noindent
In Theorem \ref{ch2,thm.finite} of
this paper the author proves that the only
Grassmannians $\Bbb{G}(k,n)$, within the range
of indices $2<k \leq {n \over 2}$, whose
homogeneous coordinate rings are of finite type
are the Grassmannians
$\Bbb{G}(3,6)$, $\Bbb{G}(3,7)$, and $\Bbb{G}(3,8)$.
As cluster algebras their coordinate rings correspond, respectively,
to the root systems $D_4$, $E_6$, and $E_8$.
The pairing between cluster variables and
almost positive roots is explicitly
worked out in each case and a 
geometric description of all cluster variables 
is presented.

\bigskip
\noindent
This paper is organized accordingly:

\bigskip
\indent Section 2: Cluster Algebra Primer

\indent \indent $\bullet$ Introduction to cluster algebras 
of geometric type and related notions

\indent \indent $\bullet$ Relationship to coordinate rings
of algebraic varieties

\indent \indent $\bullet$ Principle example: the case 
of $\Bbb{C}\Big[ \Bbb{G}(2,n) \Big]$.

\indent \indent $\bullet$ Laurent Phenomena

\bigskip
\indent Section 3: Postnikov Arrangements

\indent \indent $\bullet$ Definition of Postnikov arrangements
attached to a permutation $\pi \in S_n$

\indent \indent $\bullet$ Geometric exchanges and transitivity

\indent \indent $\bullet$ Triangulations and weakly
separated sets

\bigskip
\indent Section 4: Quadrilateral Arrangement

\indent \indent $\bullet$ Construction of quadrilateral
arrangement

\indent \indent $\bullet$ Double wiring arrangements

\bigskip
\indent Section 5: Proof of Main Result

\indent \indent $\bullet$ $\Bbb{C}\Big[ \Bbb{G}(k,n) \Big]$
is a cluster algebra of geometric type

\indent \indent $\bullet$ Clusters associated to 
Postnikov arrangements

\indent \indent $\bullet$ Geometric exchange and
seed mutation

\bigskip
\indent Section 6: Toric Charts and Positivity

\indent \indent $\bullet$ Toric embeddings associated to clusters

\indent \indent $\bullet$ Positivity of cluster variables

\bigskip
\indent Section 7: Grassmannians of finite type

\indent \indent $\bullet$ Classification of Grassmannians of finite type

\indent \indent $\bullet$ Geometric description of cluster variables

\indent \indent $\bullet$ Root correspondence 

\bigskip
\indent Section 8: Future Directions

\indent \indent $\bullet$ Laurent positivity 

\indent \indent $\bullet$ Schur polynomials

\bigskip
\indent Section 9: Acknowledgements and Bibliography

\bigskip
\section{Cluster Algebra Primer}

\bigskip
\noindent
This section presents a treatment which is a slight
alteration of the approach 
documented in \cite{CA2} and
\cite{Zelevinsky}, and considers
class of algebras known as
cluster algebras of {\it geometric type}.

\bigskip
\noindent
Unlike the first treatment of cluster algebras
given by \cite{CA1}, which develops the theory
of {\it exchange patterns} and accommodates a very general
system of coefficients defined over a semifield,
the algebras studied here are defined over
polynomial rings of scalars. 
In addition cluster algebras of geometric type
are defined for the class of {\it skew-symmetrizable} 
matrices whereas \cite{CA1} defines cluster algebras for 
the broader class of {\it sign-skew-symmetric}
matrices. The more recent approach of \cite{CA3} defines 
cluster algebras in terms of auxiliary {\it upper} and {\it lower} 
algebras and is tailored to determine which
cluster algebras are finitely generated.  

\bigskip
\noindent 
Let $\mathcal{F}$ be the field of rational functions
in $M$ independent variables with coefficient 
residing in $\Bbb{C}$. Recall that a 
{\it transcendence basis} for $\mathcal{F}$ is 
an $M$-element set of rational expressions which 
are algebraically independent. A {\it cluster algebra
of geometric type} will be defined as a subalgebra of 
$\mathcal{F}$ generated by a family of transcendence 
bases of $\mathcal{F}$ each of which is determined by
a combinatorial process called {\it mutation}. 

\section*{\it{Matrix Mutation:}}

\noindent
Let ${\bf x} \subset {\bf \tilde{x} }$ be index sets
with respective cardinalities $N = | {\bf x} |
\leq | {\bf \tilde{x} } | = M$. Let $\tilde{B} =
(b_{x,y})$ be a $M \times N$ integer matrix with rows
indexed by ${\bf \tilde{x}}$ and columns indexed
by ${\bf x}$. The $N \times N$ 
square submatrix of $\tilde{B}$ with both rows and columns 
indexed by ${\bf x}$ is called the 
{\it principle submatrix} and is denoted as $B$. 
For $z \in {\bf x}$ define 
a new $M \times N$ integer matrix $\mu_z(\tilde{B})$
with entries $b'_{x,y}$ given by:

\[ b'_{x,y} \ = \ \left\{ \begin{array}{ll}
   -b_{x,y} & \quad \text{if $x$ or $y$ equals $z$} \\
   \ \ b_{x,y} \ + \ { |b_{x,z}| \ b_{z,y} \ + \ 
b_{x,z} \ |b_{z,y}| \over 2} & \quad \text{otherwise} 
\end{array} \right. \]

\bigskip
\noindent
The matrix $\mu_z(\tilde{B})$ is called the {\it
mutation} of $\tilde{B}$ in the {\it $z$ direction}.
One can easily verify the following properties 
of matrix mutation:

\bigskip
\qquad \qquad \qquad $\bullet$ The principle submatrix
of $\mu_z(\Tilde{B})$ is $\mu_z(B)$.

\qquad \qquad \qquad $\bullet$ $\mu \Big( \mu (\tilde{B})
\Big) = \tilde{B}$.

\bigskip
\noindent
If in addition the principle submatrix $B$ is 
{\it skew-symmetrizable} - meaning that $DB$ is
skew-symmetric for some positive integer diagonal square
matrix $D$ - then

\bigskip
\qquad \qquad \qquad $\bullet$ $\mu_z(B)$ is 
skew-symmetrizable and is skew-symmetrized by

\qquad \qquad \qquad the same
diagonal matrix $D$.

\bigskip
\noindent

\section*{\it{Seeds and Clusters:}}

\noindent
Let ${\bf \tilde{x}}$ now denote a transcendence
basis of $\mathcal{F}$. Decompose ${\bf \tilde{x}}$
as a disjoint union ${\bf x} \sqcup {\bf c}$
where $|{\bf x}| = N$ and $|{\bf c}| = M-N$. Let
$\tilde{B}$ be as above and assume that the 
principle submatrix $B$ is skew-symmetrizable.
The triple ${\bf (x,c,\tilde{B}) }$ is called a
{\it seed}.

\bigskip
\noindent
For any rational expression $z \in {\bf x}$ 
define a new seed $\mu_z {\bf (x,c,\tilde{B})
= (x',c,\tilde{B}')}$ where ${\bf x' =
x - \{z\} \cup \{z'\}}$ and $z'$ is given by
the {\it exchange relation}:

\[ zz' \ = \ \doublesubscript{\prod}{x \in {\bf \tilde{x}}}
{b_{x,z}>0} 
x^{b_{x,z}} \ + \ 
\doublesubscript{\prod}{x \in {\bf \tilde{x}}}{b_{x,z}<0}
x^{-b_{x,z}}. \] 

\bigskip
\noindent
The matrix $\tilde{B}'$ is $\mu_z(\tilde{B})$
with the row and column index $z$ replaced
with $z'$. Choose an initial seed
${\bf (x_0,c_0,\tilde{B}_0)}$ and  
let $\mathcal{S}$ be the unique family of 
seeds containing the initial 
seed ${\bf (x_0,c_0,\tilde{B}_0)}$ and 
closed under mutation.
The sets ${\bf x}$ occurring in seeds
${\bf (x,c,\tilde{B})} \in \mathcal{S}$ are called
{\it clusters} and their elements are called
{\it cluster variables}; the elements $c \in
{\bf c}$ are called coefficients.

\bigskip
\begin{Def}[Cluster Algebra of Geometric Type]
The cluster algebra $\mathcal{A} = \mathcal{A}(
\mathcal{S})$ is the 
$\Bbb{C}\big[ c \ | \ c \in {\bf c} \big]$-subalgebra 
of $\mathcal{F}$ generated by all cluster variables
from $\mathcal{S}$. The number $N$ is called
the {\it rank} of the cluster algebra $\mathcal{A}$.

\end{Def}

\bigskip
\noindent
The polynomial coefficient ring used in this definition
is not used in the formulation given by \cite{CA2}
which considers a coefficient system that evolves
along with each seed.
In order to rectify this discrepancy, and to 
implement the following result of \cite{CA2}, 
we require in addition that every coefficient $c \in
{\bf c}$ occurs, by itself, in one of the monomials
of an exchange relation. This additional
assumption will be evident for
the cluster algebras considered in this paper. 

\bigskip
\noindent
Proposition 11.1 of \cite{CA2} provides
sufficient conditions for 
the coordinate ring of an algebraic variety to be a cluster
algebra. It is central to the argument in this paper.

\bigskip
\begin{Prop}\label{ch2,CA2.1}[Fomin, Zelevinsky]
Let $\mathcal{A}$ be a rank $N$ cluster algebra
and let $\mathcal{X}$ denote the set of 
all cluster variables and let
${\bf c}$ be the set of coefficients. Let $X$ be a rational
quasi-affine irreducible algebraic variety
over $\Bbb{C}$. Assume that we are given
a family of regular functions on the variety
$X$ indexed as:

\[ \Big\{ \varphi_x \ \Big| \ x \in \mathcal{X} \ 
\Big\} \cup \Big\{ \ \varphi_c \ \Big| \ c \in
{\bf c} \ \Big\}. \]

\bigskip
\noindent
Assume that the following conditions hold:

\bigskip
\qquad \qquad \qquad $\bullet \ \dim(X) = N + |{\bf c}|$

\qquad \qquad \qquad $\bullet$ The functions $\varphi_x$ 
and $\varphi_c$ generate the coordinate ring $\Bbb{C}[X]$.

\qquad \qquad \qquad $\bullet$ Every exchange relation
in the cluster algebra $\mathcal{A}$ is valid in 
$\Bbb{C}[X]$.

\bigskip
\noindent
Then the correspondences 
$x \mapsto \varphi_x$ and $c \mapsto \varphi_c$ 
extend uniquely to an algebra isomorphism of
the cluster algebra $\mathcal{A}$ and $\Bbb{C}[X]$.
\end{Prop}

\bigskip
\noindent
\section*{{\it Example}}

\noindent
The motivating example for this paper of an algebraic
variety which gives rise to a cluster algebra is 
the Grassmannian $\Bbb{G}(2,n)$. 

\bigskip
\noindent
Recall that the {\it affine cone $X(k,n)$} of any Grassmannian
$\Bbb{G}(k,n)$ is the affine 
subvariety of decomposable $k$-forms in $\bigwedge^k
\Bbb{C}^n$. The coordinates of $\bigwedge^k \Bbb{C}^n$ 
- otherwise known as {\it Pl\"ucker coordinates}- are
indexed by $k$-subsets $A$ of $[1 \dots n]$ and are 
denoted as $\Delta^A$. The Pl\"ucker coordinates restrict
to regular function on the affine cone and 
generate its coordinate ring. As regular functions on
$X(k,n)$ these coordinates obey certain determinental 
constraints known as {\it Pl\"ucker relations},
the simplest (and most relevant) of which has the form

\begin{equation}\label{ch2,short} \Delta^{Iij} \Delta^{Ist} \ = \ 
\Delta^{Iis}\Delta^{Ijt} \ + \ \Delta^{Ijs}\Delta^{Iit} 
\end{equation}

\bigskip
\noindent
where $I$ is subset of size $k-2$ disjoint from
the four indices $\{i,j,s,t\}$ and the
pairs $\{i,j\}$ and $\{s,t\}$ are {\it crossing} -
see definition below. 
The Grassmannian $\Bbb{G}(k,n)$ is obtained
as the projective quotient of the affine cone $X(k,n)$.  
The homogeneous coordinate ring $\Bbb{C}\Big[ \Bbb{G}(k,n) \Big]$
is defined as the coordinate ring of the affine cone
$X(k,n)$. This paper aims to show that this
ring is a cluster algebra.

\bigskip
\noindent
In the case of $\Bbb{G}(2,n)$ the ambient field $\mathcal{F}$ 
is identified with the field of rational functions on the Grassmannian 
$ \Bbb{C} \Big( \Bbb{G}(2,n) \Big)$. Fomin and Zelevinsky proved that
its homogeneous coordinate ring is a cluster algebra
of finite type; i.e. having only finitely many
clusters. These clusters in turn are parameterized
by triangulations on a regular $n$-gon. 

\bigskip
\noindent
Let ${\bf P}$ be a regular $n$-gon with vertices labeled 
clockwise and in order by the indices $1 \dots n$. 
By a chord $[ij]$ I will mean a line segment joining vertices
$i$ and $j$. Two chords (or two sets) $[ij]$ and $[st]$
are {\it crossing} if they intersect within the interior
of the polygon ${\bf P}$.
A triangulation of ${\bf P}$ will mean 
a maximal family of non-crossing chords $[ij]$ in ${\bf P}$.
For instance, the following diagram depicts
a triangulation of an octagon. 
  
\bigskip
\bigskip
\begin{figure}[h]
\begin{center}
\includegraphics[width=2.0in]{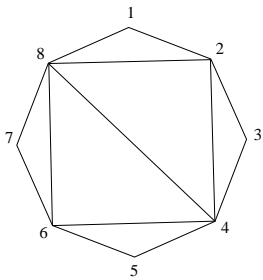}
\end{center} 
\caption{Triangulation of an Octagon}
\label{8-gon}
\end{figure}

\bigskip
\noindent
To the chord $[ij]$ in ${\bf P}$ associate the Pl\"ucker
coordinate $\Delta^{ij}$. To each triangulation
$T$ associate the collection ${\bf \tilde{x}(T)}$ of Pl\"ucker 
coordinates arising from all chords in $T$. Notice that
the external or {\it boundary} edges $[12],[2,3],\dots, [1n]$ 
are present in every triangulation, and so the matching 
Pl\"ucker coordinates $\Delta^{12}, \Delta^{23}, \dots, \Delta^{1n}$
must belong to each ${\bf \tilde{x}(T)}$; let ${\bf c}$ denote
the set of boundary Pl\"ucker coordinates and set
${\bf x(T) = \tilde{x} - c}$. We can now state
Proposition 12.5 of \cite{CA2}.

\bigskip
\begin{Prop}\label{ch2,CA2.2}[Fomin, Zelevinsky] 
For $n \geq 5$ the homogeneous coordinate ring $\Bbb{C}\Big[ 
\Bbb{G}(2,n) \Big]$ is a cluster algebra whose seeds are 
${ \bf \Big(x(T),c,B(T) \Big) }$; here $B(T)$
is the matrix with rows indexed by ${\bf \tilde{x}(T)}$ and 
columns indexed by ${\bf x(T)}$ with all entries
$0$ except $b_{[ij],[ik]}$ which equals $1$ (resp. $-1$) provided the
ordered vertices $i,j,k$ are oriented counter-clockwise
(resp. clockwise).

\end{Prop}
 
\bigskip
\noindent
In this context mutation can be geometrically interpreted
on the level of triangulations as follows: two triangulations
are mutations of one another if one can be obtained from
the other by a {\it flip} that replaces a diagonal chord
in a quadrilateral formed by two triangles of one
triangulation by the crossing diagonal in the same quadrilateral.
The picture below demonstrates this exchange for the
quadrilateral with corner vertices $2$, $4$, $6$, and
$8$ within a triangulation of
an octagon.

\bigskip
\bigskip
\begin{figure}[h]
\begin{center}
\includegraphics[width=3.0in]{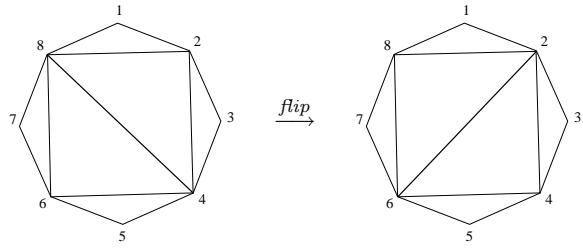}
\put(-117,43){{\bf $\stackrel{\text{{\it flip}}}{\longrightarrow} $}}
\caption{Quadrilateral Flip}
\label{flip}
\end{center} 
\end{figure} 

\bigskip
\noindent
The exchange relation represented by such a 
quadrilateral flip is a short
Pl\"ucker relation of the kind described in formula
\ref{ch2,short}. For instance, in the example above, the
flip exchanging chords $[48]$ and $[26]$ abbreviates
the relation 

\[ \Delta^{48} \Delta^{26} \ = \ 
\Delta^{28} \Delta^{46} \ + \ \Delta^{24} \Delta^{68} \]

\noindent
inside $\Bbb{C}\Big[ \Bbb{G}(2,8) \Big]$.
All mutations correspond to flips and consequently
every exchange relation is given by some short
Pl\"ucker relation.

\bigskip
\section*{ {\it Laurent Phenomena}}

\noindent
Cluster algebras were designed in part to reflect
integrality properties of certain
bilinear recurrences such as the famed
Somos 4-recurrence  

\[ y_k =
{ y_{k-3}y_{k-1} + y_{k-2}^2 \over {y_{k-4}}}\]

\bigskip
\noindent
invented by Michael Somos and subsequently
studied by Fomin and Zelevinsky.
If the initial values are all set to $1$,
recurrences of this type always
produce integer sequences despite
the fact that the recurrence relations
are rational. 
The usefulness of cluster algebras to this 
end of study stems from
a property known as the {\it Laurent
phenomena} which was first
documented and proved in \cite{CA1}
and later re-examined in \cite{Laurent} and
\cite{CA2}.
The following version of this property - Proposition
11.2 of \cite{CA2} - is used
to implement Proposition 1 in
the proof of Theorem 3 of this paper.

\bigskip
\begin{Prop}\label{ch2,laurent}[Fomin-Zelevinsky]
Let $\mathcal{A}$ be a cluster algebra
of geometric type and let ${\bf x}$ be
a cluster. Any cluster variable $x
\in \mathcal{X}$ can be uniquely
expressed as a Laurent polynomial 
in the indeterminates ${\bf x \sqcup
c}$ whose monomial denominator 
consists only of indeterminates from
${\bf x}$. 
\end{Prop}

\noindent
In the case of the Grassmannian $\Bbb{G}(2,n)$
Proposition 3 can be sharpened requiring the predicted 
Laurent expressions to have positive coefficients
in the numerator. 

\bigskip
\begin{Prop}\label{ch2,CA2.3}[Fomin-Zelevinsky]
Let ${\bf T}$ be a triangulation and let
${\bf x(T)}$ and ${\bf c}$
be as above. Let $\Delta^{ij}$ be any Pl\"ucker
coordinate for $\Bbb{G}(2,n)$. Then $\Delta^{ij}$
can be uniquely expressed as 

\[ { p \Big( {\bf \tilde{x}(T) } \Big)
\over { m \Big( {\bf x(T)} \Big) }} \]

\bigskip
\noindent
where $m \Big( {\bf x(T) } \Big)$ is a monomial of Pl\"ucker
coordinates in ${\bf x(T)}$ and $p \Big( {\bf \tilde{x}(T)}
\Big) $ is a polynomial with {\bf positive}
coefficients in the Pl\"ucker coordinates from
${\bf \tilde{x}(T) = x(T) \sqcup c}$.
\end{Prop}

\bigskip
\section{Postnikov Arrangements}

\bigskip
\noindent
This section presents a construction developed
by A. Postnikov in \cite{Postnikov} and used here to 
parameterize clusters
of Pl\"ucker coordinates inside the homogeneous
coordinate ring of the Grassmannian $\Bbb{G}(k,n)$.

\bigskip
\noindent
\begin{Def}
Let $\pi \in S_n$. Label the vertices of a convex $2n$-gon
clock-wise by the indices $1',1,2',2,\dots,n',n$. 
A Postnikov arrangement for $\pi$ is
a collection of $n$ oriented paths in the
interior of the polygon; the $i$-th path joins
the vertex $i$ with the vertex $\pi(i)'$ and is directed towards
$\pi(i)'$. The collection of paths must satisfy the
following compatibility relations:

\bigskip
\indent \qquad \qquad 1. No path intersects itself.

\bigskip
\indent \qquad \qquad 2. All path intersections are transversal.

\bigskip
\indent \qquad \qquad 3. As a path is
traversed from its source to its sink, the
paths 

\indent \qquad \qquad \ \ \ intersecting it must alternate 
in orientation cutting
the wire in 

\indent \qquad \qquad \ \ \ 
question first right, then left, right, $\dots$, left,  
and finally right.

\bigskip
\indent \qquad \qquad 4. For any two paths $i$ and $j$ the
following configuration is forbidden:

\end{Def}

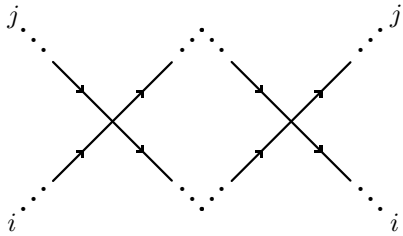
\begin{figure}[ht]

\setlength{\unitlength}{.75pt}

\begin{center}
\begin{picture}(270,50)(-50,100)
\thicklines

    \put(20,130){\line(1,-1){60}}
    \put(-3,150){$j$}
    \put(3,145){{\bf .}}
    \put(8,140){{\bf .}}
    \put(13,135){{\bf .}}
    \put(35,115){\line(0,1){3}}
    \put(35,115){\line(-1,0){3}}
    \put(65,115){\line(0,-1){3}}
    \put(65,115){\line(-1,0){3}}
    \put(65,85){\line(0,1){3}}
    \put(65,85){\line(-1,0){3}}
    \put(83,65){{\bf .}}
    \put(88,60){{\bf .}}
    \put(93,55){{\bf .}}

    \put(20,70){\line(1,1){60}}  
    \put(-3,45){$i$}
    \put(3,55){{\bf .}}
    \put(8,60){{\bf .}}
    \put(13,65){{\bf .}}
    \put(35,85){\line(-1,0){3}}
    \put(35,85){\line(0,-1){3}} 
    \put(83,135){{\bf .}}
    \put(88,140){{\bf .}}
    \put(93,145){{\bf .}}

    \put(110,130){\line(1,-1){60}}
    \put(93,145){{\bf .}}
    \put(98,140){{\bf .}}
    \put(103,135){{\bf .}}
    \put(125,115){\line(0,1){3}}
    \put(125,115){\line(-1,0){3}}
    \put(155,115){\line(0,-1){3}}
    \put(155,115){\line(-1,0){3}}
    \put(155,85){\line(0,1){3}}
    \put(155,85){\line(-1,0){3}}
    \put(173,65){{\bf .}}
    \put(178,60){{\bf .}}
    \put(183,55){{\bf .}}

    \put(110,70){\line(1,1){60}}  
    \put(93,55){{\bf .}}
    \put(98,60){{\bf .}}
    \put(103,65){{\bf .}}
    \put(125,85){\line(-1,0){3}}
    \put(125,85){\line(0,-1){3}} 
    \put(173,135){{\bf .}}
    \put(178,140){{\bf .}}
    \put(183,145){{\bf .}}
    \put(190,45){$i$}
    \put(190,150){$j$}

\end{picture}
\end{center}
\bigskip
\bigskip
\bigskip
\bigskip
\bigskip
\caption{Forbidden Crossing}
\label{crossing}
\end{figure}

\bigskip
\bigskip
\bigskip
\noindent
In the text, Postnikov arrangements attached to a permutation
$\pi$ are referred to as {\it $\pi$-diagrams}.
The following diagram is a Postnikov
arrangement for the permutation

\[ \pi =
\begin{pmatrix} 1 & 2 & 3 & 4 & 5 & 6 & 7 \\
4 & 5 & 6 & 7 & 1 & 2 & 3 \end{pmatrix}. \]

\newpage
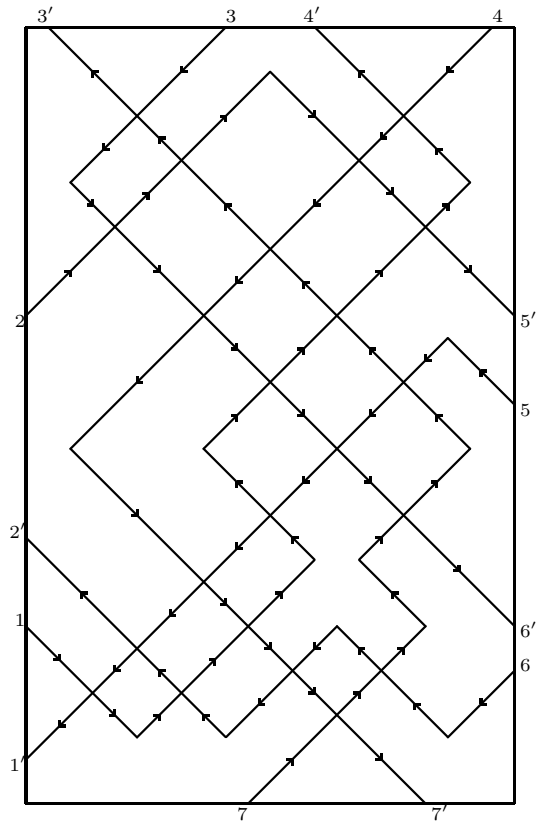
\begin{figure}[ht]

\setlength{\unitlength}{2.10pt}

\begin{center}
\begin{picture}(270,70)(-45,78)
\thicklines

    \put(0,8){\line(1,1){76}}
    \put(76,84){\line(1,-1){12}}
    \put(6,14){\line(0,1){1}}
    \put(6,14){\line(1,0){1}}    
    \put(16,24){\line(0,1){1}}
    \put(16,24){\line(1,0){1}}    
    \put(26,34){\line(0,1){1}}
    \put(26,34){\line(1,0){1}}    
    \put(38,46){\line(0,1){1}}
    \put(38,46){\line(1,0){1}}    
    \put(50,58){\line(0,1){1}}
    \put(50,58){\line(1,0){1}}    
    \put(62,70){\line(0,1){1}}
    \put(62,70){\line(1,0){1}}    
    \put(72,80){\line(0,1){1}}
    \put(72,80){\line(1,0){1}}    
    \put(82,78){\line(0,-1){1}}
    \put(82,78){\line(1,0){1}}    
   
    \put(0,88){\line(1,1){44}}
    \put(44,132){\line(1,-1){44}}
    \put(8,96){\line(-1,0){1}}
    \put(8,96){\line(0,-1){1}}
    \put(22,110){\line(-1,0){1}}
    \put(22,110){\line(0,-1){1}}
    \put(36,124){\line(-1,0){1}}
    \put(36,124){\line(0,-1){1}}
    \put(52,124){\line(0,1){1}}
    \put(52,124){\line(-1,0){1}}
    \put(66,110){\line(0,1){1}}
    \put(66,110){\line(-1,0){1}}
    \put(80,96){\line(0,1){1}}
    \put(80,96){\line(-1,0){1}}

    \put(8,112){\line(1,1){28}}
    \put(8,112){\line(1,-1){80}}
    \put(28,132){\line(1,0){1}}
    \put(28,132){\line(0,1){1}}
    \put(14,118){\line(1,0){1}}
    \put(14,118){\line(0,1){1}}
    \put(12,108){\line(0,1){1}}
    \put(12,108){\line(-1,0){1}}
    \put(24,96){\line(0,1){1}}
    \put(24,96){\line(-1,0){1}}
    \put(38,82){\line(0,1){1}}
    \put(38,82){\line(-1,0){1}}
    \put(50,70){\line(0,1){1}}
    \put(50,70){\line(-1,0){1}}
    \put(62,58){\line(0,1){1}}
    \put(62,58){\line(-1,0){1}}
    \put(78,42){\line(0,1){1}}
    \put(78,42){\line(-1,0){1}}

    \put(8,64){\line(1,1){76}}
    \put(8,64){\line(1,-1){64}}
    \put(76,132){\line(1,0){1}}
    \put(76,132){\line(0,1){1}}
    \put(64,120){\line(1,0){1}}
    \put(64,120){\line(0,1){1}}
    \put(52,108){\line(1,0){1}}
    \put(52,108){\line(0,1){1}}
    \put(38,94){\line(1,0){1}}
    \put(38,94){\line(0,1){1}}
    \put(20,76){\line(1,0){1}}
    \put(20,76){\line(0,1){1}}
    \put(20,52){\line(0,1){1}}
    \put(20,52){\line(-1,0){1}}
    \put(36,36){\line(0,1){1}}
    \put(36,36){\line(-1,0){1}}
    \put(44,28){\line(0,1){1}}
    \put(44,28){\line(-1,0){1}}
    \put(52,20){\line(0,1){1}}
    \put(52,20){\line(-1,0){1}}
    \put(64,8){\line(0,1){1}}
    \put(64,8){\line(-1,0){1}}

    \put(20,12){\line(-1,1){20}}
    \put(20,12){\line(1,1){32}}
    \put(32,64){\line(1,1){48}}
    \put(32,64){\line(1,-1){20}}
    \put(80,112){\line(-1,1){28}}
    \put(6,26){\line(0,1){1}}
    \put(6,26){\line(-1,0){1}} 
    \put(16,16){\line(0,1){1}}
    \put(16,16){\line(-1,0){1}}
    \put(24,16){\line(0,-1){1}}
    \put(24,16){\line(-1,0){1}}
    \put(34,26){\line(0,-1){1}}
    \put(34,26){\line(-1,0){1}}
    \put(46,38){\line(0,-1){1}}
    \put(46,38){\line(-1,0){1}}
    \put(48,48){\line(1,0){1}}
    \put(48,48){\line(0,-1){1}}
    \put(38,58){\line(1,0){1}}
    \put(38,58){\line(0,-1){1}}
    \put(38,70){\line(0,-1){1}}
    \put(38,70){\line(-1,0){1}}
    \put(50,82){\line(0,-1){1}}
    \put(50,82){\line(-1,0){1}}
    \put(64,96){\line(0,-1){1}}
    \put(64,96){\line(-1,0){1}}
    \put(76,108){\line(0,-1){1}}
    \put(76,108){\line(-1,0){1}}
    \put(74,118){\line(1,0){1}}
    \put(74,118){\line(0,-1){1}}
    \put(60,132){\line(1,0){1}}
    \put(60,132){\line(0,-1){1}}

    \put(0,48){\line(1,-1){36}}
    \put(56,32){\line(-1,-1){20}}
    \put(56,32){\line(1,-1){20}}
    \put(76,12){\line(1,1){12}}
    \put(10,38){\line(1,0){1}}
    \put(10,38){\line(0,-1){1}}
    \put(24,24){\line(1,0){1}}
    \put(24,24){\line(0,-1){1}}
    \put(32,16){\line(1,0){1}}
    \put(32,16){\line(0,-1){1}}
    \put(42,18){\line(0,1){1}}
    \put(42,18){\line(1,0){1}}
    \put(52,28){\line(0,1){1}}
    \put(52,28){\line(1,0){1}}
    \put(60,28){\line(0,-1){1}}
    \put(60,28){\line(1,0){1}}
    \put(70,18){\line(0,-1){1}}
    \put(70,18){\line(1,0){1}}
    \put(82,18){\line(0,1){1}}
    \put(82,18){\line(1,0){1}}

    \put(72,32){\line(-1,-1){32}}
    \put(72,32){\line(-1,1){12}}
    \put(80,64){\line(-1,-1){20}}
    \put(80,64){\line(-1,1){76}}
    \put(48,8){\line(-1,0){1}}
    \put(48,8){\line(0,-1){1}}
    \put(60,20){\line(-1,0){1}}
    \put(60,20){\line(0,-1){1}}
    \put(68,28){\line(-1,0){1}}
    \put(68,28){\line(0,-1){1}}
    \put(66,38){\line(0,-1){1}}
    \put(66,38){\line(1,0){1}}
    \put(64,48){\line(-1,0){1}}
    \put(64,48){\line(0,-1){1}}
    \put(74,58){\line(-1,0){1}}
    \put(74,58){\line(0,-1){1}}
    \put(74,70){\line(0,-1){1}}
    \put(74,70){\line(1,0){1}}
    \put(62,82){\line(0,-1){1}}
    \put(62,82){\line(1,0){1}}
    \put(50,94){\line(0,-1){1}}
    \put(50,94){\line(1,0){1}}
    \put(36,108){\line(0,-1){1}}
    \put(36,108){\line(1,0){1}} 
    \put(24,120){\line(0,-1){1}}
    \put(24,120){\line(1,0){1}}
    \put(12,132){\line(0,-1){1}}
    \put(12,132){\line(1,0){1}}

    \put(0,0){\line(1,0){88}}
    \put(0,0){\line(0,1){140}}
    \put(88,140){\line(-1,0){88}}
    \put(88,140){\line(0,-1){140}}
   
    \put(-3,6){{\bf $\p 1'$}} 
    \put(-2,32){{\bf $\p 1$}}
    \put(-3,48){{\bf $\p 2'$}}  
    \put(-2,86){{\bf $\p 2$}}
    \put(89,86){{\bf $\p 5'$}}
    \put(89,30){{\bf $\p 6'$}}
    \put(89,24){{\bf $\p 6$}}
    \put(89,70){{\bf $\p 5$}}
    \put(2,141){{\bf $\p 3'$}}
    \put(36,141){{\bf $\p 3$}}
    \put(50,141){{\bf $\p 4'$}}
    \put(84,141){{\bf $\p 4$}}    
    \put(38,-3){{\bf $\p 7$}}
    \put(73,-3){{\bf $\p 7'$}}

\end{picture}

\bigskip
\bigskip
\bigskip
\bigskip
\bigskip
\bigskip
\bigskip
\bigskip
\bigskip
\bigskip
\bigskip
\bigskip
\bigskip
\bigskip
\bigskip
\bigskip
\bigskip
\bigskip
\bigskip
\bigskip
\caption{Postnikov Arrangement}
\label{arrangement}
\end{center}
\end{figure}

\newpage
\bigskip
\noindent
Postnikov diagrams are identified up to {\it isotopy} 
- i.e. distortions of the configuration
which neither introduce nor remove crossings -
and up to the following local motion 
{\it untwisting} consecutive 
crossings of two paths:

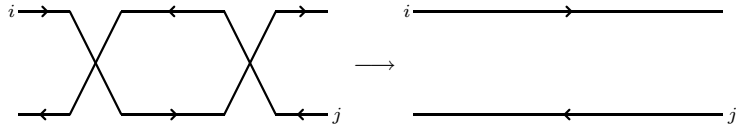
\begin{figure}[ht]

\setlength{\unitlength}{0.65pt}

\begin{center}
\begin{picture}(320,50)(50,50)

\thicklines
\put(0,0){\line(1,0){30}}
\put(60,0){\line(1,0){60}}
\put(150,0){\line(1,0){30}}

\put(0,60){\line(1,0){30}}
\put(60,60){\line(1,0){60}}
\put(150,60){\line(1,0){30}}

\put(30,60){\line(1,-2){30}}
\put(30,0){\line(1,2){30}}
\put(120,60){\line(1,-2){30}}
\put(120,0){\line(1,2){30}}

\put(13,0){\line(1,1){3}}
\put(13,0){\line(1,-1){3}}
\put(88,60){\line(1,1){3}}
\put(88,60){\line(1,-1){3}}
\put(163,0){\line(1,1){3}}
\put(163,0){\line(1,-1){3}}

\put(17,60){\line(-1,1){3}}
\put(17,60){\line(-1,-1){3}}
\put(92,0){\line(-1,1){3}}
\put(92,0){\line(-1,-1){3}}
\put(167,60){\line(-1,1){3}}
\put(167,60){\line(-1,-1){3}}

\put(-6,57){$\p{i}$}
\put(183,-3){$\p{j}$}
\put(195,25){{\bf $\longrightarrow$}}

\put(230,0){\line(1,0){180}}
\put(230,60){\line(1,0){180}}
\put(322,60){\line(-1,1){3}}
\put(322,60){\line(-1,-1){3}}
\put(318,0){\line(1,1){3}}
\put(318,0){\line(1,-1){3}}
\put(224,57){$\p{i}$}
\put(413,-3){$\p{j}$}

\end{picture}

\bigskip
\bigskip
\bigskip
\bigskip
\caption{Untwisting Move}
\label{untwisting}
\end{center}
\end{figure}
  
\bigskip
\noindent
An immediate consequence of condition 3 is that
the interior of a Postnikov arrangement is subdivided
into two types of regions: those regions having a boundary
consisting of distinct wires which, consecutively, alternate
in orientation and those whose boundary
consists of distinct oriented wires which, when
taken consecutively, patch together to
orient the entire boundary (either clockwise or
counter-clockwise) with the added
proviso that the boundary of the $2n$-gon is
given the clockwise orientation.

\bigskip
\noindent 
We will call regions of the first type {\it even}
and regions of the second type {\it odd}.
Label an even region
with the index $i$ if the $i$-th wire stays to the right
of the region; i.e. if the circuit obtained by first
traversing the $i$-th wire and then traversing the
boundary of the polygon clockwise from $\pi(i)'$ to
$i$ does not wind around the region. For example, when the
arrangement above is labeled we obtain:

\newpage
\begin{figure}[ht]

\setlength{\unitlength}{2.10pt}

\begin{center}
\begin{picture}(270,60)(-45,90)
\thicklines

    \put(0,8){\line(1,1){76}}
    \put(76,84){\line(1,-1){12}}
    \put(6,14){\line(0,1){1}}
    \put(6,14){\line(1,0){1}}    
    \put(16,24){\line(0,1){1}}
    \put(16,24){\line(1,0){1}}    
    \put(26,34){\line(0,1){1}}
    \put(26,34){\line(1,0){1}}    
    \put(38,46){\line(0,1){1}}
    \put(38,46){\line(1,0){1}}    
    \put(50,58){\line(0,1){1}}
    \put(50,58){\line(1,0){1}}    
    \put(62,70){\line(0,1){1}}
    \put(62,70){\line(1,0){1}}    
    \put(72,80){\line(0,1){1}}
    \put(72,80){\line(1,0){1}}    
    \put(82,78){\line(0,-1){1}}
    \put(82,78){\line(1,0){1}}    
   
    \put(0,88){\line(1,1){44}}
    \put(44,132){\line(1,-1){44}}
    \put(8,96){\line(-1,0){1}}
    \put(8,96){\line(0,-1){1}}
    \put(22,110){\line(-1,0){1}}
    \put(22,110){\line(0,-1){1}}
    \put(36,124){\line(-1,0){1}}
    \put(36,124){\line(0,-1){1}}
    \put(52,124){\line(0,1){1}}
    \put(52,124){\line(-1,0){1}}
    \put(66,110){\line(0,1){1}}
    \put(66,110){\line(-1,0){1}}
    \put(80,96){\line(0,1){1}}
    \put(80,96){\line(-1,0){1}}

    \put(8,112){\line(1,1){28}}
    \put(8,112){\line(1,-1){80}}
    \put(28,132){\line(1,0){1}}
    \put(28,132){\line(0,1){1}}
    \put(14,118){\line(1,0){1}}
    \put(14,118){\line(0,1){1}}
    \put(12,108){\line(0,1){1}}
    \put(12,108){\line(-1,0){1}}
    \put(24,96){\line(0,1){1}}
    \put(24,96){\line(-1,0){1}}
    \put(38,82){\line(0,1){1}}
    \put(38,82){\line(-1,0){1}}
    \put(50,70){\line(0,1){1}}
    \put(50,70){\line(-1,0){1}}
    \put(62,58){\line(0,1){1}}
    \put(62,58){\line(-1,0){1}}
    \put(78,42){\line(0,1){1}}
    \put(78,42){\line(-1,0){1}}

    \put(8,64){\line(1,1){76}}
    \put(8,64){\line(1,-1){64}}
    \put(76,132){\line(1,0){1}}
    \put(76,132){\line(0,1){1}}
    \put(64,120){\line(1,0){1}}
    \put(64,120){\line(0,1){1}}
    \put(52,108){\line(1,0){1}}
    \put(52,108){\line(0,1){1}}
    \put(38,94){\line(1,0){1}}
    \put(38,94){\line(0,1){1}}
    \put(20,76){\line(1,0){1}}
    \put(20,76){\line(0,1){1}}
    \put(20,52){\line(0,1){1}}
    \put(20,52){\line(-1,0){1}}
    \put(36,36){\line(0,1){1}}
    \put(36,36){\line(-1,0){1}}
    \put(44,28){\line(0,1){1}}
    \put(44,28){\line(-1,0){1}}
    \put(52,20){\line(0,1){1}}
    \put(52,20){\line(-1,0){1}}
    \put(64,8){\line(0,1){1}}
    \put(64,8){\line(-1,0){1}}

    \put(20,12){\line(-1,1){20}}
    \put(20,12){\line(1,1){32}}
    \put(32,64){\line(1,1){48}}
    \put(32,64){\line(1,-1){20}}
    \put(80,112){\line(-1,1){28}}
    \put(6,26){\line(0,1){1}}
    \put(6,26){\line(-1,0){1}} 
    \put(16,16){\line(0,1){1}}
    \put(16,16){\line(-1,0){1}}
    \put(24,16){\line(0,-1){1}}
    \put(24,16){\line(-1,0){1}}
    \put(34,26){\line(0,-1){1}}
    \put(34,26){\line(-1,0){1}}
    \put(46,38){\line(0,-1){1}}
    \put(46,38){\line(-1,0){1}}
    \put(48,48){\line(1,0){1}}
    \put(48,48){\line(0,-1){1}}
    \put(38,58){\line(1,0){1}}
    \put(38,58){\line(0,-1){1}}
    \put(38,70){\line(0,-1){1}}
    \put(38,70){\line(-1,0){1}}
    \put(50,82){\line(0,-1){1}}
    \put(50,82){\line(-1,0){1}}
    \put(64,96){\line(0,-1){1}}
    \put(64,96){\line(-1,0){1}}
    \put(76,108){\line(0,-1){1}}
    \put(76,108){\line(-1,0){1}}
    \put(74,118){\line(1,0){1}}
    \put(74,118){\line(0,-1){1}}
    \put(60,132){\line(1,0){1}}
    \put(60,132){\line(0,-1){1}}

    \put(0,48){\line(1,-1){36}}
    \put(56,32){\line(-1,-1){20}}
    \put(56,32){\line(1,-1){20}}
    \put(76,12){\line(1,1){12}}
    \put(10,38){\line(1,0){1}}
    \put(10,38){\line(0,-1){1}}
    \put(24,24){\line(1,0){1}}
    \put(24,24){\line(0,-1){1}}
    \put(32,16){\line(1,0){1}}
    \put(32,16){\line(0,-1){1}}
    \put(42,18){\line(0,1){1}}
    \put(42,18){\line(1,0){1}}
    \put(52,28){\line(0,1){1}}
    \put(52,28){\line(1,0){1}}
    \put(60,28){\line(0,-1){1}}
    \put(60,28){\line(1,0){1}}
    \put(70,18){\line(0,-1){1}}
    \put(70,18){\line(1,0){1}}
    \put(82,18){\line(0,1){1}}
    \put(82,18){\line(1,0){1}}

    \put(72,32){\line(-1,-1){32}}
    \put(72,32){\line(-1,1){12}}
    \put(80,64){\line(-1,-1){20}}
    \put(80,64){\line(-1,1){76}}
    \put(48,8){\line(-1,0){1}}
    \put(48,8){\line(0,-1){1}}
    \put(60,20){\line(-1,0){1}}
    \put(60,20){\line(0,-1){1}}
    \put(68,28){\line(-1,0){1}}
    \put(68,28){\line(0,-1){1}}
    \put(66,38){\line(0,-1){1}}
    \put(66,38){\line(1,0){1}}
    \put(64,48){\line(-1,0){1}}
    \put(64,48){\line(0,-1){1}}
    \put(74,58){\line(-1,0){1}}
    \put(74,58){\line(0,-1){1}}
    \put(74,70){\line(0,-1){1}}
    \put(74,70){\line(1,0){1}}
    \put(62,82){\line(0,-1){1}}
    \put(62,82){\line(1,0){1}}
    \put(50,94){\line(0,-1){1}}
    \put(50,94){\line(1,0){1}}
    \put(36,108){\line(0,-1){1}}
    \put(36,108){\line(1,0){1}} 
    \put(24,120){\line(0,-1){1}}
    \put(24,120){\line(1,0){1}}
    \put(12,132){\line(0,-1){1}}
    \put(12,132){\line(1,0){1}}

    \put(0,0){\line(1,0){88}}
    \put(0,0){\line(0,1){140}}
    \put(88,140){\line(-1,0){88}}
    \put(88,140){\line(0,-1){140}}
   
    \put(-3,6){{\bf $\p 1'$}} 
    \put(-2,32){{\bf $\p 1$}}
    \put(-3,48){{\bf $\p 2'$}}  
    \put(-2,86){{\bf $\p 2$}}
    \put(89,86){{\bf $\p 5'$}}
    \put(89,30){{\bf $\p 6'$}}
    \put(89,24){{\bf $\p 6$}}
    \put(89,70){{\bf $\p 5$}}
    \put(2,141){{\bf $\p 3'$}}
    \put(36,141){{\bf $\p 3$}}
    \put(50,141){{\bf $\p 4'$}}
    \put(84,141){{\bf $\p 4$}}    
    \put(38,-3){{\bf $\p 7$}}
    \put(73,-3){{\bf $\p 7'$}}

    \put(25,6){{\bf $\p 567 $}}
    \put(25,28){{\bf $\p 157$}}
    \put(9,28){{\bf $\p 167$}}
    \put(19,63){{\bf $\p 147$}}
    \put(77,6){{\bf $\p 456$}}
    \put(53,38){{\bf $\p 457$}}
    \put(53,75){{\bf $\p 347$}}
    \put(53,98){{\bf $\p 134$}}
    \put(77,51){{\bf $\p 345$}}
    \put(77,123){{\bf $\p 234$}}
    \put(29,98){{\bf $\p 137$}}
    \put(9,123){{\bf $\p 127$}}
    \put(41,135){{\bf $\p 123$}}

\end{picture}

\bigskip
\bigskip
\bigskip
\bigskip
\bigskip
\bigskip
\bigskip
\bigskip
\bigskip
\bigskip
\bigskip
\bigskip
\bigskip
\bigskip
\bigskip
\bigskip
\bigskip
\bigskip
\bigskip
\bigskip
\bigskip
\bigskip
\bigskip
\caption{Labeled Postnikov Arrangement}
\label{labeled arrangement}
\end{center}
\end{figure}
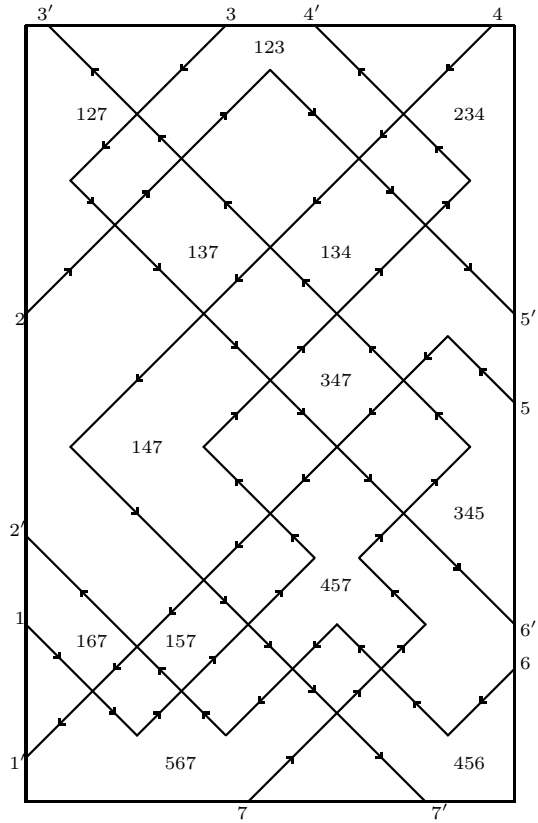

\newpage
\bigskip
\noindent
Even cells on the boundary of a Postnikov arrangement
are called {\it boundary} cells;
the interior even cells are called {\it interior}
or {\it internal}. In the case of the {\it Grassmann}
permutation, given by

\[ \pi_{k,n} \ = \ \begin{pmatrix} 
1 & \cdots & n-k & n-k+1 & \cdots & n \\
k+1 & \cdots & n & 1 & \cdots & k  \end{pmatrix}, \]

\bigskip
\noindent
the $k$-subsets labeling the boundary cells of
{\bf any} $\pi_{k,n}$-diagram are always
the intervals or {\it boundary $k$-subsets}  

\[ \big[ 1 \dots k \big] \ \ \big[ 2 \dots k+1 \big] \ \ 
\big[ 3 \dots k+2 \big] \ \ \dots \ \ 
\big[ n \dots k-1 \big] . \]

\bigskip
\noindent
In particular, any $\pi_{k,n}$-diagram contains precisely
$n$ boundary cells. Further properties of $\pi_{k,n}$-diagrams
are summarized by the next proposition of A. Postnikov.

\bigskip
\begin{Prop}\label{ch2,Post.1}[Postnikov] Let \ $\pi_{k,n}$ be the Grassmann
permutation.
 
\bigskip
\bigskip
\indent \indent 1. The number of even regions
in a Postnikov arrangement for $\pi_{k,n}$ is 

\indent \indent $k(n-k)+1$.

\bigskip
\indent \indent 2. Each even region is labeled by
exactly $k$ distinct indices from $[1 \dots n]$. 

\bigskip
\indent \indent 3. Every $k$-subset in $[1 \dots n]$
occurs as the labeling set of an even cell in
some

\indent \indent  $\pi_{k,n}$-diagram.                           
\end{Prop}

\bigskip
\noindent
Given a $\pi$-diagram ${\bf A}$ and an even 
{\bf quadrilateral} cell
inside ${\bf A}$ a new $\pi$-diagram ${\bf A'}$ is
constructed by the following local rearrangement -
called a {\it geometric exchange}:

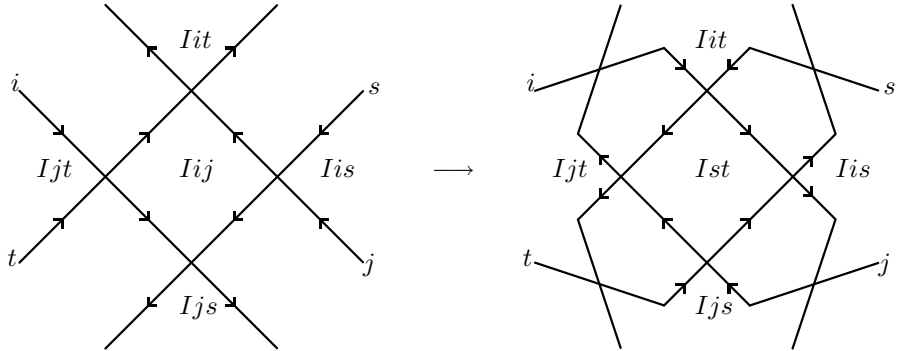
\begin{figure}[ht]

\setlength{\unitlength}{0.65pt}

\begin{center}
\begin{picture}(500,200)(50,50)

\thicklines

  \put(100,0){\line(1,1){150}}
  \put(125,25){\line(1,0){5}}
  \put(125,25){\line(0,1){5}}
  \put(175,75){\line(1,0){5}}
  \put(175,75){\line(0,1){5}}
  \put(225,125){\line(1,0){5}}
  \put(225,125){\line(0,1){5}}

  \put(50,50){\line(1,1){150}}
  \put(75,75){\line(0,-1){5}}
  \put(75,75){\line(-1,0){5}}
  \put(125,125){\line(0,-1){5}}
  \put(125,125){\line(-1,0){5}}
  \put(175,175){\line(0,-1){5}}
  \put(175,175){\line(-1,0){5}}

  \put(100,200){\line(1,-1){150}}
  \put(125,175){\line(0,-1){5}}
  \put(125,175){\line(1,0){5}}
  \put(175,125){\line(0,-1){5}}
  \put(175,125){\line(1,0){5}}
  \put(225,75){\line(0,-1){5}}
  \put(225,75){\line(1,0){5}}

  \put(50,150){\line(1,-1){150}}
  \put(75,125){\line(0,1){5}}
  \put(70,125){\line(1,0){5}}
  \put(125,75){\line(0,1){5}}
  \put(120,75){\line(1,0){5}}  
  \put(175,25){\line(0,1){5}}
  \put(170,25){\line(1,0){5}}

  \put(45,150){$i$}
  \put(250,45){$j$}
  \put(253,148){$s$}
  \put(43,48){$t$}

  \put(143,100){$Iij$}
  \put(143,175){$Iit$}
  \put(143,20){$Ijs$}
  \put(225,100){$Iis$}
  \put(60,100){$Ijt$}

\thicklines

  \put(350,50){\line(3,-1){75}}
  \put(388,88){\line(1,0){5}}
  \put(388,88){\line(0,1){5}}
  \put(425,25){\line(1,1){100}}
  \put(475,75){\line(-1,0){5}}
  \put(475,75){\line(0,-1){5}}
  \put(511,111){\line(-1,0){5}}
  \put(511,111){\line(0,-1){5}}

  \put(375,75){\line(1,-3){25}}
  \put(375,75){\line(1,1){100}}
  \put(425,125){\line(0,1){5}}
  \put(425,125){\line(1,0){5}}
  \put(500,200){\line(1,-3){25}}
  \put(463,163){\line(1,0){5}}
  \put(463,163){\line(0,1){5}} 
  \put(437,37){\line(-1,0){5}}
  \put(437,37){\line(0,-1){5}}

  \put(425,175){\line(1,-1){100}}
  \put(437,163){\line(0,1){5}}
  \put(437,163){\line(-1,0){5}}
  \put(475,125){\line(0,1){5}}
  \put(475,125){\line(-1,0){5}}
  \put(475,175){\line(3,-1){75}}
  \put(511,89){\line(0,1){5}}
  \put(511,89){\line(-1,0){5}}

  \put(375,125){\line(1,-1){100}}
  \put(475,25){\line(3,1){75}}
  \put(375,125){\line(1,3){25}}
  \put(350,150){\line(3,1){75}}
  \put(500,0){\line(1,3){25}}
  \put(388,112){\line(0,-1){5}}
  \put(388,112){\line(1,0){5}}
  \put(425,75){\line(0,-1){5}}
  \put(425,75){\line(1,0){5}}  
  \put(463,37){\line(0,-1){5}}
  \put(463,37){\line(1,0){5}}

  \put(345,150){$i$}
  \put(550,45){$j$}
  \put(553,148){$s$}
  \put(343,48){$t$}

  \put(443,100){$Ist$}
  \put(443,175){$Iit$}
  \put(443,20){$Ijs$}
  \put(525,100){$Iis$}
  \put(360,100){$Ijt$}

  \put(290,100){${\bf \longrightarrow}$}

\end{picture}

\bigskip
\bigskip
\bigskip
\bigskip
\caption{Geometric Exchange}
\label{exchange}
\end{center}
\end{figure}

\bigskip
\bigskip
\bigskip
\bigskip
\bigskip
\noindent
In the picture above, $|I|=k-2$ and $i,j,s,t$ are distinct
indices disjoint from $I$. In this context local means
that a sufficiently small euclidean neighborhood can
be taken about the nexus of paths forming the quadrilateral
and it is within this neighborhood that the
rewiring is undertaken. The geometric exchange
is involutive provided we {\it untwist} - see above -
consecutive crossing after performing the exchange.
Starting with an initial arrangement
for any permutation $\pi$, Postnikov has 
shown that the 
remaining supply of $\pi$-diagrams is 
obtained by repeated use of the exchange.

\bigskip
\begin{Prop}\label{ch2,Post.2}[Postnikov]
Let ${\bf A}$ and ${\bf A'}$ two $\pi_{k,n}$-diagrams.
Then there is a sequence of geometric exchanges transforming
${\bf A}$ into ${\bf A'}$.
\end{Prop}

\bigskip
\section*{ {\it Parameterizing $\pi_{k,n}$-Diagrams}} 

\noindent
As in Section 2, let ${\bf P}$ be a regular $n$-gon with vertices 
labeled clockwise and in order by the indices $1 \dots n$.

\bigskip
\begin{Def}
Two $k$-subsets $I$ and $J$ in $[1 \dots n]$ are said to be 
{\bf non-crossing} if no chord in ${\bf P}$ having end points labeled
by indices from $I-J$ crosses any chord with end points labeled
by indices from $J-I$.
\end{Def}

\bigskip
\noindent
This property was first introduced in \cite{Leclerc} 
under the name {\it weak separablity} and was 
studied in connection with {\it quantum flag minors}.
In a related work \cite{Scott}
determines a neccessary and sufficient condition
for pairs of quantum minors in
$\Bbb{C}_q \Big[ Mat_{k,n} \Big]$
to {\it quasi-commute} using the combinatorics of
non-crossing sets. 

\bigskip
\noindent
Scott's work in \cite{Scott}
also focuses on maximal (by inclusion) collections of 
pairwise non-crossing $k$-subsets of $[1 \dots n]$.
Scott introduces
a family of involutive transformations
- called {\it $(2,4)$-exchanges} - on the
set of all maximal families of pairwise non-crossing
$k$-subsets defined in the following manner:
given a maximal collection $\mathcal{C}$
and $k$-subsets $Iij$, $Iis$, $Isj$, $Ijt$,
and $Iit$ in $\mathcal{C}$, let $\mathcal{C}'$
be the family obtained by replacing
$Iij$ with $Ist$. The pairs $\{i,j\}$ and 
$\{s,t\}$ are required to be disjoint and
{\bf crossing} and $I$ is any subset 
disjoint from $\{i,j,s,t\}$ of size $k-2$.
Theorem 3 of \cite{Scott} proves that
$\mathcal{C}'$ is itself a maximal pairwise
non-crossing family and, in the special
cases of $k=2$ and $k=3$, Theorems 4 and 5 of
\cite{Scott} prove that any two maximal
collections can be obtained from one another
by a sequence of such transformations. 

\bigskip
\noindent
Geometric exchanges are examples of $(2,4)$-exchanges
when the maximal family is the collection of labeling
sets for a $\pi_{k,n}$-diagram. Section 4 of this
paper constructs a $\pi_{k,n}$-diagram ${\bf A_{k,n}}$
which, among other features, has the property that
the collection of its labeling sets is a maximal
non-crossing family of $k$-subsets. Every $\pi_{k,n}$-
diagram obtained from ${\bf A_{k,n}}$ through
a series of geometric exchanges also must have
the property that its labeling sets form
a maximal non-crossing family. Applying Proposition 
\ref{ch2,Post.2} of
Postnikov we arrive at the following corollary and
conjecture:

\bigskip
\begin{Cor}\label{ch2,cor.max}
The collection of labeling sets of a $\pi_{k,n}$-diagram
${\bf D}$ is a maximal (by inclusion) family
of pairwise non-crossing $k$-subsets.
\end{Cor}   

\bigskip
\begin{Conj}\label{ch2,conj.max}
Every maximal family of non-crossing 
$k$-subsets is a collection of labeling sets
for some $\pi_{k,n}$-diagram ${\bf D}$.
In particular, every maximal family consists
of $k(n-k)+1$ subsets.
\end{Conj}

\bigskip
\noindent
The assertion in Conjecture \ref{ch2,conj.max} that every maximal family
has exactly $k(n-k)+1$ members is a conjecture formulated 
and addressed at length by Scott in \cite{Scott};
there is goes under the heading {\it Purity} 
Conjecture. Proposition 1
of \cite{Scott} proves that the size of any such
family can not exceed the sharp bound of $k(n-k)+1$.
Since $(2,4)$-exchanges do not alter the cardinality
of a maximal family, one can prove purity by showing 
that $(2,4)$-exchanges act transitively on the set of all 
maximal non-crossing families. Accordingly Conjecture 2
of \cite{Scott} asserts this transitivity claim and
proves it in Theorems 4 and 5 in special
cases of $k=2$ and $k=3$.      

\bigskip
\noindent
For $k=2$ and $k=3$, one can show that every $(2,4)$-exchange which
can be performed on a family of labeling sets
for a Postnikov arrangement is realized as a geometric exchange. 
Together with Theorems 4 and 5 of \cite{Scott}, this
establishes Conjecture 1 for these special values of $k$ 
and allows us to identify (and thus parameterize)
a $\pi_{k,n}$-diagram with its collection of labeling
sets. 

\bigskip
\noindent
For $k=2$ the correspondence between maximal
families and diagrams can be made even more explicit.
In this case each maximal family is represented
as a maximal family of non-crossing chords of
a labeled regular $n$-gon ${\bf P}$; 
i.e. a triangulation of ${\bf P}$.
Given a triangulation ${\bf T}$ of ${\bf P}$ construct
a $\pi_{2,n}$-diagram ${\bf D(T)}$ by replacing every
triangle $\{a,b,c\}$ of ${\bf T}$ with the 
following local path configuration:

\begin{figure}[ht]

\setlength{\unitlength}{.70pt}

\begin{center}
\begin{picture}(450,100)(50,50)

\thicklines
\put(50,50){\line(1,0){200}}
\put(50,50){\line(1,1){100}}
\put(150,150){\line(1,-1){100}}

\put(42,42){$a$}
\put(148,152){$b$}
\put(252,42){$c$}

\put(325,75){\line(1,0){150}}
\put(347,75){\line(-1,1){3}}
\put(347,75){\line(-1,-1){3}}
\put(400,75){\line(-1,1){3}}
\put(400,75){\line(-1,-1){3}}
\put(460,75){\line(-1,1){3}}
\put(460,75){\line(-1,-1){3}}

\put(380,130){\line(1,-1){80}}
\put(390,120){\line(0,-1){3}}
\put(390,120){\line(1,0){3}}
\put(410,100){\line(0,-1){3}}
\put(410,100){\line(1,0){3}}
\put(445,65){\line(0,-1){3}}
\put(445,65){\line(1,0){3}}

\put(420,130){\line(-1,-1){80}}
\put(410,120){\line(0,1){3}}
\put(410,120){\line(1,0){3}}
\put(390,100){\line(0,1){3}}
\put(390,100){\line(1,0){3}}
\put(355,65){\line(0,1){3}}
\put(355,65){\line(1,0){3}}

\thinlines
\put(300,50){\line(1,0){200}}
\put(300,50){\line(1,1){100}}
\put(400,150){\line(1,-1){100}}

\put(292,42){$a$}
\put(398,152){$b$}
\put(502,42){$c$}

\end{picture}
\end{center}
\end{figure}

\newpage
\bigskip
\noindent
Whenever two triangles in ${\bf T}$ share an edge 
attach both pairs of oriented paths along
this edge. Label the path to the immediate right
of the vertex $i$ by $\pi_{1,n}(i)$ and 
the path to the immediate left by 
$\pi_{1,n}(i)'$. The following example illustrates
this construction:

\begin{figure}[h]
\begin{center}
\includegraphics[width=4.5in]{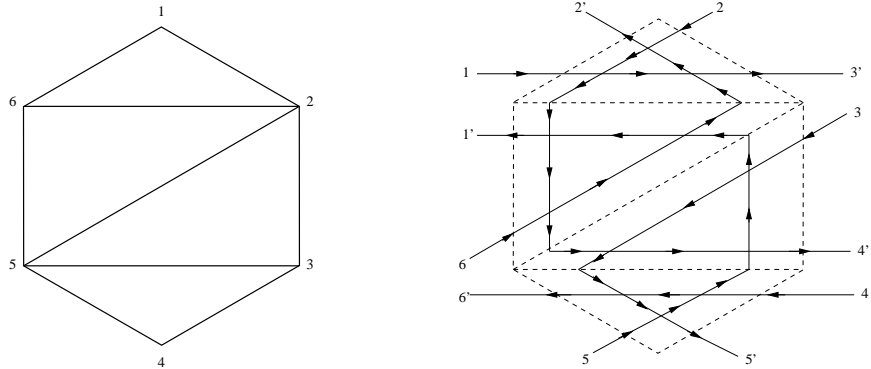}
\caption{Triangulation and Postnikov Arrangement}
\label{triangulation and arrangement}
\end{center}
\end{figure}

\noindent
Let ${\bf T_n}$ denotes the triangulation
consisting of the internal chords $[13],[14],[15],\dots
\newline [1 \ n-1]$ and boundary chords $[12],[23],\dots,[1n]$.
Let ${\bf D_n}$ denote the $\pi_{2,n}$-diagram
given by the following crossing data:
For $3 \leq i \leq n-1$ the path originating at
vertex $i$ first crosses path $1$, followed
by path $i-1$, followed by $i+1$, and finally
again by path $1$. The path originating at
vertex $2$ first crosses path $3$ and then
path $1$. Path $n$ crosses path $1$ and
then $n-1$. Path $1$ crosses path $n-1$, 
$n$, $n-2$, $n-1$, $n-3$, $n-2$, $n-4$, $n-3$,
$\dots$, $3$, $4$, $2$, and finally $3$.
If the above construction is applied to ${\bf T_n}$
then ${\bf D(T_n) = D_n}$. 
The fact that this construction actually
produces a bona fide $\pi_{2,n}$-diagram
for any triangulation ${\bf T}$ follows from the next Lemma.

\bigskip
\begin{lemma}\label{ch2,lemma.tri}
Let ${\bf T}$ and ${\bf T'}$ be triangulations
of ${\bf P}$ related by a single
{\it quadrilateral flip}. If ${\bf D(T)}$ is 
a $\pi_{2,n}$ diagram then so is ${\bf D(T')}$.
Moreover, ${\bf D(T)}$ and ${\bf D(T')}$ are
related by a single {\it geometric exchange}.
\end{lemma}

\newpage
\begin{proof}

\noindent
Let $\{i,s,j,t\}$ be the vertex labels
of the quadrilateral in which
the diagonal chord $[ij]$ in ${\bf T}$
is flipped to $[st]$ in ${\bf T'}$. Applying
the construction to the quadrilateral
in both cases yields the predicted
geometric exchange. 
 
\begin{figure}[h]
\begin{center}
\includegraphics[width=3.5in]{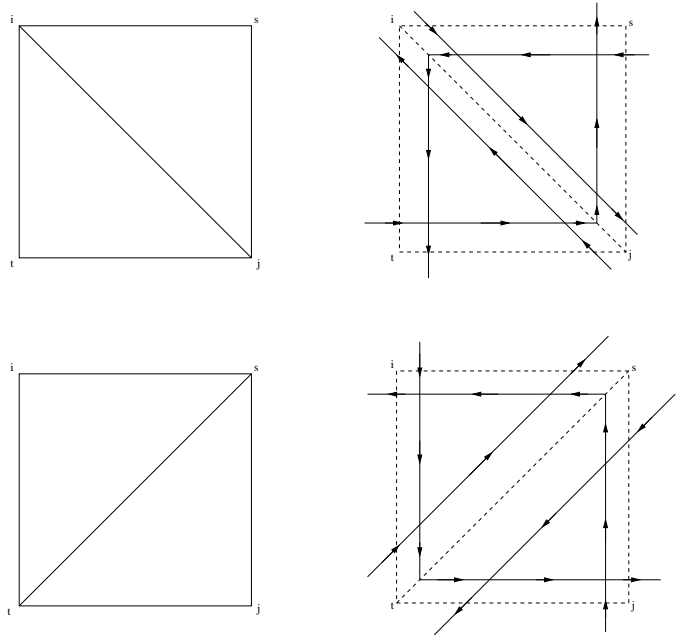}
\caption{Quadrilateral Flip and Geometric Exchange}
\label{flip and exchange}
\end{center} 
\end{figure}

\end{proof}

\bigskip
\noindent
It is well known that any two triangulations
${\bf T}$ and ${\bf T'}$ can be obtained
from one another by a sequence
of quadrilateral flips; in particular
any triangulations ${\bf T}$ can
be obtained from ${\bf T_n}$; this 
proves that ${\bf D(T)}$ 
is $\pi_{2,n}$-diagram for any triangulation
${\bf T}$.

\bigskip
\noindent
The construction has the property
that each chord $[ij]$
of a triangulation ${\bf T}$ is
the labeling set of the
quadrilateral cell in the diagram
${\bf D(T)}$ containing this chord. 
Consequently ${\bf T = T'}$ whenever
${\bf D(T)}$ and ${\bf D(T')}$
coincide. Taken together with Proposition
\ref{ch2,Post.2} of Postnikov we can deduce
the following corollary:

\bigskip
\begin{Cor}\label{ch2,cor.tri}
The construction establishes a bijection
${\bf T} \longrightarrow {\bf D(T)}$ 
between the set of all triangulations
of ${\bf P}$ and the set of all
$\pi_{2,n}$-diagrams.
\end{Cor}

\newpage
\bigskip
\section{Quadrilateral Postnikov Arrangements}

\bigskip
\noindent
This section describes the construction
of a $\pi_{k,n}$-diagram whose
internal even cells are quadrilateral.
This arrangement is neccessary for the proof of the 
chapter's principal result, Theorem \ref{ch2,thm.main}
of section 3.4. 
 
\bigskip
\noindent
In preparation for the next theorem, let ${\bf T_{k,n}}$ denote
the triangulation of the convex $n$-gon 
${\bf P}$ consisting of the following 
{\it zig-zag} configuration of internal chords:

\[ \Big\{ \ \big[ n-k+2 \ \ n-k \big], \
\big[ n-k \ \ n-k+3 \big], \ \big[n-k+3 \ \ n-k-1 \big], \
\big[ n-k-1 \ \ n-k+4 \big], \ \dots \ \Big\} \]

\bigskip
\noindent
The indices defining these chords are taken
{\bf mod $n$}. The following picture
illustrates this triangulation for $k=3$ and $n=8$.

\bigskip
\begin{figure}[h] 
\begin{center} 
\setlength{\unitlength}{2.8pt} 
\begin{picture}(60,60)(0,0) 
\thicklines 
  \multiput(0,20)(60,0){2}{\line(0,1){20}} 
  \multiput(20,0)(0,60){2}{\line(1,0){20}} 
  \multiput(0,40)(40,-40){2}{\line(1,1){20}} 
  \multiput(20,0)(40,40){2}{\line(-1,1){20}} 
 
  \multiput(20,0)(20,0){2}{\circle*{1}} 
  \multiput(20,60)(20,0){2}{\circle*{1}} 
  \multiput(0,20)(0,20){2}{\circle*{1}} 
  \multiput(60,20)(0,20){2}{\circle*{1}} 
 
\thinlines \put(0,20){\line(1,0){60}} \put(0,40){\line(1,0){60}} 
\put(20,0){\line(2,1){40}} 
\put(0,20){\line(3,1){60}}
\put(0,40){\line(2,1){40}}

\put(-4,19){$8$}
\put(-4,39){$1$} 
\put(62,19){$5$}
\put(62,39){$4$}
\put(19,-4){$7$}
\put(19,62){$2$}
\put(39,-4){$6$}
\put(39,62){$3$}

\end{picture} 
 
\end{center} 
\caption{The Triangulation {\bf $T_{3,8}$} } 
\label{snake} 
\end{figure}
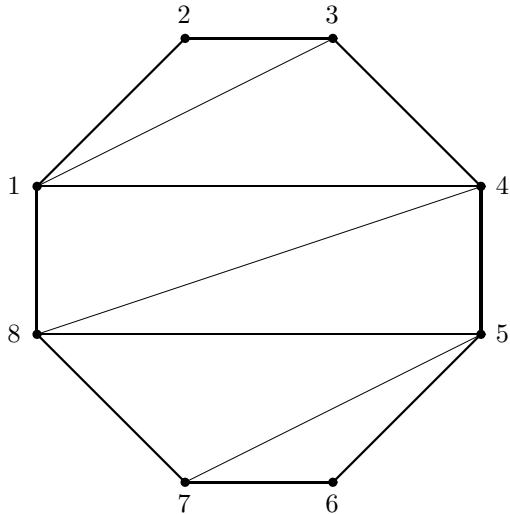

\begin{Thm}\label{ch2,thm.quad}
For positive integers $k$ and $n$ with $n \geq k+2
\geq 4$ there exists a $\pi_{k,n}$-diagram,
denoted as ${\bf A_{k,n}}$,
whose internal even cells are quadrilateral.
The collection of $k$-subset labels corresponding
to the internal cells consists of $k$-subsets
of $[1 \dots n]$ expressible as a disjoint
union $I \sqcup I'$ of {\it intervals} $I$ and $I'$
whose {\it beginning points} $i$ and $i'$ form
a chord $[ii']$ in the triangulation ${\bf T_{k,n}}$. 

\end{Thm}

\bigskip
\noindent
Intervals in $[1 \dots n]$ are subsets of the form 
$[i \dots i+ s]$ where $i+s$ is taken
{\bf mod $n$}; the index $i$ is the
beginning point. Before proving Theorem \ref{ch2,thm.quad}
let us again consider as an example
the case of $k=3$ and $n=8$. The labeled
quadrilateral arrangement is pictured below;
the $3$-subset labels have been 
partitioned by a comma to reveal the double interval.
As predicted, the beginning points of each double
interval form a chord in figure \ref{snake}.

\newpage
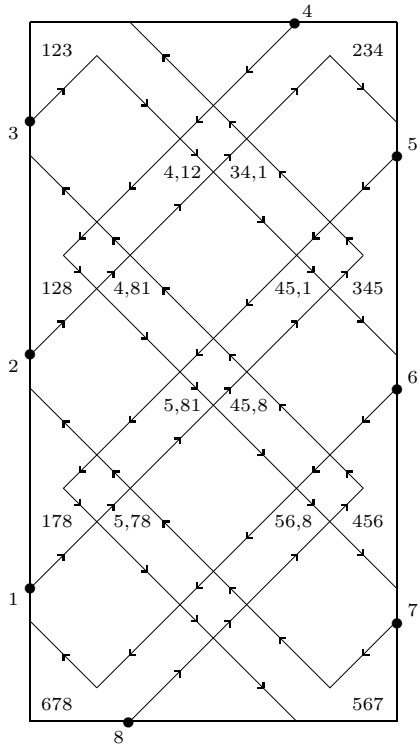
\begin{figure}[h]

\setlength{\unitlength}{2.10pt}

\begin{center}
\begin{picture}(-35,140)(50,50)
\thinlines

\put(0,0){\line(1,0){66}}
\put(0,0){\line(0,1){126}}
\put(0,126){\line(1,0){66}}
\put(66,0){\line(0,1){126}}

\put(0,18){\line(1,-1){12}}
\put(6,42){\line(1,-1){42}}
\put(0,60){\line(1,-1){54}}
\put(6,84){\line(1,-1){60}}
\put(0,102){\line(1,-1){60}}
\put(12,120){\line(1,-1){54}}
\put(18,126){\line(1,-1){42}}
\put(54,120){\line(1,-1){12}}

\put(6,42){\line(1,1){60}}
\put(0,24){\line(1,1){60}}
\put(12,6){\line(1,1){54}}
\put(18,0){\line(1,1){42}}
\put(54,6){\line(1,1){12}}
\put(0,66){\line(1,1){54}}
\put(6,84){\line(1,1){42}}
\put(0,108){\line(1,1){12}}

\put(9,45){\line(1,0){1}}
\put(9,45){\line(0,1){1}}
\put(18,54){\line(1,0){1}}
\put(18,54){\line(0,1){1}}
\put(18,12){\line(1,0){1}}
\put(18,12){\line(0,1){1}}
\put(30,24){\line(1,0){1}}
\put(30,24){\line(0,1){1}}
\put(39,33){\line(1,0){1}}
\put(39,33){\line(0,1){1}}
\put(51,45){\line(1,0){1}}
\put(51,45){\line(0,1){1}}
\put(60,54){\line(1,0){1}}
\put(60,54){\line(0,1){1}}
\put(60,12){\line(1,0){1}}
\put(60,12){\line(0,1){1}}
\put(60,96){\line(1,0){1}}
\put(60,96){\line(0,1){1}}
\put(51,87){\line(1,0){1}}
\put(51,87){\line(0,1){1}}

\put(9,87){\line(1,0){1}}
\put(9,87){\line(0,1){1}}
\put(18,96){\line(1,0){1}}
\put(18,96){\line(0,1){1}}
\put(18,54){\line(1,0){1}}
\put(18,54){\line(0,1){1}}
\put(30,66){\line(1,0){1}}
\put(30,66){\line(0,1){1}}
\put(39,75){\line(1,0){1}}
\put(39,75){\line(0,1){1}}

\put(6,30){\line(0,-1){1}}
\put(6,30){\line(-1,0){1}}
\put(15,39){\line(0,-1){1}}
\put(15,39){\line(-1,0){1}}
\put(27,51){\line(0,-1){1}}
\put(27,51){\line(-1,0){1}}
\put(27,9){\line(0,-1){1}}
\put(27,9){\line(-1,0){1}}
\put(36,18){\line(0,-1){1}}
\put(36,18){\line(-1,0){1}}
\put(48,30){\line(0,-1){1}}
\put(48,30){\line(-1,0){1}}
\put(57,39){\line(0,-1){1}}
\put(57,39){\line(-1,0){1}}
\put(36,60){\line(0,-1){1}}
\put(36,60){\line(-1,0){1}}
\put(27,93){\line(0,-1){1}}
\put(27,93){\line(-1,0){1}}
\put(6,72){\line(0,-1){1}}
\put(6,72){\line(-1,0){1}}
\put(15,81){\line(0,-1){1}}
\put(15,81){\line(-1,0){1}}
\put(48,72){\line(0,-1){1}}
\put(48,72){\line(-1,0){1}}
\put(6,114){\line(0,-1){1}}
\put(6,114){\line(-1,0){1}}
\put(27,93){\line(0,-1){1}}
\put(27,93){\line(-1,0){1}}
\put(36,102){\line(0,-1){1}}
\put(36,102){\line(-1,0){1}}

\put(6,12){\line(0,-1){1}}
\put(6,12){\line(1,0){1}}
\put(6,54){\line(0,-1){1}}
\put(6,54){\line(1,0){1}}
\put(15,45){\line(0,-1){1}}
\put(15,45){\line(1,0){1}}
\put(24,36){\line(0,-1){1}}
\put(24,36){\line(1,0){1}}
\put(36,24){\line(0,-1){1}}
\put(36,24){\line(1,0){1}}
\put(45,15){\line(0,-1){1}}
\put(45,15){\line(1,0){1}}
\put(45,57){\line(0,-1){1}}
\put(45,57){\line(1,0){1}}
\put(57,45){\line(0,-1){1}}
\put(57,45){\line(1,0){1}}
\put(6,96){\line(0,-1){1}}
\put(6,96){\line(1,0){1}}
\put(6,54){\line(0,-1){1}}
\put(6,54){\line(1,0){1}}
\put(15,87){\line(0,-1){1}}
\put(15,87){\line(1,0){1}}
\put(24,78){\line(0,-1){1}}
\put(24,78){\line(1,0){1}}
\put(36,66){\line(0,-1){1}}
\put(36,66){\line(1,0){1}}
\put(45,99){\line(0,-1){1}}
\put(45,99){\line(1,0){1}}
\put(57,87){\line(0,-1){1}}
\put(57,87){\line(1,0){1}}
\put(24,120){\line(0,-1){1}}
\put(24,120){\line(1,0){1}}
\put(36,108){\line(0,-1){1}}
\put(36,108){\line(1,0){1}}

\put(9,39){\line(0,1){1}}
\put(9,39){\line(-1,0){1}}
\put(21,27){\line(0,1){1}}
\put(21,27){\line(-1,0){1}}
\put(30,18){\line(0,1){1}}
\put(30,18){\line(-1,0){1}}
\put(42,6){\line(0,1){1}}
\put(42,6){\line(-1,0){1}}
\put(42,48){\line(0,1){1}}
\put(42,48){\line(-1,0){1}}
\put(51,39){\line(0,1){1}}
\put(51,39){\line(-1,0){1}}
\put(60,30){\line(0,1){1}}
\put(60,30){\line(-1,0){1}}
\put(60,72){\line(0,1){1}}
\put(60,72){\line(-1,0){1}}
\put(60,114){\line(0,1){1}}
\put(60,114){\line(-1,0){1}}
\put(9,81){\line(0,1){1}}
\put(9,81){\line(-1,0){1}}
\put(21,69){\line(0,1){1}}
\put(21,69){\line(-1,0){1}}
\put(30,60){\line(0,1){1}}
\put(30,60){\line(-1,0){1}}
\put(42,90){\line(0,1){1}}
\put(42,90){\line(-1,0){1}}
\put(51,81){\line(0,1){1}}
\put(51,81){\line(-1,0){1}}
\put(21,111){\line(0,1){1}}
\put(21,111){\line(-1,0){1}}
\put(30,102){\line(0,1){1}}
\put(30,102){\line(-1,0){1}}

\put(39,117){\line(0,1){1}}
\put(39,117){\line(1,0){1}}
\put(30,108){\line(0,1){1}}
\put(30,108){\line(1,0){1}}
\put(48,114){\line(0,-1){1}}
\put(48,114){\line(-1,0){1}}
\put(57,81){\line(0,-1){1}}
\put(57,81){\line(-1,0){1}}

\put(16.5,-1.3){$\bullet$}
\put(-1.2,22.7){$\bullet$}
\put(64.9,58.65){$\bullet$}
\put(64.9,16.5){$\bullet$}
\put(-1.2,64.9){$\bullet$}
\put(-1.2,106.9){$\bullet$}
\put(46.5,124.5){$\bullet$}
\put(64.9,100.65){$\bullet$}

\put(15,-4){$\p{8}$}
\put(-4,21){$\p{1}$}
\put(68,61){$\p{6}$}
\put(68,19){$\p{7}$}
\put(-4,63){$\p{2}$}
\put(-4,105){$\p{3}$}
\put(68,103){$\p{5}$}
\put(49,127){$\p{4}$}

\put(2,35){$\p{178}$}
\put(2,77){$\p{128}$}
\put(2,120){$\p{123}$}
\put(58,120){$\p{234}$}
\put(58,77){$\p{345}$}
\put(58,35){$\p{456}$}
\put(58,2){$\p{567}$}
\put(2,2){$\p{678}$}

\put(15,35){$\p{5,78}$}
\put(15,77){$\p{4,81}$}
\put(44,35){$\p{56,8}$}
\put(44,77){$\p{45,1}$}

\put(24,55.9){$\p{5,81}$}
\put(24,97.9){$\p{4,12}$}
\put(36,97.9){$\p{34,1}$}
\put(36,55.9){$\p{45,8}$}

\end{picture}
\bigskip
\bigskip
\bigskip
\bigskip
\bigskip
\bigskip
\bigskip
\bigskip
\bigskip
\bigskip
\bigskip
\bigskip
\bigskip
\caption{Quadrilateral Arrangement ${ \bf A_{3,8}}$}
\label{3-8 arrangement}
\end{center}
\end{figure}

\newpage
\bigskip
\begin{proof}
I will assume
that $k$ is odd; the construction 
for $k$ even is nearly identical except for
some minor alterations. Begin with $n=k+2$.
Visual inspection reveals that the 
internal cells of following
diagram are quadrilateral. 

\begin{figure}[ht]

\setlength{\unitlength}{2.10pt}

\begin{center}
\begin{picture}(50,30)(50,50)
\thinlines

\put(0,0){\line(1,0){99}}
\put(114,0){\line(1,0){36}}
\put(0,0){\line(0,1){63}}
\put(0,63){\line(1,0){99}}
\put(114,63){\line(1,0){36}}
\put(150,0){\line(0,1){63}}

\put(0,18){\line(1,-1){12}}
\put(6,42){\line(1,-1){42}}
\put(0,60){\line(1,-1){54}}
\put(33,57){\line(1,-1){57}}
\put(39,63){\line(1,-1){57}}
\put(75,57){\line(1,-1){24}}
\put(81,63){\line(1,-1){18}}
\put(114,18){\line(1,-1){18}}
\put(114,30){\line(1,-1){24}}
\put(117,57){\line(1,-1){33}}
\put(123,63){\line(1,-1){21}}

\put(6,42){\line(1,1){21}}
\put(0,24){\line(1,1){33}}
\put(12,6){\line(1,1){57}}
\put(18,0){\line(1,1){57}}
\put(54,6){\line(1,1){45}}
\put(60,0){\line(1,1){39}}
\put(96,6){\line(1,1){3}}
\put(114,54){\line(1,1){3}}
\put(114,24){\line(1,1){36}}
\put(114,12){\line(1,1){30}}
\put(138,6){\line(1,1){12}}

\put(9,45){\line(1,0){1}}
\put(9,45){\line(0,1){1}}
\put(18,54){\line(1,0){1}}
\put(18,54){\line(0,1){1}}
\put(18,12){\line(1,0){1}}
\put(18,12){\line(0,1){1}}
\put(30,24){\line(1,0){1}}
\put(30,24){\line(0,1){1}}
\put(39,33){\line(1,0){1}}
\put(39,33){\line(0,1){1}}
\put(51,45){\line(1,0){1}}
\put(51,45){\line(0,1){1}}
\put(60,54){\line(1,0){1}}
\put(60,54){\line(0,1){1}}
\put(60,12){\line(1,0){1}}
\put(60,12){\line(0,1){1}}
\put(72,24){\line(1,0){1}}
\put(72,24){\line(0,1){1}}
\put(81,33){\line(1,0){1}}
\put(81,33){\line(0,1){1}}
\put(93,45){\line(1,0){1}}
\put(93,45){\line(0,1){1}}
\put(123,33){\line(1,0){1}}
\put(123,33){\line(0,1){1}}
\put(135,45){\line(1,0){1}}
\put(135,45){\line(0,1){1}}
\put(144,54){\line(1,0){1}}
\put(144,54){\line(0,1){1}}
\put(144,12){\line(1,0){1}}
\put(144,12){\line(0,1){1}}

\put(6,30){\line(0,-1){1}}
\put(6,30){\line(-1,0){1}}
\put(15,39){\line(0,-1){1}}
\put(15,39){\line(-1,0){1}}
\put(27,51){\line(0,-1){1}}
\put(27,51){\line(-1,0){1}}
\put(27,9){\line(0,-1){1}}
\put(27,9){\line(-1,0){1}}
\put(36,18){\line(0,-1){1}}
\put(36,18){\line(-1,0){1}}
\put(48,30){\line(0,-1){1}}
\put(48,30){\line(-1,0){1}}
\put(57,39){\line(0,-1){1}}
\put(57,39){\line(-1,0){1}}
\put(69,51){\line(0,-1){1}}
\put(69,51){\line(-1,0){1}}
\put(69,9){\line(0,-1){1}}
\put(69,9){\line(-1,0){1}}
\put(78,18){\line(0,-1){1}}
\put(78,18){\line(-1,0){1}}
\put(90,30){\line(0,-1){1}}
\put(90,30){\line(-1,0){1}}
\put(120,18){\line(0,-1){1}}
\put(120,18){\line(-1,0){1}}
\put(132,30){\line(0,-1){1}}
\put(132,30){\line(-1,0){1}}
\put(141,39){\line(0,-1){1}}
\put(141,39){\line(-1,0){1}}

\put(6,12){\line(0,-1){1}}
\put(6,12){\line(1,0){1}}
\put(6,54){\line(0,-1){1}}
\put(6,54){\line(1,0){1}}
\put(15,45){\line(0,-1){1}}
\put(15,45){\line(1,0){1}}
\put(24,36){\line(0,-1){1}}
\put(24,36){\line(1,0){1}}
\put(36,24){\line(0,-1){1}}
\put(36,24){\line(1,0){1}}
\put(45,15){\line(0,-1){1}}
\put(45,15){\line(1,0){1}}
\put(45,57){\line(0,-1){1}}
\put(45,57){\line(1,0){1}}
\put(57,45){\line(0,-1){1}}
\put(57,45){\line(1,0){1}}
\put(66,36){\line(0,-1){1}}
\put(66,36){\line(1,0){1}}
\put(78,24){\line(0,-1){1}}
\put(78,24){\line(1,0){1}}
\put(87,15){\line(0,-1){1}}
\put(87,15){\line(1,0){1}}
\put(87,57){\line(0,-1){1}}
\put(87,57){\line(1,0){1}}
\put(120,24){\line(0,-1){1}}
\put(120,24){\line(1,0){1}}
\put(129,15){\line(0,-1){1}}
\put(129,15){\line(1,0){1}}
\put(129,57){\line(0,-1){1}}
\put(129,57){\line(1,0){1}}
\put(141,45){\line(0,-1){1}}
\put(141,45){\line(1,0){1}}

\put(9,39){\line(0,1){1}}
\put(9,39){\line(-1,0){1}}
\put(21,27){\line(0,1){1}}
\put(21,27){\line(-1,0){1}}
\put(30,18){\line(0,1){1}}
\put(30,18){\line(-1,0){1}}
\put(42,6){\line(0,1){1}}
\put(42,6){\line(-1,0){1}}
\put(42,48){\line(0,1){1}}
\put(42,48){\line(-1,0){1}}
\put(51,39){\line(0,1){1}}
\put(51,39){\line(-1,0){1}}
\put(63,27){\line(0,1){1}}
\put(63,27){\line(-1,0){1}}
\put(72,18){\line(0,1){1}}
\put(72,18){\line(-1,0){1}}
\put(84,6){\line(0,1){1}}
\put(84,6){\line(-1,0){1}}
\put(84,48){\line(0,1){1}}
\put(84,48){\line(-1,0){1}}
\put(93,39){\line(0,1){1}}
\put(93,39){\line(-1,0){1}}
\put(126,6){\line(0,1){1}}
\put(126,6){\line(-1,0){1}}
\put(126,48){\line(0,1){1}}
\put(126,48){\line(-1,0){1}}
\put(135,39){\line(0,1){1}}
\put(135,39){\line(-1,0){1}}
\put(144,30){\line(0,1){1}}
\put(144,30){\line(-1,0){1}}

\put(16.5,-1.3){$\bullet$}
\put(58.5,-1.3){$\bullet$}
\put(-1.2,22.7){$\bullet$}
\put(26,61.85){$\bullet$}
\put(68,61.85){$\bullet$}
\put(148.8,58.5){$\bullet$}
\put(148.8,16.5){$\bullet$}


\put(104,-1.1){$\cdots$}
\put(104,19.8){$\cdots$}
\put(104,40.7){$\cdots$}
\put(104,61.85){$\cdots$}

\put(14,-4){$\p{n}$}
\put(55,-4){$\p{n-1}$}
\put(-4,22){$\p{1}$}
\put(28,65){$\p{2}$}
\put(70,65){$\p{3}$}
\put(152,63){${k+3 \over 2}$}
\put(152,21){${k+5 \over 2}$}

\end{picture}
\bigskip
\bigskip
\bigskip
\bigskip
\bigskip
\bigskip
\bigskip
\bigskip
\bigskip
\bigskip
\bigskip
\bigskip
\bigskip
\caption{Quadrilateral Arrangement ${\bf A_{k,k+2}}$}
\label{quad1 arrangement}
\end{center}
\end{figure}
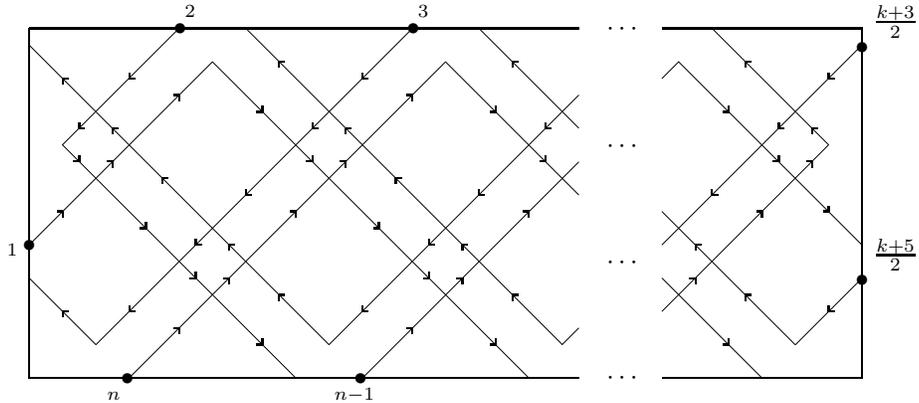

\bigskip
\noindent
It is easy to see that the conditions
stipulated in the definition of Postnikov
arrangements are satisfied. 
For instance, one sees that this diagram corresponds
to the permutation $\pi_{k,k+2}$ by noticing
that the number of sources (depicted as black nodes
on the boundary rectangle) between the ending 
position of a path and its starting position (measured
clockwise from the endpoint) is always $2$. 

\bigskip
\noindent
In order to construct the diagram ${\bf A_{k,n}}$ 
for $n>k+2$, adjoin
$n-k-2$ rectalinear {\bf bands} -
each consisting of ${k-1 \over 2} +1$ sources -
along the top of ${\bf A_{k,k+2}}$. These
bands, which come in two varieties (even and odd),
must be stacked on top of one another in 
alternation starting with even type. 
The illustrations below 
depict the cases of $n=k+3$ and $n=k+4$;
the corresponding bands are highlighted.

\newpage
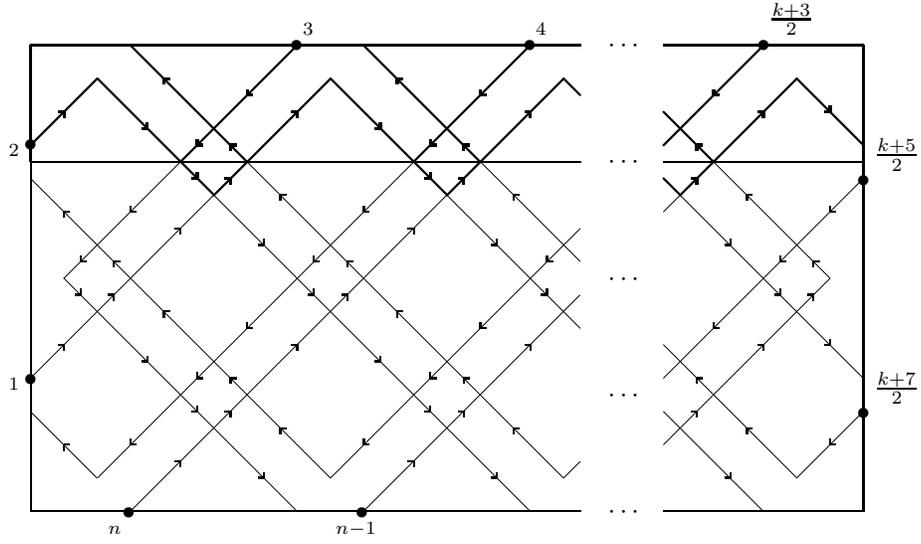
\begin{figure}[ht]

\setlength{\unitlength}{2.10pt}

\begin{center}
\begin{picture}(50,50)(50,50)
\thinlines

\put(0,0){\line(1,0){99}}
\put(114,0){\line(1,0){36}}
\put(0,0){\line(0,1){63}}
\put(0,63){\line(1,0){99}}
\put(114,63){\line(1,0){36}}
\put(150,0){\line(0,1){63}}

\put(0,18){\line(1,-1){12}}
\put(6,42){\line(1,-1){42}}
\put(0,60){\line(1,-1){54}}
\put(33,57){\line(1,-1){57}}
\put(39,63){\line(1,-1){57}}
\put(75,57){\line(1,-1){24}}
\put(81,63){\line(1,-1){18}}
\put(114,18){\line(1,-1){18}}
\put(114,30){\line(1,-1){24}}
\put(117,57){\line(1,-1){33}}
\put(123,63){\line(1,-1){21}}

\put(6,42){\line(1,1){21}}
\put(0,24){\line(1,1){33}}
\put(12,6){\line(1,1){57}}
\put(18,0){\line(1,1){57}}
\put(54,6){\line(1,1){45}}
\put(60,0){\line(1,1){39}}
\put(96,6){\line(1,1){3}}
\put(114,54){\line(1,1){3}}
\put(114,24){\line(1,1){36}}
\put(114,12){\line(1,1){30}}
\put(138,6){\line(1,1){12}}

\put(9,45){\line(1,0){1}}
\put(9,45){\line(0,1){1}}
\put(18,54){\line(1,0){1}}
\put(18,54){\line(0,1){1}}
\put(18,12){\line(1,0){1}}
\put(18,12){\line(0,1){1}}
\put(30,24){\line(1,0){1}}
\put(30,24){\line(0,1){1}}
\put(39,33){\line(1,0){1}}
\put(39,33){\line(0,1){1}}
\put(51,45){\line(1,0){1}}
\put(51,45){\line(0,1){1}}
\put(60,54){\line(1,0){1}}
\put(60,54){\line(0,1){1}}
\put(60,12){\line(1,0){1}}
\put(60,12){\line(0,1){1}}
\put(72,24){\line(1,0){1}}
\put(72,24){\line(0,1){1}}
\put(81,33){\line(1,0){1}}
\put(81,33){\line(0,1){1}}
\put(93,45){\line(1,0){1}}
\put(93,45){\line(0,1){1}}
\put(123,33){\line(1,0){1}}
\put(123,33){\line(0,1){1}}
\put(135,45){\line(1,0){1}}
\put(135,45){\line(0,1){1}}
\put(144,54){\line(1,0){1}}
\put(144,54){\line(0,1){1}}
\put(144,12){\line(1,0){1}}
\put(144,12){\line(0,1){1}}

\put(6,30){\line(0,-1){1}}
\put(6,30){\line(-1,0){1}}
\put(15,39){\line(0,-1){1}}
\put(15,39){\line(-1,0){1}}
\put(27,51){\line(0,-1){1}}
\put(27,51){\line(-1,0){1}}
\put(27,9){\line(0,-1){1}}
\put(27,9){\line(-1,0){1}}
\put(36,18){\line(0,-1){1}}
\put(36,18){\line(-1,0){1}}
\put(48,30){\line(0,-1){1}}
\put(48,30){\line(-1,0){1}}
\put(57,39){\line(0,-1){1}}
\put(57,39){\line(-1,0){1}}
\put(69,51){\line(0,-1){1}}
\put(69,51){\line(-1,0){1}}
\put(69,9){\line(0,-1){1}}
\put(69,9){\line(-1,0){1}}
\put(78,18){\line(0,-1){1}}
\put(78,18){\line(-1,0){1}}
\put(90,30){\line(0,-1){1}}
\put(90,30){\line(-1,0){1}}
\put(120,18){\line(0,-1){1}}
\put(120,18){\line(-1,0){1}}
\put(132,30){\line(0,-1){1}}
\put(132,30){\line(-1,0){1}}
\put(141,39){\line(0,-1){1}}
\put(141,39){\line(-1,0){1}}

\put(6,12){\line(0,-1){1}}
\put(6,12){\line(1,0){1}}
\put(6,54){\line(0,-1){1}}
\put(6,54){\line(1,0){1}}
\put(15,45){\line(0,-1){1}}
\put(15,45){\line(1,0){1}}
\put(24,36){\line(0,-1){1}}
\put(24,36){\line(1,0){1}}
\put(36,24){\line(0,-1){1}}
\put(36,24){\line(1,0){1}}
\put(45,15){\line(0,-1){1}}
\put(45,15){\line(1,0){1}}
\put(45,57){\line(0,-1){1}}
\put(45,57){\line(1,0){1}}
\put(57,45){\line(0,-1){1}}
\put(57,45){\line(1,0){1}}
\put(66,36){\line(0,-1){1}}
\put(66,36){\line(1,0){1}}
\put(78,24){\line(0,-1){1}}
\put(78,24){\line(1,0){1}}
\put(87,15){\line(0,-1){1}}
\put(87,15){\line(1,0){1}}
\put(87,57){\line(0,-1){1}}
\put(87,57){\line(1,0){1}}
\put(120,24){\line(0,-1){1}}
\put(120,24){\line(1,0){1}}
\put(129,15){\line(0,-1){1}}
\put(129,15){\line(1,0){1}}
\put(129,57){\line(0,-1){1}}
\put(129,57){\line(1,0){1}}
\put(141,45){\line(0,-1){1}}
\put(141,45){\line(1,0){1}}

\put(9,39){\line(0,1){1}}
\put(9,39){\line(-1,0){1}}
\put(21,27){\line(0,1){1}}
\put(21,27){\line(-1,0){1}}
\put(30,18){\line(0,1){1}}
\put(30,18){\line(-1,0){1}}
\put(42,6){\line(0,1){1}}
\put(42,6){\line(-1,0){1}}
\put(42,48){\line(0,1){1}}
\put(42,48){\line(-1,0){1}}
\put(51,39){\line(0,1){1}}
\put(51,39){\line(-1,0){1}}
\put(63,27){\line(0,1){1}}
\put(63,27){\line(-1,0){1}}
\put(72,18){\line(0,1){1}}
\put(72,18){\line(-1,0){1}}
\put(84,6){\line(0,1){1}}
\put(84,6){\line(-1,0){1}}
\put(84,48){\line(0,1){1}}
\put(84,48){\line(-1,0){1}}
\put(93,39){\line(0,1){1}}
\put(93,39){\line(-1,0){1}}
\put(126,6){\line(0,1){1}}
\put(126,6){\line(-1,0){1}}
\put(126,48){\line(0,1){1}}
\put(126,48){\line(-1,0){1}}
\put(135,39){\line(0,1){1}}
\put(135,39){\line(-1,0){1}}
\put(144,30){\line(0,1){1}}
\put(144,30){\line(-1,0){1}}

\put(16.5,-1.3){$\bullet$}
\put(58.5,-1.3){$\bullet$}
\put(-1.2,22.7){$\bullet$}
\put(148.8,58.5){$\bullet$}
\put(148.8,16.5){$\bullet$}

\put(104,-1.1){$\cdots$}
\put(104,19.8){$\cdots$}
\put(104,40.7){$\cdots$}
\put(104,61.85){$\cdots$}
\put(104,82.85){$\cdots$}


\thicklines

\put(0,63){\line(0,1){21}}
\put(0,84){\line(1,0){99}}
\put(114,84){\line(1,0){36}}
\put(150,63){\line(0,1){21}}

\put(0,66){\line(1,1){12}}
\put(27,63){\line(1,1){21}}
\put(33,57){\line(1,1){21}}
\put(69,63){\line(1,1){21}}
\put(75,57){\line(1,1){21}}
\put(114,66){\line(1,1){18}}
\put(117,57){\line(1,1){21}}

\put(12,78){\line(1,-1){21}}
\put(18,84){\line(1,-1){21}}
\put(54,78){\line(1,-1){21}}
\put(60,84){\line(1,-1){21}}
\put(96,78){\line(1,-1){3}}
\put(114,60){\line(1,-1){3}}
\put(114,72){\line(1,-1){9}}
\put(138,78){\line(1,-1){12}}

\put(6,72){\line(0,-1){1}}
\put(6,72){\line(-1,0){1}}
\put(48,72){\line(0,-1){1}}
\put(48,72){\line(-1,0){1}}
\put(90,72){\line(0,-1){1}}
\put(90,72){\line(-1,0){1}}
\put(132,72){\line(0,-1){1}}
\put(132,72){\line(-1,0){1}}
\put(36,60){\line(0,-1){1}}
\put(36,60){\line(-1,0){1}}
\put(78,60){\line(0,-1){1}}
\put(78,60){\line(-1,0){1}}
\put(120,60){\line(0,-1){1}}
\put(120,60){\line(-1,0){1}}

\put(39,75){\line(0,1){1}}
\put(39,75){\line(1,0){1}}
\put(81,75){\line(0,1){1}}
\put(81,75){\line(1,0){1}}
\put(123,75){\line(0,1){1}}
\put(123,75){\line(1,0){1}}
\put(30,66){\line(0,1){1}}
\put(30,66){\line(1,0){1}}
\put(72,66){\line(0,1){1}}
\put(72,66){\line(1,0){1}}

\put(21,69){\line(0,1){1}}
\put(21,69){\line(-1,0){1}}
\put(63,69){\line(0,1){1}}
\put(63,69){\line(-1,0){1}}
\put(144,72){\line(0,1){1}}
\put(144,72){\line(-1,0){1}}
\put(30,60){\line(0,1){1}}
\put(30,60){\line(-1,0){1}}
\put(72,60){\line(0,1){1}}
\put(72,60){\line(-1,0){1}}

\put(24,78){\line(1,0){1}}
\put(24,78){\line(0,-1){1}}
\put(66,78){\line(1,0){1}}
\put(66,78){\line(0,-1){1}}
\put(36,66){\line(1,0){1}}
\put(36,66){\line(0,-1){1}}
\put(78,66){\line(1,0){1}}
\put(78,66){\line(0,-1){1}}
\put(120,66){\line(1,0){1}}
\put(120,66){\line(0,-1){1}}

\put(-1.2,65){$\bullet$}
\put(46.7,82.85){$\bullet$}
\put(88.7,82.85){$\bullet$}
\put(130.8,82.85){$\bullet$}

\put(14,-4){$\p{n}$}
\put(55,-4){$\p{n-1}$}
\put(-4,22){$\p{1}$}
\put(-4,64){$\p{2}$}
\put(49,86){$\p{3}$}
\put(91,86){$\p{4}$}
\put(133,87.5){${k+3 \over 2}$}
\put(152,63){${k+5 \over 2}$}
\put(152,21){${k+7 \over 2}$}

\end{picture}

\bigskip
\bigskip
\bigskip
\bigskip
\bigskip
\bigskip
\bigskip
\bigskip
\bigskip
\bigskip
\bigskip
\bigskip
\bigskip
\caption{Quadrilateral Arrangement ${\bf A_{k,k+3}}$}
\label{quad arrangement 2}
\end{center}
\end{figure}

\begin{figure}[ht]
\setlength{\unitlength}{2.10pt}

\bigskip
\bigskip
\begin{center}
\begin{picture}(50,50)(50,50)
\thinlines

\put(0,0){\line(1,0){99}}
\put(114,0){\line(1,0){36}}
\put(0,0){\line(0,1){63}}
\put(150,0){\line(0,1){63}}

\put(0,18){\line(1,-1){12}}
\put(6,42){\line(1,-1){42}}
\put(0,60){\line(1,-1){54}}
\put(33,57){\line(1,-1){57}}
\put(39,63){\line(1,-1){57}}
\put(75,57){\line(1,-1){24}}
\put(81,63){\line(1,-1){18}}
\put(114,18){\line(1,-1){18}}
\put(114,30){\line(1,-1){24}}
\put(117,57){\line(1,-1){33}}
\put(123,63){\line(1,-1){21}}

\put(6,42){\line(1,1){21}}
\put(0,24){\line(1,1){33}}
\put(12,6){\line(1,1){57}}
\put(18,0){\line(1,1){57}}
\put(54,6){\line(1,1){45}}
\put(60,0){\line(1,1){39}}
\put(96,6){\line(1,1){3}}
\put(114,54){\line(1,1){3}}
\put(114,24){\line(1,1){36}}
\put(114,12){\line(1,1){30}}
\put(138,6){\line(1,1){12}}

\put(9,45){\line(1,0){1}}
\put(9,45){\line(0,1){1}}
\put(18,54){\line(1,0){1}}
\put(18,54){\line(0,1){1}}
\put(18,12){\line(1,0){1}}
\put(18,12){\line(0,1){1}}
\put(30,24){\line(1,0){1}}
\put(30,24){\line(0,1){1}}
\put(39,33){\line(1,0){1}}
\put(39,33){\line(0,1){1}}
\put(51,45){\line(1,0){1}}
\put(51,45){\line(0,1){1}}
\put(60,54){\line(1,0){1}}
\put(60,54){\line(0,1){1}}
\put(60,12){\line(1,0){1}}
\put(60,12){\line(0,1){1}}
\put(72,24){\line(1,0){1}}
\put(72,24){\line(0,1){1}}
\put(81,33){\line(1,0){1}}
\put(81,33){\line(0,1){1}}
\put(93,45){\line(1,0){1}}
\put(93,45){\line(0,1){1}}
\put(123,33){\line(1,0){1}}
\put(123,33){\line(0,1){1}}
\put(135,45){\line(1,0){1}}
\put(135,45){\line(0,1){1}}
\put(144,54){\line(1,0){1}}
\put(144,54){\line(0,1){1}}
\put(144,12){\line(1,0){1}}
\put(144,12){\line(0,1){1}}

\put(6,30){\line(0,-1){1}}
\put(6,30){\line(-1,0){1}}
\put(15,39){\line(0,-1){1}}
\put(15,39){\line(-1,0){1}}
\put(27,51){\line(0,-1){1}}
\put(27,51){\line(-1,0){1}}
\put(27,9){\line(0,-1){1}}
\put(27,9){\line(-1,0){1}}
\put(36,18){\line(0,-1){1}}
\put(36,18){\line(-1,0){1}}
\put(48,30){\line(0,-1){1}}
\put(48,30){\line(-1,0){1}}
\put(57,39){\line(0,-1){1}}
\put(57,39){\line(-1,0){1}}
\put(69,51){\line(0,-1){1}}
\put(69,51){\line(-1,0){1}}
\put(69,9){\line(0,-1){1}}
\put(69,9){\line(-1,0){1}}
\put(78,18){\line(0,-1){1}}
\put(78,18){\line(-1,0){1}}
\put(90,30){\line(0,-1){1}}
\put(90,30){\line(-1,0){1}}
\put(120,18){\line(0,-1){1}}
\put(120,18){\line(-1,0){1}}
\put(132,30){\line(0,-1){1}}
\put(132,30){\line(-1,0){1}}
\put(141,39){\line(0,-1){1}}
\put(141,39){\line(-1,0){1}}

\put(6,12){\line(0,-1){1}}
\put(6,12){\line(1,0){1}}
\put(6,54){\line(0,-1){1}}
\put(6,54){\line(1,0){1}}
\put(15,45){\line(0,-1){1}}
\put(15,45){\line(1,0){1}}
\put(24,36){\line(0,-1){1}}
\put(24,36){\line(1,0){1}}
\put(36,24){\line(0,-1){1}}
\put(36,24){\line(1,0){1}}
\put(45,15){\line(0,-1){1}}
\put(45,15){\line(1,0){1}}
\put(45,57){\line(0,-1){1}}
\put(45,57){\line(1,0){1}}
\put(57,45){\line(0,-1){1}}
\put(57,45){\line(1,0){1}}
\put(66,36){\line(0,-1){1}}
\put(66,36){\line(1,0){1}}
\put(78,24){\line(0,-1){1}}
\put(78,24){\line(1,0){1}}
\put(87,15){\line(0,-1){1}}
\put(87,15){\line(1,0){1}}
\put(87,57){\line(0,-1){1}}
\put(87,57){\line(1,0){1}}
\put(120,24){\line(0,-1){1}}
\put(120,24){\line(1,0){1}}
\put(129,15){\line(0,-1){1}}
\put(129,15){\line(1,0){1}}
\put(129,57){\line(0,-1){1}}
\put(129,57){\line(1,0){1}}
\put(141,45){\line(0,-1){1}}
\put(141,45){\line(1,0){1}}

\put(9,39){\line(0,1){1}}
\put(9,39){\line(-1,0){1}}
\put(21,27){\line(0,1){1}}
\put(21,27){\line(-1,0){1}}
\put(30,18){\line(0,1){1}}
\put(30,18){\line(-1,0){1}}
\put(42,6){\line(0,1){1}}
\put(42,6){\line(-1,0){1}}
\put(42,48){\line(0,1){1}}
\put(42,48){\line(-1,0){1}}
\put(51,39){\line(0,1){1}}
\put(51,39){\line(-1,0){1}}
\put(63,27){\line(0,1){1}}
\put(63,27){\line(-1,0){1}}
\put(72,18){\line(0,1){1}}
\put(72,18){\line(-1,0){1}}
\put(84,6){\line(0,1){1}}
\put(84,6){\line(-1,0){1}}
\put(84,48){\line(0,1){1}}
\put(84,48){\line(-1,0){1}}
\put(93,39){\line(0,1){1}}
\put(93,39){\line(-1,0){1}}
\put(126,6){\line(0,1){1}}
\put(126,6){\line(-1,0){1}}
\put(126,48){\line(0,1){1}}
\put(126,48){\line(-1,0){1}}
\put(135,39){\line(0,1){1}}
\put(135,39){\line(-1,0){1}}
\put(144,30){\line(0,1){1}}
\put(144,30){\line(-1,0){1}}

\put(16.5,-1.3){$\bullet$}
\put(58.5,-1.3){$\bullet$}
\put(-1.2,22.7){$\bullet$}
\put(148.8,58.5){$\bullet$}
\put(148.8,16.5){$\bullet$}



\put(0,63){\line(0,1){21}}
\put(0,84){\line(1,0){99}}
\put(114,84){\line(1,0){36}}
\put(150,63){\line(0,1){21}}

\put(0,66){\line(1,1){12}}
\put(27,63){\line(1,1){21}}
\put(33,57){\line(1,1){21}}
\put(69,63){\line(1,1){21}}
\put(75,57){\line(1,1){21}}
\put(114,66){\line(1,1){18}}
\put(117,57){\line(1,1){21}}

\put(12,78){\line(1,-1){21}}
\put(18,84){\line(1,-1){21}}
\put(54,78){\line(1,-1){21}}
\put(60,84){\line(1,-1){21}}
\put(96,78){\line(1,-1){3}}
\put(114,60){\line(1,-1){3}}
\put(114,72){\line(1,-1){9}}
\put(138,78){\line(1,-1){12}}

\put(6,72){\line(0,-1){1}}
\put(6,72){\line(-1,0){1}}
\put(48,72){\line(0,-1){1}}
\put(48,72){\line(-1,0){1}}
\put(90,72){\line(0,-1){1}}
\put(90,72){\line(-1,0){1}}
\put(132,72){\line(0,-1){1}}
\put(132,72){\line(-1,0){1}}
\put(36,60){\line(0,-1){1}}
\put(36,60){\line(-1,0){1}}
\put(78,60){\line(0,-1){1}}
\put(78,60){\line(-1,0){1}}
\put(120,60){\line(0,-1){1}}
\put(120,60){\line(-1,0){1}}

\put(39,75){\line(0,1){1}}
\put(39,75){\line(1,0){1}}
\put(81,75){\line(0,1){1}}
\put(81,75){\line(1,0){1}}
\put(123,75){\line(0,1){1}}
\put(123,75){\line(1,0){1}}
\put(30,66){\line(0,1){1}}
\put(30,66){\line(1,0){1}}
\put(72,66){\line(0,1){1}}
\put(72,66){\line(1,0){1}}

\put(21,69){\line(0,1){1}}
\put(21,69){\line(-1,0){1}}
\put(63,69){\line(0,1){1}}
\put(63,69){\line(-1,0){1}}
\put(144,72){\line(0,1){1}}
\put(144,72){\line(-1,0){1}}
\put(30,60){\line(0,1){1}}
\put(30,60){\line(-1,0){1}}
\put(72,60){\line(0,1){1}}
\put(72,60){\line(-1,0){1}}

\put(24,78){\line(1,0){1}}
\put(24,78){\line(0,-1){1}}
\put(66,78){\line(1,0){1}}
\put(66,78){\line(0,-1){1}}
\put(36,66){\line(1,0){1}}
\put(36,66){\line(0,-1){1}}
\put(78,66){\line(1,0){1}}
\put(78,66){\line(0,-1){1}}
\put(120,66){\line(1,0){1}}
\put(120,66){\line(0,-1){1}}

\put(-1.2,65){$\bullet$}

\thicklines

\put(0,84){\line(0,1){21}}
\put(0,105){\line(1,0){99}}
\put(114,105){\line(1,0){36}}
\put(150,105){\line(0,-1){21}}

\put(6,84){\line(1,1){21}}
\put(12,78){\line(1,1){21}}
\put(48,84){\line(1,1){21}}
\put(54,78){\line(1,1){21}}
\put(90,84){\line(1,1){9}}
\put(96,78){\line(1,1){3}}
\put(114,96){\line(1,1){3}}
\put(132,84){\line(1,1){18}}
\put(138,78){\line(1,1){6}}

\put(6,84){\line(1,-1){6}}
\put(0,102){\line(1,-1){18}}
\put(33,99){\line(1,-1){21}}
\put(39,105){\line(1,-1){21}}
\put(75,99){\line(1,-1){21}}
\put(81,105){\line(1,-1){18}}
\put(117,99){\line(1,-1){21}}
\put(123,105){\line(1,-1){21}}

\put(15,81){\line(0,-1){1}}
\put(15,81){\line(-1,0){1}}
\put(57,81){\line(0,-1){1}}
\put(57,81){\line(-1,0){1}}
\put(141,81){\line(0,-1){1}}
\put(141,81){\line(-1,0){1}}
\put(27,93){\line(0,-1){1}}
\put(27,93){\line(-1,0){1}}
\put(69,93){\line(0,-1){1}}
\put(69,93){\line(-1,0){1}}

\put(9,87){\line(0,1){1}}
\put(9,87){\line(1,0){1}}
\put(51,87){\line(0,1){1}}
\put(51,87){\line(1,0){1}}
\put(93,87){\line(0,1){1}}
\put(93,87){\line(1,0){1}}
\put(135,87){\line(0,1){1}}
\put(135,87){\line(1,0){1}}
\put(18,96){\line(0,1){1}}
\put(18,96){\line(1,0){1}}
\put(60,96){\line(0,1){1}}
\put(60,96){\line(1,0){1}}
\put(144,96){\line(0,1){1}}
\put(144,96){\line(1,0){1}}

\put(9,81){\line(0,1){1}}
\put(9,81){\line(-1,0){1}}
\put(51,81){\line(0,1){1}}
\put(51,81){\line(-1,0){1}}
\put(93,81){\line(0,1){1}}
\put(93,81){\line(-1,0){1}}
\put(135,81){\line(0,1){1}}
\put(135,81){\line(-1,0){1}}
\put(42,90){\line(0,1){1}}
\put(42,90){\line(-1,0){1}}
\put(84,90){\line(0,1){1}}
\put(84,90){\line(-1,0){1}}
\put(126,90){\line(0,1){1}}
\put(126,90){\line(-1,0){1}}

\put(15,87){\line(1,0){1}}
\put(15,87){\line(0,-1){1}}
\put(57,87){\line(1,0){1}}
\put(57,87){\line(0,-1){1}}
\put(141,87){\line(1,0){1}}
\put(141,87){\line(0,-1){1}}
\put(6,96){\line(1,0){1}}
\put(6,96){\line(0,-1){1}}
\put(45,99){\line(1,0){1}}
\put(45,99){\line(0,-1){1}}
\put(87,99){\line(1,0){1}}
\put(87,99){\line(0,-1){1}}
\put(129,99){\line(1,0){1}}
\put(129,99){\line(0,-1){1}}

\put(25.9,103.85){$\bullet$}
\put(67.9,103.85){$\bullet$}
\put(148.8,100.5){$\bullet$}

\put(104,-1.1){$\cdots$}
\put(104,19.8){$\cdots$}
\put(104,40.7){$\cdots$}
\put(104,61.85){$\cdots$}
\put(104,82.85){$\cdots$}
\put(104,103.85){$\cdots$}

\put(14,-4){$\p{n}$}
\put(55,-4){$\p{n-1}$}
\put(-4,22){$\p{1}$}
\put(-4,64){$\p{2}$}
\put(28,107){$\p{3}$}
\put(70,107){$\p{4}$}
\put(152,105){${k+3 \over 2}$}
\put(152,63){${k+5 \over 2}$}
\put(152,21){${k+7 \over 2}$}

\end{picture}

\bigskip
\bigskip
\bigskip
\bigskip
\bigskip
\bigskip
\bigskip
\bigskip
\bigskip
\bigskip
\bigskip
\bigskip
\bigskip
\caption{Quadrilateral Arrangement ${\bf A_{k,k+4}}$}
\label{quad arrangement 3}
\end{center}
\end{figure}
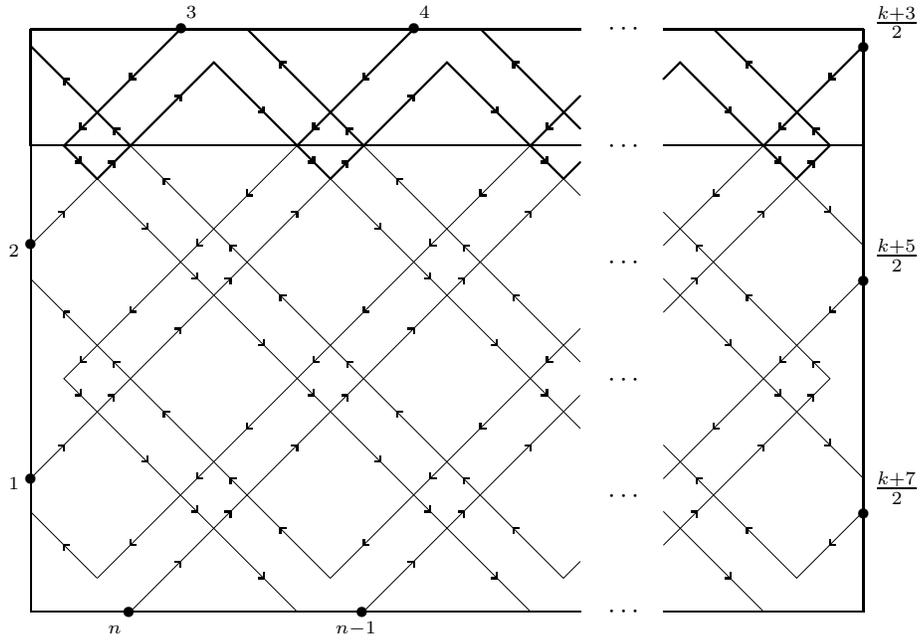

\newpage
\noindent
Notice that the paths in the 
upper band are attached to the paths in the
previously constructed diagram.
This is a Postnikov arrangement.
Most of the conditions specified
in the definition, such as the alternation condition,
are immediately verified.
As for the third condition, it is enough to check
that the paths corresponding to ``peaks'' in
the upper band do not violate the specified
crossing condition. Inspection reveals
that if one of these peaks failed to meet
the criterion then so would the neighboring
peaks underneath the band; which, by
induction, is impermissible.

\bigskip
\noindent
To check that the diagram yields the permutation $\pi_{k,n}$
it is sufficient to check that the number
of sources (black nodes) between the endpoint
of a path and its sources (measured clockwise
from the endpoint) is $n-k$ for those paths
which ``peak'' in the upper band. 
Paths of this type differ from the previous level
of paths (peaking in the original diagram)
in that their endpoints are shifted back one notch
(counter-clockwise) from their original positions -
so the number of black nodes in the interval
decreases by one. On the other hand, one extra
black node is introduced in the left hand wall
of the upper band which contributes to the interval
- thus increasing the number of black nodes back
by one. Induction on the lower portion of the
diagram predicts this number is $n-k$.

\bigskip
\noindent
Of the four paths which form the quadrilateral
perimeter of the internal rectangle situated in 
the $i$-th row and $j$-th column of ${\bf A_{k,n}}$, 
only two will contribute to its $k$-subset
label $I_{ij}$. If $p_{ij}$ and
$p_{ij}'$ denote the labels of these two paths then
$\big\{ p_{ij} \ p_{ij}' \big\}$ is given
by

\[ \Big\{ \ p_{ij} \ \ p_{ij}' \ \Big\} \ = \ \Bigg\{ \
\rho^{2 - \big\lceil {i \over 2} \big\rceil - 
\big\lceil {j \over 2} \big\rceil }
(n-k) \ \ \ \ \rho^{ \big\lfloor {i \over 2} \big\rfloor +
\big\lfloor {j \over 2} \big\rfloor} \ (n-k+2) \ \Bigg\} \] 

\[ \text{where} \]

\[ \rho \ = \ \begin{pmatrix}
1 & 2 & 3 & \cdots & n-2 & n-1 & n \\
2 & 3 & 4 & \cdots & n-1 & n & 1 \end{pmatrix} \]

\bigskip
\noindent
The chord $\big[ p_{ij} \ p_{ij}' \big]$ clearly
belongs to the ${\bf T_{k,n}}$ triangulation of ${\bf P}$ and
the corresponding $k$-subset label $I_{ij}$ 
is the disjoint union of intervals

\[ \Big[ p_{ij} \ \dots \ p_{ij}+(i-1) \ \Big] \ \sqcup  
\ \Big[p_{ij}' \ \dots \ p_{ij}' + (k-i-1) \ \Big]. \]

\end{proof}

\bigskip
\noindent
Note that two $k$-subsets of the form
$\big[i_1 \dots i_1'\big] \sqcup [i_2 \dots i_2' \big]$ and
$\big[j_1 \dots j_1'\big] \sqcup [j_2 \dots j_2' \big]$
are non-crossing if and only if the chords
$\big[i_1 \ i_2 \big]$ and $\big[j_1 \ j_2 \big]$
are non-crossing. The chords $\big[ p_{ij} \ p_{ij}' \big]$
are non-crossing because they belong to ${\bf T_{k,n}}$
and since the collection of
even cells of ${\bf A_{k,n}}$ has the requisite cardinality of 
$k(n-k)+1$, we obtain the following corollary:

\bigskip
\begin{Cor}\label{ch2,cor.quad}
The collection of labeling sets of ${\bf A_{k,n}}$ is
a maximal family of pairwise non-crossing $k$-subsets.
\end{Cor}

\bigskip
\noindent
When $n= 2k$ the diagram ${\bf A_{k,n}}$ replicates a
{\it double wiring arrangement} (cf. \cite{FZ1}) as 
introduced in section 2.4. For instance, when $k=3$,
this double wiring arrangement 
corresponds to the reduced word expression 
\ ${\bf 2 \ 1 \ 2 } \ 1 \ 2 \ 1 $ \ of
$(w_0 ,w_0) \in S_3 \times S_3$ 
and is depicted below.

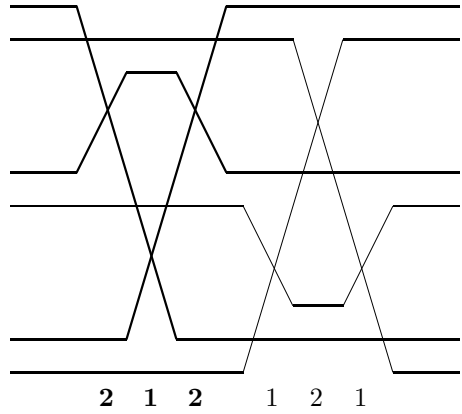
\begin{figure}[ht]

\setlength{\unitlength}{2.10pt}

\begin{center}
\begin{picture}(-15,30)(50,50)
\thinlines

\put(0,0){\line(1,0){42}}
\put(42,0){\line(18,60){18}}
\put(60,60){\line(1,0){21}}

\put(0,30){\line(1,0){42}}
\put(42,30){\line(9,-18){9}}
\put(51,12){\line(1,0){9}}
\put(60,12){\line(9,18){9}}
\put(69,30){\line(1,0){12}}

\put(0,60){\line(1,0){51}}
\put(51,60){\line(18,-60){18}}
\put(69,0){\line(1,0){12}}

\thicklines

\put(0,6){\line(1,0){21}}
\put(21,6){\line(18,60){18}}
\put(39,66){\line(1,0){42}}

\put(0,36){\line(1,0){12}}
\put(12,36){\line(9,18){9}}
\put(21,54){\line(1,0){9}}
\put(30,54){\line(9,-18){9}}
\put(39,36){\line(1,0){42}}

\put(0,66){\line(1,0){12}}
\put(12,66){\line(18,-60){18}}
\put(30,6){\line(1,0){51}}

\put(16,-6){${\bf 2}$}
\put(24,-6){${\bf 1}$}
\put(32,-6){${\bf 2}$}

\put(46,-6){$1$}
\put(54,-6){$2$}
\put(62,-6){$1$}

\end{picture}

\bigskip
\bigskip
\bigskip
\bigskip
\bigskip
\bigskip
\bigskip
\bigskip
\bigskip
\bigskip
\bigskip
\bigskip
\bigskip
\bigskip
\caption{Double Wiring Arrangement}
\label{double wiring}
\end{center}
\end{figure}

\newpage
\bigskip
\noindent
To recover the Postnikov diagram
simply add a boundary and assign
alternating orientations at the
endpoints of the
``wires'' of the double wiring arrangement
as follows:

\begin{figure}[ht]

\setlength{\unitlength}{2.10pt}

\begin{center}
\begin{picture}(-15,30)(50,50)

\thicklines
\put(0,-6){\line(0,1){78}}
\put(0,72){\line(1,0){81}}
\put(81,72){\line(0,-1){78}}
\put(0,-6){\line(1,0){81}}

\thinlines
\put(0,0){\line(1,0){42}}
\put(42,0){\line(18,60){18}}
\put(60,60){\line(1,0){21}}

\put(0,30){\line(1,0){42}}
\put(42,30){\line(9,-18){9}}
\put(51,12){\line(1,0){9}}
\put(60,12){\line(9,18){9}}
\put(69,30){\line(1,0){12}}

\put(0,60){\line(1,0){51}}
\put(51,60){\line(18,-60){18}}
\put(69,0){\line(1,0){12}}


\put(0,6){\line(1,0){21}}
\put(21,6){\line(18,60){18}}
\put(39,66){\line(1,0){42}}

\put(0,36){\line(1,0){12}}
\put(12,36){\line(9,18){9}}
\put(21,54){\line(1,0){9}}
\put(30,54){\line(9,-18){9}}
\put(39,36){\line(1,0){42}}

\put(0,66){\line(1,0){12}}
\put(12,66){\line(18,-60){18}}
\put(30,6){\line(1,0){51}}

\put(-1.15,4.75){$\bullet$}
\put(-1.15,34.75){$\bullet$}
\put(-1.15,64.75){$\bullet$}
\put(79.8,58.75){$\bullet$}
\put(79.8,28.75){$\bullet$}
\put(79.8,-1.25){$\bullet$}

\put(-4,4.5){$1$}
\put(-4,34.5){$2$}
\put(-4,64.5){$3$}
\put(83,58.5){$4$}
\put(83,28.5){$5$}
\put(83,-1.5){$6$}

\put(9,6){\line(-1,1){1}}
\put(9,6){\line(-1,-1){1}}
\put(39,6){\line(-1,1){1}}
\put(39,6){\line(-1,-1){1}}
\put(57,6){\line(-1,1){1}}
\put(57,6){\line(-1,-1){1}}
\put(75,6){\line(-1,1){1}}
\put(75,6){\line(-1,-1){1}}

\put(9,36){\line(-1,1){1}}
\put(9,36){\line(-1,-1){1}}
\put(27,54){\line(-1,1){1}}
\put(27,54){\line(-1,-1){1}}
\put(48,36){\line(-1,1){1}}
\put(48,36){\line(-1,-1){1}}
\put(75,36){\line(-1,1){1}}
\put(75,36){\line(-1,-1){1}}
\put(9,66){\line(-1,1){1}}
\put(9,66){\line(-1,-1){1}}
\put(75,66){\line(-1,1){1}}
\put(75,66){\line(-1,-1){1}}

\put(8,0){\line(1,1){1}}
\put(8,0){\line(1,-1){1}}
\put(74,0){\line(1,1){1}}
\put(74,0){\line(1,-1){1}}

\put(8,30){\line(1,1){1}}
\put(8,30){\line(1,-1){1}}
\put(38,30){\line(1,1){1}}
\put(38,30){\line(1,-1){1}}
\put(56,12){\line(1,1){1}}
\put(56,12){\line(1,-1){1}}
\put(74,30){\line(1,1){1}}
\put(74,30){\line(1,-1){1}}
\put(8,60){\line(1,1){1}}
\put(8,60){\line(1,-1){1}}
\put(26,60){\line(1,1){1}}
\put(26,60){\line(1,-1){1}}
\put(47,60){\line(1,1){1}}
\put(47,60){\line(1,-1){1}}
\put(74,60){\line(1,1){1}}
\put(74,60){\line(1,-1){1}}

\end{picture}

\bigskip
\bigskip
\bigskip
\bigskip
\bigskip
\bigskip
\bigskip
\bigskip
\bigskip
\bigskip
\bigskip
\bigskip
\bigskip
\bigskip
\caption{Corresponding Postnikov Arrangement}
\label{corresponding arrangement}
\end{center}
\end{figure}
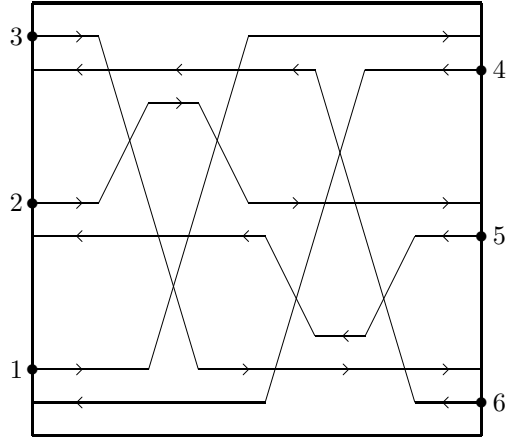

\bigskip
\noindent
For $k=4,5,6$, the reduced word expressions of $(w_0,w_0)
\in S_k \times S_k$ corresponding to the Postnikov
diagrams ${\bf A_{k,2k}}$ are

\[ {\bf 2 \ 3 \ 2 \ 1 \ 2 \ 3 } \ 1 \ 2 \ 3 \ 2 \ 1 \ 2 \] 
\[ {\bf 3 \ 2 \ 3 \ 4 \ 3 \ 2\ 1 \ 2 \ 3 \ 4 }
\ 1\ 2 \ 3 \ 4 \ 3 \ 2 \ 1 \ 2 \ 3 \ 2 \]
\[ {\bf 3 \ 4 \ 3 \ 2 \ 3 \ 4 \ 5 \ 4 \ 3 \ 2 \ 1 
\ 2 \ 3 \ 4 \ 5 } \ 1 \ 2 \ 3 \ 4 \ 5 \ 4
\ 3 \ 2 \ 1 \ 2 \ 3 \ 4 \ 3 \ 2 \ 3 .\]

\bigskip
\noindent
To obtain the reduced word $\mathcal{W}$ for general $k$, take
the single reduced word expression $\mathcal{R}$
for $w_0 \in S_k$ given by

\[ 1 \ \cdots \ k-1 \ \cdots \ 1 \ \cdots \ k-2 \ \cdots
\ 2 \ \cdots \ k-3 \ \cdots \ 3 \ \cdots \ k-4 \ \cdots
\ \Big\lfloor {k \over 2} \Big\rfloor.  \]

\bigskip
\noindent
Let $\mathcal{R}'$ be the reduced word obtained
by replacing each index $i$ occurring in $\mathcal{R}$
by $k-i$ and then flipping the entire expression.
The double reduced word $\mathcal{W}$ is ${\bf \mathcal{R}'} 
\ \mathcal{R}$.

\newpage
\bigskip
\section{Main Result}

\bigskip
\noindent
Two even regions of a Postnikov arrangement, 
with labels $I$ and $J$, are said to be
{\it neighbors} in a Postnikov arrangement if locally they
are situated as:

\begin{figure}[ht]

\setlength{\unitlength}{.75pt}

\begin{center}
\begin{picture}(270,60)(-80,100)
\thicklines

    \put(20,130){\line(1,-1){60}}
    \put(3,145){{\bf .}}
    \put(8,140){{\bf .}}
    \put(13,135){{\bf .}}
    \put(35,115){\line(0,1){3}}
    \put(35,115){\line(-1,0){3}}
    \put(65,115){\line(0,-1){3}}
    \put(65,115){\line(-1,0){3}}
    \put(65,85){\line(0,1){3}}
    \put(65,85){\line(-1,0){3}}
    \put(83,65){{\bf .}}
    \put(88,60){{\bf .}}
    \put(93,55){{\bf .}}

    \put(20,70){\line(1,1){60}}  
    \put(3,55){{\bf .}}
    \put(8,60){{\bf .}}
    \put(13,65){{\bf .}}
    \put(35,85){\line(-1,0){3}}
    \put(35,85){\line(0,-1){3}} 
    \put(83,135){{\bf .}}
    \put(88,140){{\bf .}}
    \put(93,145){{\bf .}}

    \put(15,95){$I$}
    \put(75,95){$J$}

\end{picture}
\end{center}
\end{figure}

\bigskip
\bigskip
\bigskip
\bigskip
\noindent
In this case $I$ is said to be
{\it oriented towards} $J$; dually $J$ said to be is {\it 
oriented away} from $I$. Given a $\pi_{k,n}$-diagram
${\bf A}$ let ${\bf \tilde{B}(A)}$ be the integer
matrix with rows indexed by the $k$-subset labels
of ${\bf A}$ and columns indexed by the corresponding
interior $k$-subsets labels and possessing the entries
$b_{\s I,J}$ defined by

\bigskip
\[ b_{\s I,J} \ = \
\left\{ \begin{array}{ll}
  \ \ 1 & \quad \text{if $I$ is oriented towards $J$} \\
  -1 & \quad \text{if $I$ is oriented away from $J$} \\
  \ \ 0 & \quad \text{otherwise.} 
\end{array} \right. \]

\bigskip
\noindent
The principal submatrix ${\bf B(A)}$ is clearly skew-symmetric.
Consider now the field $\mathcal{F}$
of rational functions generated by indeterminates $[K]$ for
$k$-subsets labels $K$ arising in the
arrangement ${\bf A_{k,n}}$. Set ${\bf x \big(A_{k,n}\big)}$ to be the
set of indeterminates corresponding to interior 
$k$-subset labels in ${\bf A_{k,n}}$
and let ${\bf c}$ be the set of indeterminates $[K_1],\dots
[K_n]$ corresponding
to the remaining $k$-subset labels which are boundaries.
Let $\mathcal{A}_{k,n}$ denote the cluster algebra
generated by the seed ${ \bf \Big( x \big( A_{k,n} \big),
c, \tilde{B} \big( A_{k,n} \big) \Big) }$ inside $\mathcal{F}$.

\bigskip
\begin{Thm}\label{ch2,thm.clusters}
Each $\pi_{k,n}$-diagram ${\bf A}$ gives rise to
a seed ${\bf \Big( x(A), c , \tilde{B}(A) \Big) }$
in $\mathcal{A}_{k,n}$ - whose cluster variables
are indexed by the interior $k$-subsets of ${\bf A}$ -
with the property that if ${\bf A'}$
is obtained from \ ${\bf A}$ by a single
geometric exchange through a quadrilateral cell
labeled $K$ in ${\bf A}$, 
then \ $\mu_{K} \Big( {\bf \tilde{B}\big( A \big)}
\Big) = {\bf \tilde{B} \big( A' \big) }$.
\end{Thm}

\bigskip
\noindent
Recall that a cluster algebra enjoys an
algebra structure over the ring of polynomials 
in indeterminates $c \in {\bf c}$. In the
present case ${ \bf c} = \big\{ 
[K_1], \dots, [K_n] \big\}$ where
$K_1, \dots, K_n$ are the boundary $k$-subsets
of $[1 \dots n]$. The homogeneous
coordinate ring $\Bbb{C} \Big[ \Bbb{G}(k,n) \Big]$
of the Grassmannian 
can also be consider an algebra over
the polynomial ring $\Bbb{C} \Big[[K_1],\dots,[K_n]\Big]$;
simply identify the scalars $[K_1], \dots,
[K_n]$ with
the Pl\"ucker coordinates $\Delta^{K_1}, \dots, 
\Delta^{K_n}$.

\bigskip 
\begin{Thm}\label{ch2,thm.main}
There is a isomorphism 
$\varphi: \mathcal{A}_{k,n}
\longrightarrow \Bbb{C} \Big[ \Bbb{G}(k,n) \Big]$
of \ $\Bbb{C} \Big[[K_1],\dots,[K_n]\Big]$- algebras
with the property that $[K] \stackrel{\varphi}{\mapsto} \Delta^K$ 
for every $k$-subsets $K$ of $[1 \dots n]$.
\end{Thm}

\bigskip
\section*{ {\it Proof of Theorem \ref{ch2,thm.main} :}}

\noindent
This argument proceeds by establishing the conditions
neccessary to implement Proposition \ref{ch2,CA2.1}
for the cluster algebra $\mathcal{A}_{k,n}$
and the affine cone $X(k,n)$ of the Grassmannian.
We will use Theorem \ref{ch2,thm.clusters} and
return to its proof later.

\bigskip
\noindent
Let $\mathcal{X}$ denote the set of cluster variables
of $\mathcal{A}_{k,n}$.
By Proposition \ref{ch2,laurent} every cluster 
variable $x \in \mathcal{X}$ can be expressed as 
a unique Laurent polynomial 

\[ { 
P \ \Big( [K] \ ; 
\ [K] \in {\bf x \big(A_{k,n}\big) \sqcup c} \ \Big)
\over 
{ m \Big( [K] \ ; \ [K] \in {\bf x \big(A_{k,n} \big) \Big)} } 
} \]

\bigskip
\noindent
where $P$ is a polynomial and $m$ is monomial. 
By this expression, the element $x$ determines
a rational function 

\[ \varphi_x := { P \ \Big( \Delta^K \ ; \ [K] 
\in {\bf x \big(A_{k,n}\big) \sqcup c} \ \Big)
\over { m \Big( \Delta^K \ ; \ [K] \in {\bf x \big(A_{k,n} \big) \Big)} 
} } \]

\bigskip
\noindent
on the Grassmannian $\Bbb{G}(k,n)$. According
to Theorem \ref{ch2,thm.clusters}, every $k$-subset 
$K$ occurring as a label in some $\pi_{k,n}$-diagram
${\bf A}$ corresponds to a cluster variable 
$[K] \in \mathcal{A}_{k,n}$. 

\bigskip
\begin{Claim}\label{ch2,claim.plucker}
$\varphi_{\s [K] } = \Delta^K$ for $k$-subsets
$K$ occurring as labels in $\pi_{k,n}$-diagrams.
\end{Claim}

\begin{proof}
By construction, $\varphi_{\s [K]} = \Delta^K$
for all $k$-subsets $K$ in ${\bf A_{k,n}}$.
According to Proposition \ref{ch2,Post.1}
all $\pi_{k,n}$-diagrams ${\bf A}$ can be 
obtained from ${\bf A_{k,n}}$ by a sequence
of geometric exchanges. Therefore, in
order to establish the claim it is enough to
prove the following: whenever the claim holds for  
the $k$-subsets of 
a particular $\pi_{k,n}$-diagram ${\bf A}$ 
then the claim must hold for $k$-subsets
of $\pi_{k,n}$-diagrams ${\bf A'}$ obtained from
${\bf A}$ by a single geometric exchange.

\bigskip
\noindent
With this aim, take an interior $k$-subset
$K$ of a $\pi_{k,n}$-diagram ${\bf A}$ -
a diagram for which the claim is valid. 
Assume that $K$ labels a quadrilateral
cell in ${\bf A}$. In this case there exists
a subset $I \subset [1 \dots n]$ of size $|I| = k-2$
and four distinct indices $i,j,s,t$ disjoint from $I$
such that $K = Iij$ and the following local configuration
is present in the diagram ${\bf A_{k,n}}$:

\newpage
\begin{figure}[ht]

\setlength{\unitlength}{0.65pt}

\begin{center}
\begin{picture}(200,170)(50,50)

\thicklines

  \put(100,0){\line(1,1){150}}
  \put(125,25){\line(1,0){5}}
  \put(125,25){\line(0,1){5}}
  \put(175,75){\line(1,0){5}}
  \put(175,75){\line(0,1){5}}
  \put(225,125){\line(1,0){5}}
  \put(225,125){\line(0,1){5}}

  \put(50,50){\line(1,1){150}}
  \put(75,75){\line(0,-1){5}}
  \put(75,75){\line(-1,0){5}}
  \put(125,125){\line(0,-1){5}}
  \put(125,125){\line(-1,0){5}}
  \put(175,175){\line(0,-1){5}}
  \put(175,175){\line(-1,0){5}}

  \put(100,200){\line(1,-1){150}}
  \put(125,175){\line(0,-1){5}}
  \put(125,175){\line(1,0){5}}
  \put(175,125){\line(0,-1){5}}
  \put(175,125){\line(1,0){5}}
  \put(225,75){\line(0,-1){5}}
  \put(225,75){\line(1,0){5}}

  \put(50,150){\line(1,-1){150}}
  \put(75,125){\line(0,1){5}}
  \put(70,125){\line(1,0){5}}
  \put(125,75){\line(0,1){5}}
  \put(120,75){\line(1,0){5}}  
  \put(175,25){\line(0,1){5}}
  \put(170,25){\line(1,0){5}}

  \put(45,150){$i$}
  \put(250,45){$j$}
  \put(253,148){$s$}
  \put(43,48){$t$}

  \put(143,100){$Iij$}
  \put(143,175){$Iit$}
  \put(143,20){$Ijs$}
  \put(225,100){$Iis$}
  \put(60,100){$Ijt$}

\end{picture}
\end{center}
\end{figure}

\bigskip
\bigskip
\bigskip
\bigskip
\noindent
Let ${\bf A'}$ be the $\pi_{k,n}$-diagram
obtained by executing the geometric
exchange about the quadrilateral cell
labeled by $Iij$. 
Theorem \ref{ch2,thm.clusters} predicts
that the cluster variables $[Iij] \in
{\bf x \big( A \big)}$ and
$[Ist] \in {\bf x \big( A' \big)}$ are related by the following
exchange relation

\[ [Ist] \ [Iij] \ = \ [Iit][Ijs] \ + \ [Ijt][Iis]. \]
 
\bigskip
\noindent
By Corollary \ref{ch2,cor.max} of section 3, the
five $k$-subset $Iij$, $Iit$, $Ijs$, $Ijt$, and $Iis$ 
are pairwise non-crossing and thus this exchange
relation corresponds to the short Pl\"ucker relation:

\[ \Delta^{Ist} \Delta^{Iij} \ = \ \Delta^{Iit} \Delta^{Ijs}
\ + \ \Delta^{Ijt} \Delta^{Iis} \] 

\bigskip
\noindent
It follows immediately that $\varphi_{\s [Ist]}$ is the Pl\"ucker
coordinate $\Delta^{Ist}$. Aside of $Iij$ and $Ist$, 
the remaining $k$-subsets of ${\bf A'}$ are the same as 
${\bf A}$'s; and thus the claim holds for all
$k$-subsets of ${\bf A'}$.

\end{proof}

\bigskip
\noindent
Proposition \ref{ch2,Post.1}, section 3, stipulates that
every $k$-subset $K \in [1 \dots n]$
occurs as a labeling set in some $\pi_{k,n}$-diagram
${\bf A}$. From this we may conclude that

\[ [K] \stackrel{\varphi}{\mapsto} \Delta^K \]

\bigskip
\noindent
is valid for all $k$-subsets in $[1 \dots n]$;
in particular $\varphi_{\s [K]}$ is regular
on $X(k,n)$. Proposition \ref{ch2,CA2.1} requires
that $\varphi_x$ must be regular for each cluster
variable $x$.   

\bigskip
\begin{Prop}\label{ch2,prop.reg} For each $x \in \mathcal{X}$ the 
corresponding rational function $\varphi_x$ is
in fact regular on the affine cone $X(k,n)$
of the Grassmannian $\Bbb{G}(k,n)$.
\end{Prop}

\begin{proof}
Every interior $k$-subset of ${\bf A_{k,n}}$
is the label of a quadrilateral cell. As indicated
in the figure above,
such a $k$-subset label is of the form $Iij$ where
$I \subset [1 \dots n]$ is a subset of size $k-2$
and $i,j,s,t$ 
are four distinct indices disjoint from $I$.
Let ${ \bf A_{k,n}^{I {\s ij} } }$ be the
$\pi_{k,n}$-diagram obtained from ${\bf A_{k,n}}$
by the geometric exchange replacing $Iij$ with
$Ist$. 

\bigskip
\noindent
Let $U$ be the open subvariety of
$X(k,n)$ consisting of all points for which $\Delta^K
\ne 0$ for each interior $k$-subset label $K$ in
${\bf A_{k,n}}$. For an interior $k$-subset $Iij$  
let $U^{I {\s ij}}$ be the open subvariety of 
$X(k,n)$ consisting of points for which $\Delta^K
\ne 0$ for all interior $k$-subsets in ${\bf A_{k,n}}$ 
with $K \ne Iij$ \ and such that $\Delta^{Ist} \ne 0$. Finally
let $\mathcal{U}$ denote the union of $U$ and
all $U^{I {\s ij}}$. The coordinates $\Delta^{Iij}$ and $\Delta^{Ist}$
are prime to one another and so 
the codimension of the intersection of their
zero-loci is 2 or greater. Consequently the 
codimention of the complement of 
$\mathcal{U}$ is $2$ or greater.

\bigskip
\noindent
Take $x \in \mathcal{X}$. If we can establish that
$\varphi_x$ is regular on $\mathcal{U}$ then it
will follow that $\varphi_x$ is regular on all
$X(k,n)$; this is in virtue of the codimension 
of $\mathcal{U}$'s complement.
The Laurent Phenomena detailed in Proposition \ref{ch2,laurent}
of section 2 guarantees that $x$ can be 
simultaneously expressed as 
a Laurent polynomial in the cluster variables from
${\bf x\big(A_{k,n}\big)}$ and as a Laurent
polynomial in the cluster variables of ${\bf x
\big( A_{k,n}^{I {\s ij}} \big)}$ for all
interior $k$-subsets $Iij$ \ of ${\bf A_{k,n}}$.
From this we may conclude that the
denominator of $\varphi_x$ - with
respect to the ${\bf x\big(A_{k,n}\big)}$ cluster
- is a monomial in the Pl\"ucker
coordinates $\Delta^K$ for interior
$k$-subsets $K$ of ${\bf A_{k,n}}$. This implies
that $\varphi_x$ is regular on $U$.
On the other hand, the denominator of 
$\varphi_x$ - with respect to the  
${\bf x \big( A_{k,n}^{I {\s ij}} \big)}$ cluster -
is a monomial in the Pl\"ucker coordinates
$\Delta^{I {\s s t}}$ and 
$\Delta^K$ for $k$-subsets $K \ne Iij$ in
${\bf A_{k,n}}$; consequently we may surmise that 
$\varphi_x$ is regular on each $U^{I {\s ij}}$
and thus regular on all of $\mathcal{U}$.
\end{proof}

\bigskip
\noindent
In view of Proposition \ref{ch2,Post.1}, section 3, it follows that
$\dim \Big( X(k,n) \Big) =  \rank \Big( \mathcal{A}_{k,n} \Big) + 
| {\bf c } |$. Any exchange relation shared among the cluster
variables $x \in \mathcal{X}$ remains valid for
the elements $\varphi_x$; this is because each
$\varphi_x$ is constructed in 
$\Bbb{C}\Big[ \Bbb{G}(k,n) \Big]$ by means of
a Laurent expansion with respect to the initial
cluster ${\bf x\big( A_{k,n} \big)}$.
Finally
Proposition \ref{ch2,CA2.1} requires  
of $x \in \mathcal{X}$ and $c \in {\bf
\frak{c}}$ that the 
corresponding regular functions $\varphi_x$ and $\varphi_c$ generate
$\Bbb{C}\Big[ \Bbb{G}(k,n) \Big]$;
this is true because the Pl\"ucker coordinates
generate the coordinate ring of the Grassmannian.

\bigskip
\section*{ {\it Proof of Theorem \ref{ch2,thm.clusters} :}}

\noindent
By Proposition \ref{ch2,Post.2} any
$\pi_{k,n}$-diagram ${\bf A}$ can
be obtained from ${\bf A_{k,n}}$ by
a sequence of geometric exchanges.
Consequently it is enough
to show that if ${\bf \Big( x(A),c, \tilde{B}(A) \Big) }$
is a seed and $K$ is a labeling set
of a quadrilateral cell in ${\bf A}$ then  

\[ \mu_{K} \Big({\bf  \tilde{B}(A) } \Big) \ = 
\ {\bf \tilde{B}(A')}. \]

\noindent
where ${\bf A'}$ is the $\pi_{k,n}$-diagram
obtained by executing the geometric
exchange associated to $K$ in ${\bf A}$.
The pertinent identity to be checked is

\[ b_{X,Y}' \ = \
b_{X,Y} +
{ \Big| \ b_{X,K} \Big|
  \ b_{K,Y} \ + \
  b_{X,K} \ \Big| \
  b_{K,Y} \Big| \over 2 } \]

\bigskip
\noindent
where $X$ and $Y$ are the labeling sets of any two of the 
four neighboring cells adjacent to $K$. 
Up to dihedral symmetry
there are two cases to consider.
Let $\{i,j,s,t,u,v\}$ be six distinct
indices disjoint from an index set $I \subset 
[1 \dots n]$ of cardinality $|I|=k-3$ 
(in case 1) and $|I|=k-4$ (in case 2).

\newpage
\bigskip
\noindent
{\bf Case 1: $K=Iijv$ and $X=Iisv$}

\begin{figure}[ht]

\setlength{\unitlength}{1.0pt}

\begin{center}
\begin{picture}(150,200)(100,50)

\thicklines

  \put(100,0){\line(1,1){200}}
  \put(125,25){\line(1,0){5}}
  \put(125,25){\line(0,1){5}}
  \put(175,75){\line(1,0){5}}
  \put(175,75){\line(0,1){5}}
  \put(238,138){\line(1,0){5}}
  \put(238,138){\line(0,1){5}}
  \put(287,187){\line(1,0){5}}
  \put(287,187){\line(0,1){5}}

  \put(50,50){\line(1,1){150}}
  \put(75,75){\line(0,-1){5}}
  \put(75,75){\line(-1,0){5}}
  \put(125,125){\line(0,-1){5}}
  \put(125,125){\line(-1,0){5}}
  \put(175,175){\line(0,-1){5}}
  \put(175,175){\line(-1,0){5}}

  \put(100,200){\line(1,-1){200}}
  \put(125,175){\line(0,-1){5}}
  \put(125,175){\line(1,0){5}}
  \put(175,125){\line(0,-1){5}}
  \put(175,125){\line(1,0){5}}
  \put(238,62){\line(0,-1){5}}
  \put(238,62){\line(1,0){5}}
  \put(287,13){\line(0,-1){5}}
  \put(287,13){\line(1,0){5}}

  \put(50,150){\line(1,-1){150}}
  \put(75,125){\line(0,1){5}}
  \put(70,125){\line(1,0){5}}
  \put(125,75){\line(0,1){5}}
  \put(120,75){\line(1,0){5}}  
  \put(175,25){\line(0,1){5}}
  \put(170,25){\line(1,0){5}}
  
  \put(250,0){\line(1,1){50}}
  \put(287,37){\line(-1,0){5}}
  \put(287,37){\line(0,-1){5}} 
  \put(263,13){\line(-1,0){5}}
  \put(263,13){\line(0,-1){5}}

  \put(250,200){\line(1,-1){50}}
  \put(287,163){\line(-1,0){5}}
  \put(287,163){\line(0,1){5}}
  \put(263,187){\line(-1,0){5}}
  \put(263,187){\line(0,1){5}}  

  \put(45,150){$i$}
  \put(43,48){$t$}
  \put(302,202){$s$}
  \put(302,-5){$j$}
  \put(244,202){$u$}
  \put(244,-5){$v$}

  \put(141,100){$Iijv$}
  \put(141,175){$Iitv$}
  \put(141,20){$Ijsv$}
  \put(223,100){$Iisv$}
  \put(62,100){$Ijtv$}
  \put(264,202){$Iiuv$}
  \put(264,-5){$Iijs$}

\end{picture}
\end{center}
\end{figure}

\bigskip
\bigskip
\bigskip
\bigskip
\bigskip
\bigskip
\bigskip
\noindent
The submatrix
of ${\bf \tilde{B}(A)}$ whose rows and
columns are indexed by the (ordered)
collection of $k$-subsets $[Iisv],[Iijv],
[Iitv],[Ijsv],[Ijtv]$ is
given by:

\[ \begin{pmatrix} 
\ \ 0 & \ \ 1 & \ \ 0 & \ \ 0 & \ \ 0 \\
-1 & \ \ 0 & \ \ 1 & \ \ 1 & -1 \\ 
\ \ 0 & -1 & & & \\
\ \ 0 & -1 & & \ \ * & \\ 
\ \ 0 & \ \ 1 & & & \end{pmatrix}. \]

\bigskip
\noindent
The same submatrix in $\mu_{I {\s ijv} } \Big(
{\bf \tilde{B}(A)} \Big)$ is

\[ \begin{pmatrix} 
\ \ 0 & -1 & \ \ 1 & \ \ 1 & \ \ 0 \\ 
\ \ 1 & \ \ 0 & -1 & -1 & \ \ 1 \\ 
-1 & \ \ 1 & & & \\ 
-1 & \ \ 1 & & \ \ * & \\ 
\ \ 0 & -1 & & &  \end{pmatrix}. \]

\bigskip
\noindent
This is identical to the submatrix 
of ${\bf \tilde{B}(A')}$ obtained 
after applying the geometric exchange about
the $Iijv$ cell. The local configuration
of ${\bf A'}$ is depicted below.

\newpage
\begin{figure}[ht]

\setlength{\unitlength}{1.0pt}

\begin{center}
\begin{picture}(200,60)(20,200)
\thicklines

  \put(-30,0){\line(1,1){210}}
  \put(0,30){\line(-1,0){3}}
  \put(0,30){\line(0,-1){3}}
  \put(60,90){\line(-1,0){3}}
  \put(60,90){\line(0,-1){3}}
  \put(120,150){\line(-1,0){3}}
  \put(120,150){\line(0,-1){3}}
  \put(165,195){\line(-1,0){3}}
  \put(165,195){\line(0,-1){3}}

  \put(-30,240){\line(1,-1){210}}
  \put(0,210){\line(0,1){3}}
  \put(0,210){\line(-1,0){3}}
  \put(60,150){\line(0,1){3}}
  \put(60,150){\line(-1,0){3}}
  \put(120,90){\line(0,1){3}}
  \put(120,90){\line(-1,0){3}}
  \put(165,45){\line(0,1){3}}
  \put(165,45){\line(-1,0){3}}
   
  \put(210,240){\line(1,-1){90}}
  \put(225,225){\line(0,1){3}}
  \put(225,225){\line(-1,0){3}}
  \put(270,180){\line(0,1){3}}
  \put(270,180){\line(-1,0){3}}
   
  \put(210,0){\line(1,1){90}} 
  \put(225,15){\line(-1,0){3}}
  \put(225,15){\line(0,-1){3}} 
  \put(270,60){\line(-1,0){3}}
  \put(270,60){\line(0,-1){3}}

  \put(-60,90){\line(1,1){150}}
  \put(90,240){\line(1,-1){90}}
  \put(180,150){\line(1,1){90}}
  \put(0,150){\line(0,1){3}}
  \put(0,150){\line(1,0){3}}
  \put(60,210){\line(0,1){3}}
  \put(60,210){\line(1,0){3}}
  \put(210,180){\line(0,1){3}}
  \put(210,180){\line(1,0){3}}
  \put(255,225){\line(0,1){3}}
  \put(255,225){\line(1,0){3}}
  \put(120,210){\line(0,-1){3}}
  \put(120,210){\line(1,0){3}}
  \put(165,165){\line(0,-1){3}}
  \put(165,165){\line(1,0){3}}

  \put(-60,150){\line(1,-1){150}}
  \put(90,0){\line(1,1){90}}
  \put(180,90){\line(1,-1){90}}  
  \put(-45,105){\line(0,1){3}}
  \put(-45,105){\line(1,0){3}}  
  \put(-45,135){\line(1,0){3}}
  \put(-45,135){\line(0,-1){3}}
  \put(0,90){\line(0,-1){3}}
  \put(0,90){\line(1,0){3}}
  \put(60,30){\line(0,-1){3}}
  \put(60,30){\line(1,0){3}} 
  \put(210,60){\line(0,-1){3}}
  \put(210,60){\line(1,0){3}}
  \put(255,15){\line(0,-1){3}}
  \put(255,15){\line(1,0){3}}
  \put(120,30){\line(0,1){3}}
  \put(120,30){\line(1,0){3}}
  \put(165,75){\line(0,1){3}}
  \put(165,75){\line(1,0){3}}

  \put(-35,242){$i$}
  \put(-35,-5){$t$}
  \put(204,-5){$v$}
  \put(204,242){$u$}
  \put(272,-5){$j$}
  \put(272,242){$s$}

  \put(20,117){$Istv$}
  \put(20,225){$Iitv$}  
  \put(20,18){$Ijsv$}
  \put(-60,117){$Ijtv$}
  \put(230,242){$Iiuv$}
  \put(230,-5){$Iijs$}
  \put(172,117){$Iisv$}

\end{picture}
\end{center}
\end{figure}

\newpage
\bigskip
\noindent
{\bf Case 2: $K=Iijuv$ and $X=Iisuv$}

\begin{figure}[ht]

\setlength{\unitlength}{1.0pt}

\begin{center}
\begin{picture}(200,60)(20,200)
\thicklines

  \put(-30,0){\line(1,1){210}}
  \put(0,30){\line(1,0){3}}
  \put(0,30){\line(0,1){3}}
  \put(60,90){\line(1,0){3}}
  \put(60,90){\line(0,1){3}}
  \put(120,150){\line(1,0){3}}
  \put(120,150){\line(0,1){3}}
  \put(165,195){\line(1,0){3}}
  \put(165,195){\line(0,1){3}}

  \put(-30,240){\line(1,-1){210}}
  \put(0,210){\line(0,-1){3}}
  \put(0,210){\line(1,0){3}}
  \put(60,150){\line(0,-1){3}}
  \put(60,150){\line(1,0){3}}
  \put(120,90){\line(0,-1){3}}
  \put(120,90){\line(1,0){3}}
  \put(165,45){\line(0,-1){3}}
  \put(165,45){\line(1,0){3}}
   
  \put(210,240){\line(1,-1){90}}
  \put(225,225){\line(0,-1){3}}
  \put(225,225){\line(1,0){3}}
  \put(270,180){\line(0,-1){3}}
  \put(270,180){\line(1,0){3}}

  \put(120,30){\line(1,-1){30}}
  \put(120,30){\line(1,1){60}} 
  \put(180,90){\line(1,-1){90}}
  \put(135,45){\line(-1,0){3}}
  \put(135,45){\line(0,-1){3}}
  \put(135,15){\line(1,0){3}}
  \put(135,15){\line(0,-1){3}}
  \put(165,75){\line(-1,0){3}}
  \put(165,75){\line(0,-1){3}}
  \put(225,45){\line(0,1){3}}
  \put(225,45){\line(-1,0){3}}

  \put(-60,90){\line(1,1){150}}
  \put(90,240){\line(1,-1){90}}
  \put(180,150){\line(1,1){90}}
  \put(0,150){\line(0,-1){3}}
  \put(0,150){\line(-1,0){3}}
  \put(60,210){\line(0,-1){3}}
  \put(60,210){\line(-1,0){3}}
  \put(210,180){\line(0,-1){3}}
  \put(210,180){\line(-1,0){3}}
  \put(255,225){\line(0,-1){3}}
  \put(255,225){\line(-1,0){3}}
  \put(120,210){\line(0,1){3}}
  \put(120,210){\line(-1,0){3}}
  \put(165,165){\line(0,1){3}}
  \put(165,165){\line(-1,0){3}}

  \put(-60,150){\line(1,-1){120}}
  \put(30,0){\line(1,1){30}}
  \put(45,15){\line(1,0){3}}
  \put(45,15){\line(0,1){3}}
  \put(-45,105){\line(0,-1){3}}
  \put(-45,105){\line(-1,0){3}}  
  \put(-45,135){\line(-1,0){3}}
  \put(-45,135){\line(0,1){3}}
  \put(0,90){\line(0,1){3}}
  \put(0,90){\line(-1,0){3}}
  \put(45,45){\line(0,1){3}}
  \put(45,45){\line(-1,0){3}}

  \put(18,117){$Iijuv$}
  \put(18,205){$Iituv$}  
  \put(18,28){$Ijsuv$}
  \put(-62,117){$Ijtuv$}
  \put(228,235){$Iistv$}
  \put(138,28){$Iijst$}
  \put(138,117){$Iisuv$}

  \put(-65,150){$i$}
  \put(-65,85){$t$}
  \put(151,-5){$v$}
  \put(181,25){$j$}
  \put(181,210){$s$}
  \put(301,145){$u$}
\end{picture}
\end{center}
\end{figure}
   
\bigskip
\bigskip
\bigskip
\bigskip
\bigskip
\bigskip
\bigskip
\bigskip
\bigskip
\bigskip  
\bigskip
\bigskip
\bigskip
\bigskip
\bigskip
\bigskip
\bigskip
\bigskip
\bigskip
\bigskip
\bigskip
\bigskip
\bigskip
\bigskip
\bigskip
\noindent
The submatrix
of ${\bf \tilde{B}(A)}$ whose rows and
columns are indexed by the (ordered)
collection of $k$-subsets $[Iisuv],[Iijuv],
[Iituv],[Ijsuv],[Ijtuv]$ is
given by:

\[ \begin{pmatrix} 
\ \ 0 & \ \ 1 & -1 & \ \ 0& \ \ 0 \\
-1 & \ \ 0 & \ \ 1 & \ \ 1 & -1 \\ 
\ \ 1 & -1 & & & \\
\ \ 0 & -1 & & \ \ * & \\ 
\ \ 0 & \ \ 1 & & & \end{pmatrix}. \]

\bigskip
\noindent
The same submatrix in $\mu_{I {\s ijuv}} \Big(
{\bf \tilde{B}(A)} \Big)$ is

\[ \begin{pmatrix} 
\ \ 0 & -1 & \ \ 0 & \ \ 1 & \ \ 0 \\ 
\ \ 1 & \ \ 0 & -1 & -1 & \ \ 1 \\ 
\ \ 0 & \ \ 1 & & & \\ 
-1 & \ \ 1 & & \ \ * & \\ 
\ \ 0 & -1 & & &  \end{pmatrix}. \]

\newpage
\bigskip
\noindent
This is identical to the submatrix 
of ${\bf \tilde{B}(A')}$ obtained 
after applying the geometric exchange about
the $Iijuv$ cell. The local configuration
of ${\bf A'}$ is depicted below.

\begin{figure}[ht]

\setlength{\unitlength}{1.0pt}

\begin{center}
\begin{picture}(200,150)(20,100)
\thicklines

  \put(-30,0){\line(1,1){210}}
  \put(0,30){\line(-1,0){3}}
  \put(0,30){\line(0,-1){3}}
  \put(60,90){\line(-1,0){3}}
  \put(60,90){\line(0,-1){3}}
  \put(120,150){\line(-1,0){3}}
  \put(120,150){\line(0,-1){3}}
  \put(165,195){\line(-1,0){3}}
  \put(165,195){\line(0,-1){3}}

  \put(-30,240){\line(1,-1){210}}
  \put(0,210){\line(0,1){3}}
  \put(0,210){\line(-1,0){3}}
  \put(60,150){\line(0,1){3}}
  \put(60,150){\line(-1,0){3}}
  \put(120,90){\line(0,1){3}}
  \put(120,90){\line(-1,0){3}}
  \put(165,45){\line(0,1){3}}
  \put(165,45){\line(-1,0){3}}

  \put(-60,90){\line(1,1){120}}
  \put(45,195){\line(1,0){3}}
  \put(45,195){\line(0,1){3}}
  \put(45,225){\line(-1,0){3}}
  \put(45,225){\line(0,1){3}}
  \put(30,240){\line(1,-1){30}}
  \put(120,210){\line(1,-1){60}}
  \put(120,210){\line(1,1){30}}
  \put(135,225){\line(-1,0){3}}
  \put(135,225){\line(0,-1){3}}
  \put(135,195){\line(0,-1){3}}
  \put(135,195){\line(1,0){3}}

  \put(180,150){\line(1,1){90}}
  \put(0,150){\line(0,1){3}}
  \put(0,150){\line(1,0){3}}
  \put(210,180){\line(0,1){3}}
  \put(210,180){\line(1,0){3}}
  \put(255,225){\line(0,1){3}}
  \put(255,225){\line(1,0){3}}  
  \put(165,165){\line(0,-1){3}}
  \put(165,165){\line(1,0){3}}

  \put(210,0){\line(1,1){90}}
  \put(225,15){\line(0,-1){3}}
  \put(225,15){\line(-1,0){3}}
  \put(270,60){\line(0,-1){3}}
  \put(270,60){\line(-1,0){3}}   

  \put(-60,150){\line(1,-1){150}}
  \put(90,0){\line(1,1){90}}
  \put(180,90){\line(1,-1){90}}  
  \put(-45,105){\line(0,1){3}}
  \put(-45,105){\line(1,0){3}}  
  \put(-45,135){\line(1,0){3}}
  \put(-45,135){\line(0,-1){3}}
  \put(0,90){\line(0,-1){3}}
  \put(0,90){\line(1,0){3}}
  \put(60,30){\line(0,-1){3}}
  \put(60,30){\line(1,0){3}} 
  \put(210,60){\line(0,-1){3}}
  \put(210,60){\line(1,0){3}}
  \put(255,15){\line(0,-1){3}}
  \put(255,15){\line(1,0){3}}
  \put(120,30){\line(0,1){3}}
  \put(120,30){\line(1,0){3}}
  \put(165,75){\line(0,1){3}}
  \put(165,75){\line(1,0){3}}

  \put(-35,242){$i$}
  \put(-35,-5){$t$}
  \put(204,-5){$v$}
  \put(272,-5){$j$}
  \put(272,242){$u$}
  \put(24,242){$s$}

  \put(18,117){$Istuv$}
  \put(18,208){$Iituv$}  
  \put(18,24){$Ijsuv$}
  \put(-62,117){$Ijtuv$}
  \put(138,208){$Iistv$}
  \put(228,-5){$Iijst$}
  \put(170,117){$Iisuv$}

\end{picture}
\end{center}
\end{figure}

\newpage
\bigskip
\section{Toric Charts and Positivity}

\bigskip
\noindent
An important corollary 
of the cluster algebra structure
explained in Theorem \ref{ch2,thm.main} is that
each cluster gives
rise to a system of coordinates -
known as a {\it toric chart} - parameterizing
an algebraic torus $\big( \Bbb{C}^* \big)^{k(n-k)+1}$
of maximal dimension within the affine cone $X(k,n)$
of the Grassmannian.  
 
\bigskip
\noindent
To see this, set $N = (k-1)(n-k-1)$ and let ${\bf x} = \big\{
x_1, \dots, x_{\s N} \big\}$ be any cluster of the Grassmannian
$\Bbb{G}(k,n)$. Let $U$ be the open subvariety of 
$X(k,n)$ obtained by omitting the zero loci of
the cluster variables $x_1, \dots, x_{\s N}$
and the Pl\"ucker coordinates $\Delta^K$ whose
index $K$ is a boundary $k$-subset. 
Define a regular map 

\[ \Phi_{\bf x}: U \longrightarrow 
\Big( \Bbb{C}^* \Big)^{\s{N}} \times
\Big( \Bbb{C}^* \Big)^n \ \ \ \text{by} \]

\[ \Phi_{\bf x}(p) \ = \ \Big( \ x_1(p), \dots,
x_{\s N}(p) \ ; \ \Delta^{K_1}(p),
\dots, \Delta^{K_n}(p) \  \Big). \]

\bigskip
\noindent
For $1 \leq i \leq n$, the boundary $k$-subset
$K_i$ is the interval $[i \dots i+k]$ 
taken {\bf mod $n$}.   

\bigskip
\begin{Thm}\label{ch2,thm.toric} For any cluster ${\bf x}$ the
map $\Phi_{\bf x}$ is biregular.
\end{Thm}

\begin{proof}
By Proposition \ref{ch2,laurent}, any 
cluster variable $x$ in $\Bbb{C}\Big[
\Bbb{G}(k,n) \Big]$ can be expressed uniquely,
with respect to the cluster ${\bf x}$, as 

\[ { P \ \Big( x_1,\dots,x_{\s{N}} \ ;
\ \Delta^{K_1}, \dots, \Delta^{K_n} \ \Big)
\over {x_1^{i_1} \dots x_{\s{N}}^{i_{\s{N}}} }} \]

\bigskip
\noindent
where $P$ is polynomial, the exponents $i_1, \dots, i_{\s{N}}$
are integers, and $P$ is prime to the monomial 
$x_1^{i_1} \dots x_{\s{N}}^{i_{\s{N}}}$. 
In particular the value of each Pl\"ucker coordinate,
is recovered from such a Laurent expansion. This
implies that $\Phi_{\bf x}$ is injective.

\bigskip
\noindent
Given a point $({\bf z},{\bf w}) \in 
\Big( \Bbb{C}^* \Big)^{\s{N}} \times \Big( \Bbb{C}^* \Big)^n$
let $[{\bf z},{\bf w}]$ be the vector in 
$\bigwedge^k \Bbb{C}^n$ whose
Pl\"ucker coordinates are determined using
the Laurent expansions for the cluster ${\bf x}$.
This point will lie in $X(k,n)$ if 
its coordinates satisfy every Pl\"ucker relation. 
As a cluster algebra, the coordinate ring
$\Bbb{C}\Big[ \Bbb{G}(k,n) \Big]$ is
identified with the abstract cluster algebras
$\mathcal{A}_{\bf x}$ determined by the
seed corresponding to the cluster ${\bf x}$.
Recall that $\mathcal{A}_{\bf x}$ is a
subalgebra of the field of rational
functions in the indeterminates $x \in {\bf \tilde{x}}$,
hence the Laurent polynomials $L_K$ expressing 
the Pl\"ucker coordinates $\Delta^K$ reside
within $\mathcal{A}_{\bf x}$. Since $\mathcal{A}_{\bf x}$
and $\Bbb{C}\Big[ \Bbb{G}(k,n) \Big]$
are isomorphic it follows that any 
polynomial relation satisfied 
by the Pl\"ucker coordinates $\Delta^K$
in $\Bbb{C}\Big[ \Bbb{G}(k,n) \Big]$
must also be satisfied by the corresponding
Laurent polynomials $L_K$. Thus we may
conclude that $[{\bf z},{\bf w}] \in X(k,n)$
and that the map $\Phi_{\bf x}$ is surjective.  
\end{proof}

\bigskip
\noindent
Recall that a point $p$ in the affine cone $X(k,n)$
of the real Grassmannian 
is {\it totally positive} if every Pl\"ucker coordinate
$\Delta^K(p)$ is positive; let $X(k,n)^+$ denote
the variety of positive points. A regular function $f$ in
the homogeneous coordinate ring $\Bbb{C}\Big[ \Bbb{G}(k,n) \Big]$
is {\it positive} if it takes positive values
for each $p \in X(k,n)^+$; in particular all
Pl\"ucker coordinates are positive regular functions.
Cluster variables are also positive owing to
the fact that the value of any cluster
variable is obtained through a sequence
of successive exchange relations -
each of which is subtraction free -
from any fixed cluster consisting entirely of Pl\"ucker
coordinates.

\bigskip
\noindent
If, for a point $p \in X(k,n)$, it so
happens that $x_i(p)>0$ for all $x_i \in {\bf x}$
and $\Delta^K(p)>0$ for all boundary $k$-subsets, then
it must be the case that $p$ is totally positive.
Again this is a consequence of the fact that
all exchange relations are subtraction free and
hence the value of each Pl\"ucker coordinate
is computed as an iterated subtraction free
expression of positive real numbers. Taken
is combination with Corollary \ref{ch2,thm.toric} above we can
surmise: 

\bigskip
\begin{Cor}\label{ch2,cor.positive}
For any cluster ${\bf x}$ in $\Bbb{C}\Big[ \Bbb{G}(k,n) \Big]$
the embedding $\varphi_{\bf x}$ restricts to the embedding

\[ \Phi_{\bf x}: X(k,n)^+ \longrightarrow
\Bbb{R}_{>0}^{\s{N}} \times
\Bbb{R}_{>0}^n . \]

\end{Cor}

\bigskip
\section{Grassmannians of Finite Type}

\bigskip
\noindent
A cluster algebra is of {\it finite type} if
it possess only finitely many cluster variables.
We have already seen an example, namely the
homogeneous coordinate ring of the Grassmannian
$\Bbb{G}(2,n)$. The cluster variables in this 
case are Pl\"ucker coordinates $\Delta^{ij}$
for $1 \leq i<j \leq n$. The results of the
section 5 have shown that the
Pl\"ucker coordinates are cluster
variables of the homogeneous coordinate 
ring of $\Bbb{G}(k,n)$. However, with the
exception of $\Bbb{G}(2,n)$ - and its
complementary Grassmannian $\Bbb{G}(n-2,n)$ -
there will be more cluster variables than
Pl\"ucker coordinates.  
Moreover, of all the Grassmannians $\Bbb{G}(k,n)$ - with
the range $2 < k \leq {n \over 2}$ - only $\Bbb{G}(3,6)$,
$\Bbb{G}(3,7)$, and $\Bbb{G}(3,8)$ have
homogeneous coordinate rings of finite type. 

\bigskip
\noindent
\cite{CA2} establishes a classification of 
finite type cluster algebras which models
the Cartan-Killing classification of
semi-simple Lie algebras. This classification
scheme relies on the following graph-theoretic
construction. 
The principle submatrix
$B$ of an exchange matrix $\tilde{B}$ is skew-symmetric and
it is convenient to visualize it, and its
subsequent mutations, as the
incidence matrix of an
weighted oriented graph $\Gamma \big( B \big)$
with vertices indexed and
labeled by the set of cluster variables ${\bf x}$.
The graph possess an oriented edge from $x$
to $y$, with weight $b_{xy}^2$, if and only if
$b_{xy} > 0$. Using this
graph-theoretic formalism, matrix mutations
of $B$ can be interpreted geometrically by
local rearrangements of $\Gamma(B)$. 
The relationship between this graph and
the type of the corresponding cluster algebra
is summarized by the following rendition
of a result in \cite{CA2}.
 
\bigskip
\begin{Prop}\label{ch2,CA2.4}[Fomin \& Zelevinsky]
A cluster algebra is of finite type if and only
if there exists a seed $\Big( {\bf x}, {\bf c}, 
\tilde{B} \Big)$ such that the underlying
graph (without orientations and weights)
of $\Gamma \big( B \big)$ is a Dynkin diagram.
Moreover, if a cluster algebra possesses
a seed for which the underlying graph of
$\Gamma \big( B \big)$ contains - as
an induced subgraph - an acyclic affine Dynkin diagram
then the cluster algebra is of infinite type.

\end{Prop}

\bigskip
\begin{Thm}\label{ch2,thm.finite}  
$\Bbb{G}(3,6)$, $\Bbb{G}(3,7)$, and
$\Bbb{G}(3,8)$ are the {\bf only} Grassmannians
$\Bbb{G}(k,n)$, within the range $2<k \leq {n \over 2}$,
whose homogeneous coordinate rings are 
cluster algebras of finite type.
\end{Thm}
  
\begin{proof}

\noindent
Consider first the remaining Grassmannians
in the range $2<k \leq {n \over 2}$
Let ${\bf A_{k,n}}$ be the Postnikov
arrangement constructed earlier and let $\Gamma_{k,n}=
\Gamma \big( B_{k,n} \big)$ be the
oriented graph of the
corresponding seed. The vertices of 
and edges of $\Gamma_{k,n}$ form  
a $(k-1) \times (n-k-1)$ square lattice.
This graph contains, as induced subgraph,
either

\newpage
\begin{figure}[ht]

\setlength{\unitlength}{1.0pt}

\begin{center}
\begin{picture}(90,90)(0,100)

\thinlines


\thinlines
\put(0,100){\line(1,0){80}}
\put(0,180){\line(1,0){80}}
\put(40,100){\line(0,1){80}}
\put(22,100){\line(-1,1){3}}
\put(22,100){\line(-1,-1){3}}
\put(58,100){\line(1,1){3}}
\put(58,100){\line(1,-1){3}}
\put(22,180){\line(-1,1){3}}
\put(22,180){\line(-1,-1){3}}
\put(58,180){\line(1,1){3}}
\put(58,180){\line(1,-1){3}}
\put(40,122){\line(1,-1){3}}
\put(40,122){\line(-1,-1){3}}
\put(40,158){\line(1,1){3}}
\put(40,158){\line(-1,1){3}}

\put(0,100){\line(0,1){80}}
\put(80,100){\line(0,1){80}}
\put(0,140){\line(1,0){80}}
\put(18,140){\line(1,1){3}}
\put(18,140){\line(1,-1){3}}
\put(62,140){\line(-1,-1){3}}
\put(62,140){\line(-1,1){3}}
\put(0,118){\line(1,1){3}}
\put(0,118){\line(-1,1){3}}
\put(80,118){\line(1,1){3}}
\put(80,118){\line(-1,1){3}}
\put(0,162){\line(-1,-1){3}}
\put(0,162){\line(1,-1){3}}
\put(80,162){\line(-1,-1){3}}
\put(80,162){\line(1,-1){3}}


\put(-40,0){\line(1,0){160}}
\put(0,0){\line(0,1){40}}
\put(80,0){\line(0,1){40}}
\put(-22,0){\line(1,1){3}}
\put(-22,0){\line(1,-1){3}}
\put(22,0){\line(-1,-1){3}}
\put(22,0){\line(-1,1){3}}
\put(58,0){\line(1,1){3}}
\put(58,0){\line(1,-1){3}}
\put(102,0){\line(-1,-1){3}}
\put(102,0){\line(-1,1){3}}
\put(0,18){\line(1,1){3}}
\put(0,18){\line(-1,1){3}}
\put(80,18){\line(1,1){3}}
\put(80,18){\line(-1,1){3}}

\put(35,70){\bf \text{or}}

\put(-40,40){\line(1,0){160}}
\put(-40,0){\line(0,1){40}}
\put(40,0){\line(0,1){40}}
\put(120,0){\line(0,1){40}}
\put(-18,40){\line(-1,1){3}}
\put(-18,40){\line(-1,-1){3}}
\put(18,40){\line(1,1){3}}
\put(18,40){\line(1,-1){3}}
\put(62,40){\line(-1,1){3}}
\put(62,40){\line(-1,-1){3}}
\put(98,40){\line(1,1){3}}
\put(98,40){\line(1,-1){3}}
\put(-40,22){\line(-1,-1){3}}
\put(-40,22){\line(1,-1){3}}
\put(40,22){\line(-1,-1){3}}
\put(40,22){\line(1,-1){3}}
\put(120,22){\line(-1,-1){3}}
\put(120,22){\line(1,-1){3}}

\end{picture} 
\end{center}
\end{figure}

\bigskip
\bigskip
\bigskip
\bigskip
\bigskip
\bigskip
\bigskip
\bigskip
\bigskip
\bigskip
\bigskip
\bigskip
\noindent
These two subgraphs contain, as induced subgraphs,
a copy of the affine Dynkin diagram 
$\widehat{D_6}$; this subgraph is highlighted
in both cases below:

\bigskip
\begin{figure}[ht]

\setlength{\unitlength}{1.0pt}

\begin{center}
\begin{picture}(90,50)(0,130)

\thinlines


\thinlines
\put(0,100){\line(1,0){80}}
\put(0,180){\line(1,0){80}}
\put(40,100){\line(0,1){80}}
\put(22,100){\line(-1,1){3}}
\put(22,100){\line(-1,-1){3}}
\put(58,100){\line(1,1){3}}
\put(58,100){\line(1,-1){3}}
\put(22,180){\line(-1,1){3}}
\put(22,180){\line(-1,-1){3}}
\put(58,180){\line(1,1){3}}
\put(58,180){\line(1,-1){3}}
\put(40,122){\line(1,-1){3}}
\put(40,122){\line(-1,-1){3}}
\put(40,158){\line(1,1){3}}
\put(40,158){\line(-1,1){3}}

\thicklines
\put(0,100){\line(0,1){80}}
\put(80,100){\line(0,1){80}}
\put(0,140){\line(1,0){80}}
\put(18,140){\line(1,1){3}}
\put(18,140){\line(1,-1){3}}
\put(62,140){\line(-1,-1){3}}
\put(62,140){\line(-1,1){3}}
\put(0,118){\line(1,1){3}}
\put(0,118){\line(-1,1){3}}
\put(80,118){\line(1,1){3}}
\put(80,118){\line(-1,1){3}}
\put(0,162){\line(-1,-1){3}}
\put(0,162){\line(1,-1){3}}
\put(80,162){\line(-1,-1){3}}
\put(80,162){\line(1,-1){3}}


\thicklines
\put(-40,0){\line(1,0){160}}
\put(0,0){\line(0,1){40}}
\put(80,0){\line(0,1){40}}
\put(-22,0){\line(1,1){3}}
\put(-22,0){\line(1,-1){3}}
\put(22,0){\line(-1,-1){3}}
\put(22,0){\line(-1,1){3}}
\put(58,0){\line(1,1){3}}
\put(58,0){\line(1,-1){3}}
\put(102,0){\line(-1,-1){3}}
\put(102,0){\line(-1,1){3}}
\put(0,18){\line(1,1){3}}
\put(0,18){\line(-1,1){3}}
\put(80,18){\line(1,1){3}}
\put(80,18){\line(-1,1){3}}


\thinlines
\put(-40,40){\line(1,0){160}}
\put(-40,0){\line(0,1){40}}
\put(40,0){\line(0,1){40}}
\put(120,0){\line(0,1){40}}
\put(-18,40){\line(-1,1){3}}
\put(-18,40){\line(-1,-1){3}}
\put(18,40){\line(1,1){3}}
\put(18,40){\line(1,-1){3}}
\put(62,40){\line(-1,1){3}}
\put(62,40){\line(-1,-1){3}}
\put(98,40){\line(1,1){3}}
\put(98,40){\line(1,-1){3}}
\put(-40,22){\line(-1,-1){3}}
\put(-40,22){\line(1,-1){3}}
\put(40,22){\line(-1,-1){3}}
\put(40,22){\line(1,-1){3}}
\put(120,22){\line(-1,-1){3}}
\put(120,22){\line(1,-1){3}}

\end{picture} 
\end{center}
\end{figure}

\bigskip
\bigskip
\bigskip
\bigskip
\bigskip
\bigskip
\bigskip
\bigskip
\bigskip
\bigskip
\bigskip
\bigskip
\bigskip
\bigskip
\bigskip
\noindent
Proposition \ref{ch2,CA2.4} above
implies that the coordinate ring
of $\Bbb{G}(k,n)$ must be of 
infinite type, with the exceptions of
$\Bbb{G}(3,6)$, $\Bbb{G}(3,7)$, and
$\Bbb{G}(3,8)$.  
In the remaining three cases we begin by mutating
the ${\bf A_{k,n}}$ seed until the 
underlying graph of $\Gamma$ is
a Dynkin diagram. The requisite sequence of mutations
performed in each case is not necessarily
the most efficient route, nevertheless each
Dynkin diagram obtained
has the added feature that it is 
oriented as a bi-partide graph; each vertex
is either a sink or a source.  
This property is essential in working
out the route correspondence later in this
section.  

\bigskip
\noindent
In order to keep track of matrix mutation
the vertices are also labeled with auxiliary
integer indices.

\bigskip
\noindent
{\bf $\Bbb{C}\Big[ \Bbb{G}(3,6) \Big]$ case:}
The oriented graph $\Gamma_{3,6}$ associated with the 
${\bf A_{3,6}}$ cluster -
labeled with the corresponding
cluster of Pl\"ucker coordinates -
is depicted below.

\begin{figure}[ht]

\setlength{\unitlength}{1.0pt}

\begin{center}
\begin{picture}(90,90)(0,130)

\thinlines
  
\thicklines

\put(0,100){\line(1,0){100}}
\put(45,100){\line(1,1){3}}
\put(45,100){\line(1,-1){3}}

\put(0,200){\line(1,0){100}}
\put(50,200){\line(-1,1){3}}
\put(50,200){\line(-1,-1){3}}

\put(0,100){\line(0,1){100}}
\put(0,150){\line(-1,-1){3}}
\put(0,150){\line(1,-1){3}}
\put(100,150){\line(-1,1){3}}
\put(100,150){\line(1,1){3}}

\put(100,100){\line(0,1){100}}

\put(-5,92){$4$} \put(5,105){$\Delta^{356}$}
\put(-5,202){$3$} \put(5,188){$\Delta^{136}$}
\put(100,202){$1$} \put(75,188){$\Delta^{236}$}
\put(100,92){$2$} \put(75,105){$\Delta^{346}$}

\end{picture}
\end{center}
\end{figure}

\bigskip
\bigskip
\bigskip
\bigskip
\bigskip
\bigskip
\noindent
Below is the $D_4$ graph obtained 
from the ${\bf A_{3,6}}$ 
cluster by performing
the sequence of mutations $4, 2, 4, 1$.

\begin{figure}[ht]

\setlength{\unitlength}{1.0pt}

\begin{center}
\begin{picture}(90,90)(50,130)

\thinlines
  
\thicklines

\put(0,100){\line(1,1){100}}
\put(50,150){\line(-1,0){3}}
\put(50,150){\line(0,-1){3}}

\put(100,200){\line(0,-1){100}}
\put(100,145){\line(-1,-1){3}}
\put(100,145){\line(1,-1){3}}

\put(100,200){\line(1,-1){100}}
\put(150,150){\line(0,-1){3}}
\put(150,150){\line(1,0){3}}

\put(-7,90){$1$}  
\put(98,90){$3$}  
\put(203,90){$4$} 
\put(98,203){$2$} 

\end{picture}
\end{center}
\end{figure}

\newpage
\noindent
{\bf $\Bbb{C} \Big[ \Bbb{G}(3,7) \Big]$ case:}
The cluster-labeled oriented graph in this case
is

\begin{figure}[ht]

\setlength{\unitlength}{1.0pt}

\begin{center}
\begin{picture}(90,90)(50,110)

\thinlines
  
\thicklines

\put(0,100){\line(1,0){100}}
\put(45,100){\line(1,1){3}}
\put(45,100){\line(1,-1){3}}

\put(0,200){\line(1,0){100}}
\put(50,200){\line(-1,1){3}}
\put(50,200){\line(-1,-1){3}}

\put(0,100){\line(0,1){100}}
\put(0,150){\line(-1,-1){3}}
\put(0,150){\line(1,-1){3}}
\put(100,150){\line(-1,1){3}}
\put(100,150){\line(1,1){3}}

\put(100,100){\line(0,1){100}}
\put(100,200){\line(1,0){100}}
\put(150,200){\line(1,1){3}}
\put(150,200){\line(1,-1){3}}

\put(200,200){\line(0,-1){100}}
\put(200,150){\line(1,-1){3}}
\put(200,150){\line(-1,-1){3}}

\put(100,100){\line(1,0){100}}
\put(150,100){\line(-1,1){3}}
\put(150,100){\line(-1,-1){3}}

\put(-5,91){$5$} \put(5,105){$\Delta^{267}$}
\put(-5,202){$3$} \put(5,188){$\Delta^{126}$}
\put(97,202){$2$} \put(75,188){$\Delta^{125}$}
\put(97,91){$4$} \put(75,105){$\Delta^{256}$}
\put(200,91){$1$} \put(175,188){$\Delta^{235}$}
\put(200,202){$6$} \put(175,105){$\Delta^{356}$}

\end{picture}
\end{center}
\end{figure}

\bigskip
\bigskip
\bigskip
\noindent
After mutating the ${\bf A_{3,7}}$ cluster
through the sequence of mutations corresponding
to the vertices $2, 4, 3, 5, 6, 5, 1$ one
obtains a cluster whose associated graph
is $E_6$. This cluster is shown below:

\begin{figure}[ht]

\setlength{\unitlength}{1.0pt}

\begin{center}
\begin{picture}(90,90)(50,130)

\thinlines
  
\thicklines

\put(0,200){\line(0,-1){100}}
\put(0,150){\line(-1,-1){3}}
\put(0,150){\line(1,-1){3}}

\put(0,100){\line(1,1){100}}
\put(50,150){\line(0,-1){3}}
\put(50,150){\line(-1,0){3}}

\put(100,200){\line(0,-1){100}}
\put(100,150){\line(-1,-1){3}}
\put(100,150){\line(1,-1){3}}

\put(100,200){\line(1,-1){100}}
\put(150,150){\line(0,-1){3}}
\put(150,150){\line(1,0){3}}

\put(200,100){\line(0,1){100}}
\put(200,150){\line(-1,-1){3}}
\put(200,150){\line(1,-1){3}}

\put(-2,202){$1$}

\put(-3,91){$3$}

\put(97,91){$2$}

\put(97,202){$4$}

\put(198,202){$6$}

\put(198,91){$5$}

\end{picture}
\end{center}
\end{figure}

\bigskip
\bigskip
\bigskip
\bigskip
\bigskip
\noindent
{\bf $\Bbb{C}\Big[ \Bbb{G}(3,8) \Big]$ case:}
We begin with the initial cluster arising from the
$\pi_{3,8}$-diagram described in the previous section. The 
incidence graph for this cluster, labeled
with its cluster of Pl\"ucker coordinates is depicted below:

\newpage
\begin{figure}[ht]

\setlength{\unitlength}{1.0pt}

\begin{center}
\begin{picture}(90,90)(100,120)

\thinlines
  
\thicklines

\put(0,100){\line(1,0){300}}
\put(45,100){\line(1,1){3}}
\put(45,100){\line(1,-1){3}}
\put(150,100){\line(-1,-1){3}}
\put(150,100){\line(-1,1){3}}
\put(245,100){\line(1,1){3}}
\put(245,100){\line(1,-1){3}}

\put(0,200){\line(1,0){300}}
\put(50,200){\line(-1,1){3}}
\put(50,200){\line(-1,-1){3}}
\put(145,200){\line(1,-1){3}}
\put(145,200){\line(1,1){3}}
\put(250,200){\line(-1,1){3}}
\put(250,200){\line(-1,-1){3}}

\put(0,100){\line(0,1){100}}
\put(0,150){\line(-1,-1){3}}
\put(0,150){\line(1,-1){3}}
\put(200,150){\line(-1,-1){3}}
\put(200,150){\line(1,-1){3}}
\put(100,150){\line(-1,1){3}}
\put(100,150){\line(1,1){3}}
\put(300,150){\line(-1,1){3}}
\put(300,150){\line(1,1){3}}

\put(100,100){\line(0,1){100}}
\put(200,100){\line(0,1){100}}
\put(300,100){\line(0,1){100}}

\put(-5,92){$8$} \put(5,105){$\Delta^{124}$}
\put(-5,202){$2$} \put(5,188){$\Delta^{134}$}
\put(98,202){$1$} \put(105,188){$\Delta^{145}$}
\put(98,92){$3$} \put(105,105){$\Delta^{148}$}
\put(198,92){$4$} \put(175,105){$\Delta^{158}$}
\put(198,202){$7$} \put(175,188){$\Delta^{458}$}
\put(300,202){$6$} \put(275,188){$\Delta^{568}$}
\put(300,92){$5$} \put(275,105){$\Delta^{578}$}

\end{picture}
\end{center}
\end{figure}

\bigskip
\bigskip
\bigskip
\noindent
The $E_8$ graph can be obtain by mutating the
incidence graph above through the sequence
of nodes (starting with) $1$ then $3$, $7$, $6$, $5$,
$2$, $4$, $3$, $8$, $7$, $6$, and finally $8$. 
It is decorated by the following cluster of functions:

\begin{figure}[ht]

\setlength{\unitlength}{1.0pt}

\begin{center}
\begin{picture}(90,90)(100,130)

\thinlines
  
\thicklines

\put(0,200){\line(0,-1){100}}
\put(0,150){\line(-1,-1){3}}
\put(0,150){\line(1,-1){3}}

\put(0,100){\line(1,1){100}}
\put(50,150){\line(0,-1){3}}
\put(50,150){\line(-1,0){3}}

\put(100,200){\line(0,-1){100}}
\put(100,150){\line(-1,-1){3}}
\put(100,150){\line(1,-1){3}}

\put(100,200){\line(1,-1){100}}
\put(150,150){\line(0,-1){3}}
\put(150,150){\line(1,0){3}}

\put(200,100){\line(0,1){100}}
\put(200,150){\line(-1,-1){3}}
\put(200,150){\line(1,-1){3}}

\put(200,200){\line(1,-1){100}}
\put(250,150){\line(0,-1){3}}
\put(250,150){\line(1,0){3}}

\put(300,100){\line(0,1){100}}
\put(300,150){\line(-1,-1){3}}
\put(300,150){\line(1,-1){3}}

\put(-2,202){$1$}

\put(-3,91){$3$}

\put(97,91){$2$}

\put(97,202){$4$}

\put(198,202){$6$}

\put(198,91){$5$}

\put(298,91){$7$}

\put(298,202){$8$}

\end{picture}
\end{center}
\end{figure}

\end{proof}

\newpage
\section*{\it{Projective Plane Geometry and Cluster Variables:}}

\noindent
In order to geometrically interpret the cluster
variables attached to the Grassmannians 
$\Bbb{G}(3,6)$, $\Bbb{G}(3,7)$, and $\Bbb{G}(3,8)$,
it is convenient to represent points of these Grassmannians
as ordered configurations of, respectively, six,
seven, and eight points in the real projective 
plane $\Bbb{R}\Bbb{P}^2$. Express each point
in $\Bbb{G}(3,n)$ as a $3 \times n$ matrix 
$\big( v_1 \ \dots \ v_n \big)$ where $v_i$ is
a vector residing in $\Bbb{R}^3 - {\bf 0}$.
To this matrix, associate the ordered configuration
$[v_1], \dots, [v_n]$ of points in $\Bbb{R}\Bbb{P}^2$.

\bigskip
\noindent
Recall that if $[v]$ and $[v']$ are points in the projective
plane, then line $[l]$ joining them, written projectively as
a homogeneous linear form, is $[v \times v']$ where
$\times$ is the cross product in $\Bbb{R}^3$. Dually,
if $[l]$ and $[l']$ are projective lines then their
point of intersection is also given by $[l \times l']$.
The following elementary cross product identities
will be useful in this analysis:

\begin{equation}   u \cdot \big(v \times w \big) \ = \
   \big(u \times v \big) \cdot w \ = \  
   \det \Big( u \ \ v \ \ w \Big) \end{equation}

\begin{equation}   \big(u \times v \big) \cdot \big(
   w \times z \big) \ = \
   \begin{pmatrix} u \cdot w & u \cdot z \\
   v \cdot w & v \cdot z \end{pmatrix} \end{equation}

\bigskip
\noindent
In view of equation (2), three projective points
$[v_i]$, $[v_j]$ and $[v_k]$ drawn from a
configuration $[v_1], \dots, [v_n]$ are colinear
if and only if $\det \big(v_i \ v_j \ v_k \big) = 0$;
in other words the Pl\"ucker coordinate $\Delta^{ijk}$
vanishes on those points in the Grassmannian
$\Bbb{G}(3,n)$ for which the $i$-th, $j$-th,
and $k$-th members of the associated configuration
of projective points are colinear.

\bigskip
\noindent
The aim here is to identify 
the cluster variables which are not Pl\"ucker
coordinates and to describe their vanishing 
loci. In cataloging these functions,
it is useful to take into account
certain combinatorial symmetries which
are present. The dihedral group $D_n$ acts on
the affine cone $X(k,n)$ of the Grassmannian by simply
permuting the $n$ columns of a $k \times n$ matrix
in a manner consistent with the action
of $D_n$ on $[1 \dots n]$. In the cases studied
below, all cluster variables have the 
property that their dihedral translates
are, up to sign, also cluster variables;
so its is enough to determine representatives
of the dihedral orbits.

\bigskip
\begin{Thm}\label{ch2,thm.G6}
As a cluster algebra, $\Bbb{C} \Big[ \Bbb{G}(3,6) \Big]$
possesses sixteen cluster variables. Fourteen of these are
the Pl\"ucker coordinates $\Delta^{ijk}$ - where
$\big\{i,j,k \big\} \subset [1 \dots 6]$ is an internal $3$-subset - 
and the remaining two, denoted as 
$X^{123456}$ and $Y^{123456}$, are quadratic regular functions. 
The zero locus of $X^{123456}$
consists of those configurations of six points
$[v_1], \dots, [v_6]$ in $\Bbb{R}\Bbb{P}^2$ for 
which the three lines - joining $[v_1]$ with $[v_2]$, $[v_3]$
with $[v_4]$, and $[v_5]$ with $[v_6]$ - have a common
intersection point, as depicted below:

\newpage
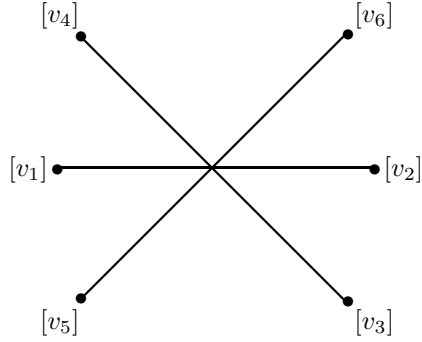
\begin{figure}[ht]

\setlength{\unitlength}{1.0pt}

\begin{center}
\begin{picture}(100,100)(0,130)

\thinlines
  
\thicklines
\put(50,170){\line(1,1){50}}
\put(50,170){\line(-1,-1){50}}
\put(50,170){\line(1,-1){50}}
\put(50,170){\line(-1,1){50}}
\put(50,170){\line(-1,0){60}}
\put(50,170){\line(1,0){60}}

\put(-11,167){$\bullet$}
\put(-27,167){$[v_1]$}

\put(109,167){$\bullet$}
\put(115,167){$[v_2]$}

\put(99,218){$\bullet$}
\put(104,225){$[v_6]$}

\put(-2,118){$\bullet$}
\put(-15,109){$[v_5]$}

\put(-2,217){$\bullet$}
\put(-15,225){$[v_4]$}

\put(99,117){$\bullet$}
\put(104,109){$[v_3]$}

\end{picture}

\bigskip
\bigskip
\bigskip
\bigskip
\caption{Vanishing Locus of $X^{123456}$}
\label{G6 locus}
\end{center}
\end{figure}

\noindent
The cluster variable $Y^{123456}$ satisfies the
identity 

\[ Y^{123456} \big( v_1, v_2, v_3, v_4, v_5, v_6 \big) \ = \
X^{123456} \big( v_6, v_1, v_2,  v_3, v_4, v_5 \big) \] 

\bigskip
\noindent 
and consequently its vanishing locus is the
same as $X^{123456}$'s except that the indices
are cyclically shifted.

\end{Thm}

\begin{proof}

\bigskip
\noindent
We begin with the cluster corresponding to the 
$D_4$ graph constructed in
Theorem \ref{ch2,thm.finite}; hereafter called the $D_4$-cluster. 
If one carefully executes all of the exchange 
relations prescribed by 
the listed sequence of mutations from the 
initial ${\bf A_{3,6}}$-cluster, then
the cluster variable attached to node \#2 of $D_4$,
now denoted as $Y^{123456}$, is defined by
the exchange relation 
   
\[ \Delta^{346} \ Y^{123456} \ = \ \Delta^{146} \Delta^{236} 
\Delta^{345} \ + \ \Delta^{136} \Delta^{234} 
\Delta^{456}. \]

\bigskip
\noindent
The Pl\"ucker coordinate $\Delta^{346}$ can
be factored from the right side using
a combination of short Pl\"ucker relations
and after cancellation one obtains the following
explicit formula for $Y^{123456}$:

\[ Y^{123456} \ = \ \Delta^{236} \Delta^{145} \ - \
\Delta^{123} \Delta^{456} \ = \
\det \begin{pmatrix} \Delta^{236} & \Delta^{456}
\\ \Delta^{123} & \Delta^{145} \end{pmatrix} \]

\bigskip
\noindent
Using cross product identities (2) and (3) one can re-express
$Y^{123456}$ as the compound determinant

\[ \det \Big( v_6 \times v_1 \ \ \ v_2 \times v_3 \ \ \ 
   v_4 \times v_5  \Big). \]

\bigskip
\noindent
This determinental formula implies that
whenever $[v_1],\dots [v_6]$ are distinct,
the zero locus of $Y^{123456}$ coincides with
point configurations for which the lines joining
$[v_1]$ with $[v_6]$, $[v_2]$ with $[v_3]$,
and $[v_4]$ with $[v_5]$
having a common intersection point. 

\bigskip
\noindent
Aside of $Y^{123456}$, there is only one
remaining cluster function which
is not a Pl\"ucker coordinate. One can
show that this
function, denoted $X^{123456}$, also can
be expressed as a compound determinant, namely

   \begin{eqnarray} X^{123456} 
   & \ = \ & \Delta^{134} \Delta^{256} -
   \Delta^{156} \Delta^{234} \\
   & \ = \ & \det \Big( v_1 \times v_2 \ \ \ v_3 \times v_4 \ \ \ 
   v_5 \times v_6  \Big) . \end{eqnarray}

\noindent
For distinct projective points $[v_1],\dots,[v_6]$,
the zero locus of $X^{123456}$ coincides
with point configurations for which
the three lines formed by joining
$[v_1]$ with $[v_2]$, $[v_3]$ with $[v_4]$,
and $[v_5]$ with $[v_6]$ respectively
have a common intersection point.

\end{proof}

\bigskip
\noindent
In what follows assume $n > 6$ and that $I$ is
a subset of $[1 \dots n]$ whose complement
has cardinality six. Let
$X^{[1 \dots n] - I}$ denote the composition of
$X^{123456}$ with the projection 
the $\Bbb{G}(3,n) \stackrel{I}{\longrightarrow}
\Bbb{G}(3,6)$ which forgets the 
columns indexed by $I$ of the $3 \times n$ matrix representing
a point in $\Bbb{G}(3,n)$; define  $Y^{[1 \dots n] -I}$
likewise.

\bigskip
\begin{Thm}\label{ch2,thm.G7} The cluster algebra
$\Bbb{C} \Big[ \Bbb{G}(3,7) \Big]$ possesses
forty two cluster variables. Twenty eight
of these are the Pl\"ucker coordinates 
$\Delta^{ijk}$ - where $\big\{ i,j,k \big\} \subset
[1 \dots 7]$ is an internal $3$-subset - and the
remaining fourteen are the quadratic regular
functions $X^{[1 \dots 7] - i}$ and $Y^{[1 \dots 7] -i}$
defined above for $i \in [1 \dots 7]$. 

\end{Thm}  

\begin{proof}

\bigskip
\noindent
Using exchange relations one computes the
cluster variables attached to the 
nodes of the $E_6$-cluster from
Theorem \ref{ch2,thm.finite}. Specifically, the cluster variables 
corresponding to nodes \#4 and \#5 are
precisely the regular functions  
$X^{123567}$ and $X^{123467}$.

\end{proof}

\bigskip
\begin{Thm}\label{ch2,thm.G8}
The cluster algebra $\Bbb{C}\Big[ \Bbb{G}(3,8) \Big]$
possesses 128 cluster variables. Of these forty eight
are Pl\"ucker coordinates $\Delta^{ijk}$ -
where $\big\{ i,j,k \big\} \subset [1 \dots n]$ 
is an internal $3$-subset. Fifty six cluster variables are the
quadratic regular functions $X^{[1 \dots 8] - \{ij\}}$
and $Y^{[1 \dots 8] - \{ij\}}$ - for $1 \leq i < j \leq 8$ -
inherited from $\Bbb{G}(3,6)$. The remaining twenty four 
are dihedral translates of two cubic regular functions,
denoted as $A$ and $B$. The zero locus of $A$ consists
of configuration of eight projective points 
$[v_1], \dots, [v_8]$ for which the points
$p$, $q$ and $[v_5]$ are colinear (see illustration below).

\newpage
\begin{figure}[ht]

\setlength{\unitlength}{1.0pt}

\begin{center}
\begin{picture}(90,90)(20,100)

\thinlines
  
\thinlines

\put(15,0){\line(1,2){75}}
\put(30,0){\line(0,1){180}}
\put(0,90){\line(1,0){130}}
\put(0,180){\line(1,-1){130}}

\thicklines
\put(0,0){\line(1,1){130}}

\put(60,59){$\bullet$}
\put(64,53){$[v_5]$}

\put(27.5,87.5){$\bullet$}
\put(34,96){$[v_1]$}

\put(27.5,177.5){$\bullet$}
\put(34,178){$[v_2]$}

\put(127.5,87.5){$\bullet$}
\put(134,96){$[v_8]$}

\put(-2,177.5){$\bullet$}
\put(-18,178){$[v_6]$}

\put(127.5,48.5){$\bullet$}
\put(134,53){$[v_7]$}

\put(13,-2){$\bullet$}
\put(9,-10){$[v_4]$}

\put(88,148){$\bullet$}
\put(83,156){$[v_3]$}

\put(87.8,87.5){$\bullet$}
\put(89,96){$q$}

\put(27.5,27.5){$\bullet$}
\put(35,23){$p$}

\end{picture}

\bigskip
\bigskip
\bigskip
\bigskip
\bigskip
\bigskip
\bigskip
\bigskip
\bigskip
\bigskip
\bigskip
\bigskip
\bigskip
\caption{Vanishing Locus of $A$}
\label{G8 locus}
\end{center}
\end{figure}
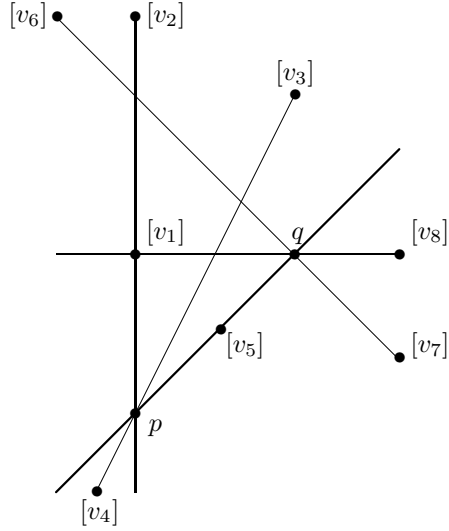

\noindent
By definition, $p$ is the intersection point of the line
joining $[v_1]$ with $[v_2]$ and the line
joining $[v_3]$ with $[v_4]$; similarly $q$ is the 
the intersection point of the line joining
$[v_1]$ with $[v_8]$ and the line joining $[v_6]$
with $[v_7]$. 

\end{Thm}

\bigskip
\noindent
where $Y^{234678}$ and $Y^{135678}$ denoted the function
$Y$ inherited from the projections
$\Bbb{G}(3,8) \stackrel{15,24}{\longrightarrow}
\Bbb{G}(3,6)$ which, respectively, forget columns
$1 \ \& \ 2$ and columns $2 \ \& \ 4$ of a $3 \times 8$ matrix.
For $i,j \in [1 \dots 8]$, 
each projection $\Bbb{G}(3,8) \stackrel{ij}{\longrightarrow}
\Bbb{G}(3,6)$ gives rise to the
pair of cluster functions $X^{[1 \dots 8] - ij}$
and $Y^{[1 \dots 8] - ij}$.
The functions $A$ and $B$ cubic and are described implicitly by
the following exchange relations:

\[ \Delta^{578} A \ = \ \Delta^{178} \Delta^{567} X^{123458} \
+ \ \Delta^{158} \Delta^{678} X^{123457} \] 

\[ \Delta^{158} B \ = \ \Delta^{128} \Delta^{567} 
   X^{123458} \ + \ \Delta^{258} A \] 

\bigskip
\noindent
Explicit formulas for $A$ and $B$ are
obtained using the cross product identities mentioned in the 
case of $\Bbb{G}(3,6)$:

\[ A \ = \ \Delta^{134} \Big( \Delta^{258} \Delta^{167}
   - \Delta^{678} \Delta^{125} \Big) - \Delta^{158}
   \Delta^{234} \Delta^{167} \] 

\[ B \ = \ \Delta^{258} \Delta^{134} \Delta^{267} -
   \Delta^{234} \Big( \Delta^{128} \Delta^{567} +
   \Delta^{258} \Delta^{167} \Big) \]
 
\bigskip
\noindent
The formula for $B$ can be rewritten using short Pl\"ucker
relations to obtain

\begin{eqnarray} B & \ = \ & \Delta^{258} \Delta^{134} \Delta^{267} -
   \Delta^{234} \Big( \Delta^{158} \Delta^{267} +
   \Delta^{678} \Delta^{125} \Big) \\
   & \ = \ & F_1 \Big( A \Big) \   
\text{for the derivation $F_1 \in \frak{sl}_8$.} \end{eqnarray} 

\bigskip
\noindent
As before, view a point in the Grassmannian $\Bbb{G}(3,8)$ as
an ordered configuration of eight points $v_1 \dots
v_8$ in $\Bbb{RP}^2$. Let $p$ be the intersection
point of the line joining the points $v_1$ $v_2$ and the
line joining points $v_3$ and $v_4$. If the each point $v_i$
is represented by a vector in $\Bbb{R}^3 - {\bf 0}$ then
as a vector 

\[ p \ = \ \Big( v_1 \times v_2 \Big) \times 
   \Big( v_3 \times v_4 \Big). \] 

\bigskip
\noindent
Similarly, let $q$ be the intersection point
of the line joining points $v_1$ and $v_8$ and 
the line joining points $v_6$ and $v_7$. The point $q$ is 
given represented by

\[ q \ = \ \Big( v_1 \times v_8 \Big) \times 
   \Big( v_6 \times v_7 \Big). \] 

\bigskip
\begin{Claim}
The cluster function $A$ vanishes whenever $v_5$ lies
on the line joining $p$ and $q$.
Moreover, $A$ can be expressed as the compound
determinant

\[ \det \Bigg( \big(v_1 \times v_2 \big) \times 
   \big( v_3 \times v_4 \big) \ \
   \big( v_1 \times v_8 \big) \times 
   \big( v_6 \times v_7 \big) \ \
   v_5 \Bigg). \]

\end{Claim}

\noindent
Interpreted geometrically, this claim
stipulates that the zero locus of $A$ 
consists of point configurations 
in $\Bbb{R}\Bbb{P}^2$ having the following
intersection type:

\bigskip
\noindent
Moreover, any auxiliary configuration degeneracies 
predicted by the vanishing of 
$\det \big( \ p \ \  q \ \ v_5 \ \big) $ -
e.g. the line joining $[v_1]$ with $[v_2]$ coinciding
with the line joining $[v_3]$ with $[v_4]$ -
contribute to $A's$ zero locus.

\begin{proof}
Since $A$ varies linearly with respect to the 
vector $v_5$ it is enough to verify that
$A$ vanishes whenever $v_5$ coincides with
$p$ or $q$. Setting $v_5=p$ in the expression
for $A$ we obtain:

\[ \Delta^{134} \Bigg( \Delta^{167} 
\det \Big( v_2 \ \ (v_1 \times v_2) \times 
(v_3 \times v_4) \ \ v_8 \Big) -
\Delta^{678} \det \Big( v_1 \ \ v_2 \ \ 
(v_1 \times v_2) \times (v_3 \times v_4) \Big) \
\Bigg) \]

\[ - \Delta^{234} \Delta^{167} 
\det \Big( v_1 \ \ (v_1 \times v_2) \times
(v_3 \times v_4) \ \ v_8 \Big) \]

\[ = \Delta^{134} \Bigg( -\Delta^{167}
\Big( \ (v_2 \times v_8) \cdot \Big[ 
(v_1 \times v_2) \times (v_3 \times v_4) \Big] \ \Big)
- \Delta^{678} \Big( \ (v_1 \times v_2) \cdot
\Big[ (v_1 \times v_2) \times (v_3 \times v_4) \Big] \ \Big)
\ \Bigg) \]

\[ \Delta^{234} \Delta^{167} \Big( \ (v_1 \times v_8)
\cdot \Big[ (v_1 \times v_2) \times (v_3 \times v_4) 
\Big] \ \Big) \]

\[ = \ - \Delta^{134} \Delta^{167} 
\det \begin{pmatrix} 0 & \Delta^{234} \\ \Delta^{128}
& \Delta^{348} \end{pmatrix} \ + \
\Delta^{234} \Delta^{167} \begin{pmatrix}
0 & \Delta^{134} \\ \Delta^{128} & \Delta^{348} \end{pmatrix} \
= \ 0 \] 

\bigskip
\noindent
The proof that $A$ vanishes for $v_5 = q$ is nearly
identical to the calculation above. The determinant
$\det \big(\ p \ \ q \ \ v_5  \ \big) $ also measures
colinearity of the points $p$, $q$, and $[v_5]$
and so its zero locus must be the same as $A$'s.
It follows that $A$ has
the prescribed determinental form.
\end{proof}

\bigskip
\noindent
The point configurations for which the
cluster function $B$ vanishes are identical
to $A$'s except the roles of $[v_1]$
and $[v_2]$ are interchanged. This is 
a consequence of the formula
$B = F_1(A)$.

\section*{\it{Root Correspondence:}}

\noindent
In \cite{CA2} Fomin and Zelevinsky construct
a bijection between cluster variables of a cluster algebra
of finite type and the {\it almost positive} roots of
the associated root system. The bijection is constructed
by first choosing a cluster whose incidence graph is 
the Dynkin diagram of the associated cluster algebra
type; call this a {\it Dynkin cluster}. Let $x_i$
denote the variable in this cluster attached to the $i$-th
node of the Dynkin graph.
In virtue of the Laurent property, every cluster
variable $X$ can be expanded as a unique Laurent polynomial
of the form

\[ { P \Big( x_1, \dots, x_{\s{N}}  \ ; \ c \in {\bf c} \Big) \over
{ x_1^{i_1} \cdots x_{\s{N}}^{i_{\s{N}}} }} \] 

\bigskip
\noindent
where $N$ is the rank of the cluster algebra (i.e.
rank of the corresponding root system), $P$ is a polynomial and the 
exponents $i_1 \dots i_N$ are integers. The bijection pairs
the cluster function $X$ with the almost
positive root 

\[ i_1 \alpha_1 \ + \ \dots \ + \ i_{\s{N}} \alpha_{\s{N}} \]

\bigskip
\noindent
associated to the vector of exponents
in the monomial denominator of the Laurent expansion above.
In particular, the variable $x_i$ of the Dynkin cluster
is paired with the negative simple root $-\alpha_i$.

\bigskip
\noindent
This bijection relies on two piece-wise
linear involutions $\tau_+$ and $\tau_-$,
developed in \cite{Ysystems},
of the root lattice $Q$ obtained 
by {\it tropicalizing} the exchange relations
for the Dynkin cluster. For $\alpha \in Q$
let $[\alpha:\alpha_i$ be the coefficient
of $\alpha_i$ in the expansion of $\alpha$
in the basis of simple roots. 
The coefficients of the element $\tau_{\pm} \alpha$ of $Q$ are
given by:

\[ \big[ \tau_{\pm} \alpha : \alpha_i \big] \ = 
\ \left\{ \begin{array}{ll}
-\big[\alpha:\alpha_i\big] -
\sum_{j \ne i} \ a_{ij} \max\Big(\big[\alpha:\alpha_j\big],0 \Big)
& \quad \text{if $\epsilon(i) = \pm$} \\ \\
\ \ \big[ \alpha : \alpha_i \big] & \quad \text{if $\epsilon(i) = \mp$} 
\end{array} \right. \]

\bigskip
\noindent
where $a_{ij}$ are the entries of the Cartan Matrix
of the root system and
where $\epsilon(i)$ is $+$ (resp. $-$) if $i$ is a sink
(resp. source) of the Dynkin graph. 
Direct computation shows that $\tau_+$ and $\tau_-$
are involutive and that they preserve the set of 
almost positive roots. A more subtle attribute enjoyed
by these involutions is that every positive root can be obtained
by applying, in alternation, $\tau_+$ and $\tau_-$
to a unique negative simple root. Thus the
set of negative simple roots may be viewed
as a fundamental domain with respect to the action of
the $\tau$'s; the orbit attached to a negative
simple root $-\alpha_i$ is called the $\tau$-{\it chain}
through $-\alpha_i$.

\bigskip
\noindent
The bijection $X: \Phi_{\geq -1} \longrightarrow 
\mathcal{X}$ between almost positive roots and 
cluster variable is designed with the property that,
after specializing the coefficients $c \in {\bf c}$ of
the cluster algebra to $1$,
the involutions $\tau_+$ and $\tau_-$ extend to
algebra automorphisms of the cluster algebra
defined by the recipe $\tau_{\pm} \ X[ \alpha ] =
X[ \tau_{\pm} \ \alpha ]$ along with the initial
conditions $X[-\alpha_i] = x_i$ and
$X[\alpha_i] = x_i '$ where $x_i'$ is the cluster variable 
that $x_i$ exchanges with.

\bigskip
\noindent
Consider the example
of $\Bbb{G}(3,8)$ and $-\alpha_4$ of the 
$E_8$ root system. 
In this case $X[ -\alpha_4] = B$ and
from the exchange relation  

\[ B \cdot B^{\sigma} \ = \ \Delta^{258} \ Y^{234678} \ A 
\ + \ \Big(\Delta^{128}\Big)^2 \Delta^{567} \Delta^{678} \Delta^{234}
\Delta^{345}. \]

\bigskip
\noindent
we have $X[ \alpha_4 ] = B^{\sigma} $. Notice that the cluster
variables $\Delta^{258}$, $Y^{234678}$, and 
$A$ occurring on the 
right hand side of the above exchange relation
are, respectively, $X[-\alpha_2]$, $X[-\alpha_3]$,
and $X[-\alpha_5]$; these cluster variables exchange with
$X[\alpha_2] = X^{123467}$, $X[\alpha_3] =
X^{123458}$, and $X[\alpha_5] = Y^{235678}$.
Specialize the boundary Pl\"ucker coordinates to $1$ and
apply the automorphism $\tau_-$ to the exchange relation.

\[ \tau_- \Big( \ X[-\alpha_4] \ X[\alpha_4] \ \Big) \ = \
\tau_- \Big( \ X[-\alpha_2] \ X[-\alpha_3] \ X[-\alpha_5] 
\ + \ 1 \ \Big) \]   

\[ X[-\alpha_4] \ X \big[ \alpha_2 + \alpha_3 + \alpha_4
+ \alpha_5 \big] 
\ = \ X[\alpha_2] \ X[\alpha_3] \ X[\alpha_5] \ + \ 1 \]

\begin{equation} B \ X \big[ \alpha_2 + \alpha_3 + \alpha_4
+ \alpha_5 \big]
\ = \ X^{123467} \ X^{123458} \ Y^{235678} \ + \ 1 \end{equation}

\bigskip
\noindent
Any cluster variable $X$ of the Grassmannian $\Bbb{G}(k,n)$ 
is homogeneous
with respect to the left action of the 
maximal torus $H \subset SL_n(\Bbb{C})$; meaning
that for any diagonal $n \times n$ matrix
$h = \diag \big(h_1, \dots, h_n\big)$ and any
point $p$, viewed as a $k \times n$ matrix,
$X ( p \cdot h) \ = \ h_1^{i_1} \cdots h_n^{i_n} \ X(p)$.
The vector $(i_1,\dots,i_n)$ of exponents 
is called the {\it toral weight} of $X$.
The toral weights of $B$, $X^{123467}$, $X^{123458}$,
and $Y^{235678}$ are, respectively,

\[ \begin{array}{cc}
& (1,2,1,1,1,1,1,1) \\
& (1,1,1,1,0,1,1,0) \\
& (1,1,1,1,1,0,0,1) \\
& (0,1,1,0,1,1,1,1) \end{array} \]

\bigskip
\noindent
Adding the weights of the left hand side of equation (6) and
subtracting the toral weight of $B$  
it follows that the toral weight of the cluster
variable $X \big[ \alpha_2 + \alpha_3 + \alpha_4
+ \alpha_5 \big]$ must be $(1,1,2,1,1,1,1,1)$.
There are only three cluster variables which
have this particular weight: $A^{\rho^2}$,
$B^{\rho}$, and $B^{\sigma \rho^5}$.

\bigskip
\noindent
Applying the Lie derivation $E_6$ to the left and right
hand sides of equation (6) we obtain

\[ B \ E_6 \Big( X \big[ \alpha_2 + \alpha_3 + \alpha_4
+ \alpha_5 \big] \Big) \ + \ X \big[ \alpha_2 + \alpha_3 + \alpha_4
+ \alpha_5 \big] \ E_6 \big( B \big) \ = \]

\bigskip
\[ = \ \begin{array}{cc}
& E_6 \big( X^{123467} \big) \ X^{123458} \ Y^{235678} \\
& + \\
& E_6 \big( X^{123458} \big) \ X^{123467} \ Y^{235678} \\
& + \\   
& E_6 \big( Y^{235678} \big) \ X^{123467} \ X^{123458} \end{array} 
\]

\bigskip
\noindent
The cluster variables $B$, $X^{123467}$, and $Y^{235678}$
are annihilated by $E_6$ because, as functions of
the column vectors $v_6$ and $v_7$, all three
depend only on the plane (line projectively) containing
$v_6$ and $v_7$. The remaining term $X^{123458}$ is also
killed by $E_6$ since it is altogether independent
of $v_6$ and $v_7$. These considerations 
force $E_6 \Big( X \big[ \alpha_2 + \alpha_3 + \alpha_4
+ \alpha_5 \big] \Big)$ to vanish. Of the three possibilities,
only $B^{\sigma \rho^5}$ satisfies this condition,
and hence $X \big[ \alpha_2 + \alpha_3 + \alpha_4
+ \alpha_5 \big] = B^{\sigma \rho^5}$. These rest
of the correspondence is worked out in this 
Lie-theoretic manner.

\bigskip
\begin{Prop}\label{ch2,prop.G6} The following four tables record
the correspondence
between the {\it almost positive roots} of
$D_4$, $E_6$, and $E_8$ and
the cluster variables of $\Bbb{G}(3,6)$,
$\Bbb{G}(3,7)$, and $\Bbb{G}(3,8)$.
Each almost positive root (left hand column)
encodes the vector of exponents of the
monomial denominator of the corresponding
cluster variable (right hand column)
when expressed as a Laurent polynomial in
the Dynkin cluster. 
\end{Prop}

\[ \begin{array}{cc} \begin{tabular}{|l|c|} 
\hline & \\ 
$-\alpha_1$ & $\Delta^{145}$ \\ 
\hline & \\
$\ \ \alpha_1$ & $\Delta^{236}$ \\
\hline & \\
$\ \ \alpha_1 + \alpha_2$ & $\Delta^{356}$ \\
\hline & \\
$\ \ \alpha_2 + \alpha_3 + \alpha_4$ & 
$\Delta^{124}$ \\
\hline \end{tabular} 
& \begin{tabular}{|l|l|} 
\hline & \\
$ -\alpha_2$ & $Y^{123456}$ \\
\hline & \\
$\ \ \alpha_2$ & $\Delta^{135}$ \\
\hline & \\
$\ \ \alpha_1 + \alpha_2
+ \alpha_3 + \alpha_4$ & $\Delta^{246}$ \\
\hline & \\
$\ \ \alpha_1 + 2 \alpha_2
+ \alpha_3 + \alpha_4$ & $X^{123456}$ \\
\hline \end{tabular} \\ \\ \\
\begin{tabular}{|l|c|} 
\hline & \\
$ -\alpha_3$ & $\Delta^{136}$ \\
\hline & \\ 
$\ \ \alpha_3$ & $\Delta^{245}$ \\
\hline & \\
$\ \ \alpha_2 + \alpha_3$ & $\Delta^{125}$ \\
\hline & \\
$\ \ \alpha_1 + \alpha_2 + \alpha_4 $ &
$\Delta^{346}$ \\
\hline \end{tabular} 
& \begin{tabular}{|l|c|} 
\hline & \\
$ -\alpha_4$ & $\Delta^{235}$ \\
\hline & \\
$\ \ \alpha_4$ & $\Delta^{146}$ \\
\hline & \\
$\ \ \alpha_2 + \alpha_4$ & $\Delta^{134}$ \\
\hline & \\
$\ \ \alpha_1 + \alpha_2 + \alpha_3$ 
& $\Delta^{256}$ \\
\hline \end{tabular} \end{array}
\]
\begin{figure}[h]
\begin{center}
\caption{$D_4$ Correspondence}
\label{D1}
\end{center}
\end{figure}

\[
\begin{array}{cc}
\begin{tabular}{|l|c|}
\hline & \\
$ -\alpha_1$ & $\Delta^{237}$ \\
\hline & \\
$ \ \ \alpha_1$ & $\Delta^{156}$ \\
\hline & \\
$ \ \ \alpha_1 + \alpha_3$ & $\Delta^{124}$ \\
\hline & \\
$ \ \ \alpha_3 + \alpha_4$ & $\Delta^{467}$ \\
\hline & \\
$ \ \ \alpha_2 + \alpha_4 + \alpha_5$ & $\Delta^{356}$ \\
\hline & \\
$ \ \ \alpha_2 + \alpha_4 + \alpha_5 + \alpha_6$ 
& $\Delta^{235}$ \\
\hline & \\
$ \ \ \alpha_3 + \alpha_4 + \alpha_5 + \alpha_6$ 
& $\Delta^{457}$ \\
\hline & \\
$ \ \ \alpha_1 + \alpha_3 + \alpha_4 + \alpha_5$ 
& $\Delta^{134}$ \\
\hline & \\
$ \ \ \alpha_1 + \alpha_2 + \alpha_3 + \alpha_4$ 
& $\Delta^{126}$ \\
\hline & \\
$ \ \ \alpha_2 + \alpha_3 + \alpha_4$ 
& $\Delta^{267}$ \\
\hline & \\
$ \ \ \alpha_4 + \alpha_5$ 
& $\Delta^{137}$ \\
\hline & \\
$ \ \ \alpha_5 + \alpha_6$ 
& $\Delta^{157}$ \\
\hline & \\
$ \ \ \alpha_6$ 
& $\Delta^{245}$ \\
\hline & \\
$ -\alpha_6$ 
& $\Delta^{346}$ \\
\hline
\end{tabular} 
& \begin{tabular}{|l|l|}
\hline & \\
$- \alpha_3$ & $Y^{123567}$ \\
\hline & \\
$\ \ \alpha_3$ & $\Delta^{247}$ \\
\hline & \\
$\ \ \alpha_1 + \alpha_3 + \alpha_4$ 
& $\Delta^{146}$ \\
\hline & \\
$\ \ \alpha_1 + \alpha_2 + \alpha_3
+ \alpha_4 + \alpha_5$ & $X^{123456}$ \\
\hline & \\
$\ \ \alpha_2 + \alpha_3 + 2 \alpha_4 
+ \alpha_5 + \alpha_6$ & $X^{234567}$ \\
\hline & \\
$\ \ \alpha_2 + \alpha_3 + 2 \alpha_4 
+ 2 \alpha_5 + \alpha_6$ & $\Delta^{357}$ \\
\hline & \\
$\ \ \alpha_1 + \alpha_2 + \alpha_3 
+ 2 \alpha_4 + 2 \alpha_5 + 
\alpha_6$ & $\Delta^{135}$ \\
\hline & \\
$\ \ \alpha_1 + \alpha_2 + 2 \alpha_3 
+ 2 \alpha_4 + \alpha_5 + 
\alpha_6$ & $X^{124567}$ \\
\hline & \\
$\ \ \alpha_1 + \alpha_2 + 2 \alpha_3 
+ 2 \alpha_4 + \alpha_5$ & $X^{123467}$ \\
\hline & \\
$\ \ \alpha_1 + \alpha_2 + \alpha_3 
+ 2 \alpha_4 + \alpha_5$ & $\Delta^{136}$ \\
\hline & \\
$\ \ \alpha_2 + \alpha_3 + \alpha_4 
+ \alpha_5 + \alpha_6$ & $\Delta^{257}$ \\
\hline & \\
$\ \ \alpha_4 
+ \alpha_5 + \alpha_6$ & $Y^{123457}$ \\
\hline & \\
$\ \ \alpha_5$ & $Y^{134567}$ \\
\hline & \\
$ -\alpha_5$ & $\Delta^{246}$ \\
\hline
\end{tabular} \\ \\ \\ 
\begin{tabular}{|l|c|}
\hline & \\
$- \alpha_2$ & $\Delta^{147}$ \\
\hline & \\
$\ \ \alpha_2$ & $\Delta^{256}$ \\
\hline & \\
$\ \ \alpha_2 + \alpha_4$ & $\Delta^{236}$ \\
\hline & \\
$\ \ \alpha_3 + \alpha_4 + \alpha_5$ 
& $\Delta^{347}$ \\
\hline & \\
$\ \ \alpha_1 + \alpha_3 + \alpha_4
+ \alpha_5 + \alpha_6$ & $\Delta^{145}$ \\
\hline & \\
$\ \ \alpha_1 + \alpha_2 + \alpha_3 + \alpha_4
+ \alpha_5 + \alpha_6$ & $\Delta^{125}$ \\
\hline & \\
$\ \ \alpha_2 + \alpha_3 + 2 \alpha_4
+ \alpha_5$ & $\Delta^{367}$ \\
\hline 
\end{tabular}
& \begin{tabular}{|l|c|}
\hline & \\
$- \alpha_4$ & $Y^{124567}$ \\
\hline & \\
$\ \ \alpha_4$ & $Y^{123467}$ \\
\hline & \\
$\ \ \alpha_2 + \alpha_3 + \alpha_4
+ \alpha_5$ & $Y^{234567}$ \\
\hline & \\
$\ \ \alpha_1 + \alpha_2 + \alpha_3 
+ 2 \alpha_4 + \alpha_5 + \alpha_6$ 
& $Y^{123456}$ \\
\hline & \\
$\ \ \alpha_1 + \alpha_2 + 2\alpha_3
+ 2\alpha_4 + 2\alpha_5 + \alpha_6$ & $X^{123457}$ \\
\hline & \\
$\ \ \alpha_1 + \alpha_2 + 2\alpha_3 + 3\alpha_4
+ 2\alpha_5 + \alpha_6$ & $X^{134567}$ \\
\hline & \\
$\ \ \alpha_1 + 2\alpha_2 + 2\alpha_3
+ 3\alpha_4 + 2\alpha_5 + \alpha_6$ & $X^{123567}$ \\
\hline \end{tabular} \end{array} \]

\begin{center}
\begin{tabular}{|l|c|}
\hline & \\
$ -\alpha_1$ & $\Delta^{348}$ \\
\hline & \\
$\ \ \alpha_1$ & $\Delta^{267}$ \\
\hline & \\
$\ \ \alpha_1 + \alpha_3$ & $\Delta^{125}$ \\
\hline & \\
$\ \ \alpha_3 + \alpha_4$ & $\Delta^{158}$ \\
\hline & \\
$\ \ \alpha_2 + \alpha_4 + \alpha_5$ & $\Delta^{367}$ \\
\hline & \\
$\ \ \alpha_2 + \alpha_4 + \alpha_5 + \alpha_6$ & $\Delta^{347}$ \\
\hline & \\
$\ \ \alpha_3 + \alpha_4 + \alpha_5 + \alpha_6 
+ \alpha_7$ & $\Delta^{458}$ \\
\hline & \\
$\ \ \alpha_1 + \alpha_3 + \alpha_4
+ \alpha_5 + \alpha_6 + \alpha_7 + \alpha_8$ & $\Delta^{256}$ \\
\hline & \\
$\ \ \alpha_1 + \alpha_2 + \alpha_3 + \alpha_4 + \alpha_5 
+ \alpha_6 + \alpha_7 + \alpha_8$ & $\Delta^{126}$  \\
\hline & \\
$\ \ \alpha_2 + \alpha_3 + 2 \alpha_4 
+ \alpha_5 + \alpha_6 + \alpha_7$ & $\Delta^{148}$ \\
\hline & \\
$\ \ \alpha_2 + \alpha_3 + 2 \alpha_4 
+ 2 \alpha_5 + \alpha_6$ & $\Delta^{378}$ \\
\hline & \\
$\ \ \alpha_1 + \alpha_2 + \alpha_3 +
2 \alpha_4 + 2 \alpha_5 + \alpha_6$ & $\Delta^{237}$ \\
\hline & \\
$\ \ \alpha_1 + \alpha_2 + 2 \alpha_3 +
2 \alpha_4 + \alpha_5 + \alpha_6 + \alpha_7$ & $\Delta^{145}$ \\
\hline & \\
$\ \ \alpha_1 + \alpha_2 + 2 \alpha_3 +
2 \alpha_4 + \alpha_5 + \alpha_6 + 
\alpha_7 + \alpha_8$ & $\Delta^{156}$ \\
\hline & \\
$\ \ \alpha_1 + \alpha_2 + \alpha_3 +
2 \alpha_4 + 2 \alpha_5 + \alpha_6 + 
\alpha_7 + \alpha_8$ & $\Delta^{236}$ \\
\hline & \\ 
$\ \ \alpha_2 + \alpha_3 +
2 \alpha_4 + 2 \alpha_5 + 2 \alpha_6 + 
\alpha_7$ & $\Delta^{478}$ \\
\hline \end{tabular}
\end{center}
\begin{figure}[h]
\begin{center}
\caption{$E_8$ Correspondence}
\label{E1}
\end{center}
\end{figure}

\begin{center}
\begin{tabular}{|l|c|}
\hline & \\
$ -\alpha_2$ & $\Delta^{258}$ \\
\hline & \\
$\ \ \alpha_2$ & $X^{123467}$ \\
\hline & \\
$\ \ \alpha_2 + \alpha_4$ & $Y^{134678}$ \\
\hline & \\
$\ \ \alpha_3 + \alpha_4 + \alpha_5$ & $\Delta^{358}$ \\
\hline & \\
$\ \ \alpha_1 + \alpha_3 + \alpha_4 +
\alpha_5 + \alpha_6$ & $\Delta^{257}$ \\
\hline & \\
$\ \ \alpha_1 + \alpha_2 + \alpha_3 + \alpha_4 +
\alpha_5 + \alpha_6 + \alpha_7$ & $X^{124567}$ \\
\hline & \\
$\ \ \alpha_2 + \alpha_3 + 2 \alpha_4 + \alpha_5 
+ \alpha_6 + \alpha_7 + \alpha_8$ & $Y^{134568}$ \\
\hline & \\
$\ \ \alpha_2 + \alpha_3 + 2 \alpha_4
+ 2 \alpha_5 + \alpha_6 + \alpha_7 + \alpha_8$ & $\Delta^{368}$ \\
\hline & \\
$\ \ \alpha_1 + \alpha_2 + \alpha_3 + 2 \alpha_4 + 2 \alpha_5 
+ 2 \alpha_6 + \alpha_7$ & $\Delta^{247}$  \\
\hline & \\
$\ \ \alpha_1 + \alpha_2 + 2 \alpha_3 
+ 2 \alpha_4 + 2 \alpha_5 + 2 \alpha_6 +
\alpha_7$ & $X^{124578}$ \\
\hline & \\
$\ \ \alpha_1 + \alpha_2 + 2 \alpha_3 
+ 3 \alpha_4 + 2 \alpha_5 + \alpha_6 +
\alpha_7 + \alpha_8$ & $Y^{123568}$ \\
\hline & \\
$\ \ \alpha_1 + 2 \alpha_2 + 2 \alpha_3 +
3 \alpha_4 + 2 \alpha_5 + \alpha_6 +
\alpha_7 + \alpha_8$ & $\Delta^{136}$ \\
\hline & \\
$\ \ \alpha_1 + 2 \alpha_2 + 2 \alpha_3 +
3 \alpha_4 + 2 \alpha_5 + 2 \alpha_6 + \alpha_7$ & $\Delta^{147}$ \\
\hline & \\
$\ \ \alpha_1 + \alpha_2 + 2 \alpha_3 +
3 \alpha_4 + 3 \alpha_5 + 2 \alpha_6 + 
\alpha_7$ & $X^{234578}$ \\
\hline & \\
$\ \ \alpha_1 + \alpha_2 + 2 \alpha_3 +
3 \alpha_4 + 3 \alpha_5 + 2 \alpha_6 + 
\alpha_7 + \alpha_8$ & $X^{235678}$ \\
\hline & \\ 
$\ \ \alpha_1 + 2 \alpha_2 +
2 \alpha_3 + 3 \alpha_4 + 2 \alpha_5 + 
2 \alpha_6 + 2 \alpha_7 + \alpha_8$ & $\Delta^{146}$ \\
\hline \end{tabular}
\end{center}
\begin{figure}[h]
\begin{center}
\caption{$E_8$ Correspondence}
\label{E2}
\end{center}
\end{figure}

\begin{center}
\begin{tabular}{|l|c|}
\hline & \\
$ -\alpha_3$ & $Y^{234678}$ \\
\hline & \\
$\ \ \alpha_3$ & $X^{123458}$ \\
\hline & \\
$\ \ \alpha_1 + \alpha_3 + \alpha_4$ & $Y^{125678}$ \\
\hline & \\
$\ \ \alpha_1 + \alpha_2 + \alpha_3 +
\alpha_4 + \alpha_5$ & $X^{123567}$ \\
\hline & \\
$\ \ \alpha_2 + \alpha_3 + 2 \alpha_4 +
\alpha_5 + \alpha_6$ & $Y^{134578}$ \\
\hline & \\
$\ \ \alpha_2 + \alpha_3 + 2 \alpha_4 + 2 \alpha_5 +
\alpha_6 + \alpha_7$ & $Y^{345678}$ \\
\hline & \\
$\ \ \alpha_1 + \alpha_2 +  \alpha_3 + 2 \alpha_4 
+ 2 \alpha_5 + 2 \alpha_6 + \alpha_7 + \alpha_8$ & $Y^{234567}$ \\
\hline & \\
$\ \ \alpha_1 + \alpha_2 + 2 \alpha_3
+ 2 \alpha_4 + 2 \alpha_5 + 2 \alpha_6 + 2 \alpha_7 +
\alpha_8$ & $X^{124568}$ \\
\hline & \\
$\ \ \alpha_1 + \alpha_2 + 2 \alpha_3 + 3 \alpha_4 + 2 \alpha_5 
+ 2 \alpha_6 + 2 \alpha_7 + \alpha_8$ & $Y^{124568}$  \\
\hline & \\
$\ \ \alpha_1 + 2 \alpha_2 + 2 \alpha_3 
+ 3 \alpha_4 + 3 \alpha_5 + 2 \alpha_6 +
\alpha_7 + \alpha_8$ & $X^{123678}$ \\
\hline & \\
$\ \ \alpha_1 + 2 \alpha_2 + 2 \alpha_3 
+ 4 \alpha_4 + 3 \alpha_5 + 2 \alpha_6 +
\alpha_7 + \alpha_8$ & $Y^{123478}$ \\
\hline & \\
$\ \ \alpha_1 + 2 \alpha_2 + 3 \alpha_3 +
4 \alpha_4 + 3 \alpha_5 + 2 \alpha_6 +
\alpha_7$ & $X^{134578}$ \\
\hline & \\
$\ \ 2 \alpha_1 + 2 \alpha_2 + 3 \alpha_3 +
4 \alpha_4 + 3 \alpha_5 + 2 \alpha_6 + \alpha_7 +
\alpha_8$ & $Y^{123567}$ \\
\hline & \\
$\ \ 2 \alpha_1 + 2 \alpha_2 + 3 \alpha_3 +
4 \alpha_4 + 3 \alpha_5 + 2 \alpha_6 + 
2 \alpha_7 + \alpha_8$ & $Y^{123456}$ \\
\hline & \\
$\ \ \alpha_1 + 2 \alpha_2 + 3 \alpha_3 +
4 \alpha_4 + 3 \alpha_5 + 3 \alpha_6 + 
2 \alpha_7 + \alpha_8$ & $X^{145678}$ \\
\hline & \\ 
$\ \ \alpha_1 + 2 \alpha_2 +
2 \alpha_3 + 4 \alpha_4 + 4 \alpha_5 + 
3 \alpha_6 + 2 \alpha_7 + \alpha_8$ & $X^{234678}$ \\
\hline \end{tabular}
\end{center}
\begin{figure}[h]
\begin{center}
\caption{$E_8$ Correspondence}
\label{E3}
\end{center}
\end{figure}

\begin{center}
\begin{tabular}{|l|c|}
\hline & \\
$ -\alpha_4$ & $B$ \\
\hline & \\
$\ \ \alpha_4$ & $B^\sigma$ \\
\hline & \\
$\ \ \alpha_2 + \alpha_3 + \alpha_4 +
\alpha_5$ & $B^{\sigma \rho^5}$ \\
\hline & \\
$\ \ \alpha_1 + \alpha_2 + \alpha_3 +
2 \alpha_4 + \alpha_5 + \alpha_6$ & $B^{\rho^5}$ \\
\hline & \\
$\ \ \alpha_1 + \alpha_2 + 2 \alpha_3 +
2 \alpha_4 + 2 \alpha_5 +
\alpha_6 + \alpha_7$ & $B^{\rho^3}$ \\
\hline & \\
$\ \ \alpha_1 + \alpha_2 + 2 \alpha_3 + 3 \alpha_4 +
2 \alpha_5 + 2 \alpha_6 + \alpha_7 + \alpha_8$ & 
$B^{\sigma \rho^3}$ \\
\hline & \\
$\ \ \alpha_1 + 2 \alpha_2 +  2 \alpha_3 + 3 \alpha_4 
+ 3 \alpha_5 + 2 \alpha_6 + 2 \alpha_7 + \alpha_8$ & 
$B^{\sigma \rho^2}$ \\
\hline & \\
$\ \ \alpha_1 + 2 \alpha_2 + 2 \alpha_3
+ 4 \alpha_4 + 3 \alpha_5 + 3 \alpha_6 + 2 \alpha_7 +
\alpha_8$ & $B^{\rho^2}$ \\
\hline & \\
$\ \ \alpha_1 + 2 \alpha_2 + 3 \alpha_3 + 4 \alpha_4 + 4 \alpha_5 
+ 3 \alpha_6 + 2 \alpha_7 + \alpha_8$ & $B^{\rho^6}$  \\
\hline & \\
$\ \ 2 \alpha_1 + 2 \alpha_2 + 3 \alpha_3 
+ 5 \alpha_4 + 4 \alpha_5 + 3 \alpha_6 +
2 \alpha_7 + \alpha_8$ & $B^{\sigma \rho^6}$ \\
\hline & \\
$\ \ 2 \alpha_1 + 3 \alpha_2 + 4 \alpha_3 
+ 5 \alpha_4 + 4 \alpha_5 + 3 \alpha_6 +
2 \alpha_7 + \alpha_8$ & $B^{\sigma \rho^7}$ \\
\hline & \\
$\ \ 2 \alpha_1 + 3 \alpha_2 + 4 \alpha_3 +
6 \alpha_4 + 4 \alpha_5 + 3 \alpha_6 +
2 \alpha_7 + \alpha_8$ & $B^{\rho^7}$ \\
\hline & \\
$\ \ 2 \alpha_1 + 3 \alpha_2 + 4 \alpha_3 +
6 \alpha_4 + 5 \alpha_5 + 3 \alpha_6 + 2 \alpha_7 +
\alpha_8$ & $B^{\rho}$ \\
\hline & \\
$\ \ 2 \alpha_1 + 3 \alpha_2 + 4 \alpha_3 +
6 \alpha_4 + 5 \alpha_5 + 4 \alpha_6 + 
2 \alpha_7 + \alpha_8$ & $B^{\sigma \rho}$ \\
\hline & \\
$\ \ 2 \alpha_1 + 3 \alpha_2 + 4 \alpha_3 +
6 \alpha_4 + 5 \alpha_5 + 4 \alpha_6 + 
3 \alpha_7 + \alpha_8$ & $B^{\sigma \rho^4}$ \\
\hline & \\ 
$\ \ 2 \alpha_1 + 3 \alpha_2 +
4 \alpha_3 + 6 \alpha_4 + 5 \alpha_5 + 
4 \alpha_6 + 3 \alpha_7 + 2 \alpha_8$ & $B^{\rho^4}$ \\
\hline \end{tabular}
\end{center}
\begin{figure}[h]
\begin{center}
\caption{$E_8$ Correspondence}
\label{E4}
\end{center}
\end{figure}

\begin{center}
\begin{tabular}{|l|c|}
\hline & \\
$ -\alpha_5$ & $A$ \\
\hline & \\
$\ \ \alpha_5$ & $Y^{235678}$ \\
\hline & \\
$\ \ \alpha_4 + \alpha_5 + \alpha_6$ & $Y^{234578}$ \\
\hline & \\
$\ \ \alpha_2 + \alpha_3 + \alpha_4 +
2 \alpha_5 + \alpha_6 + \alpha_7$ & $A^{\rho^3}$ \\
\hline & \\
$\ \ \alpha_1 + \alpha_2 + \alpha_3 +
2 \alpha_4 + \alpha_5 +
\alpha_6 + \alpha_7 + \alpha_8$ & $A^{\rho^5}$ \\
\hline & \\
$\ \ \alpha_1 + \alpha_2 + 2 \alpha_3 + 2 \alpha_4 +
2 \alpha_5 + \alpha_6 + \alpha_7 + \alpha_8$ & 
$X^{123568}$ \\
\hline & \\
$\ \ \alpha_1 + \alpha_2 +  2 \alpha_3 + 3 \alpha_4 
+ 2 \alpha_5 + 2 \alpha_6 + \alpha_7$ & 
$Y^{124578}$ \\
\hline & \\
$\ \ \alpha_1 + 2 \alpha_2 + 2 \alpha_3
+ 3 \alpha_4 + 3 \alpha_5 + 2 \alpha_6 + \alpha_7$ & $A^{\rho^6}$ \\
\hline & \\
$\ \ \alpha_1 + 2 \alpha_2 + 2 \alpha_3 + 4 \alpha_4 + 3 \alpha_5 
+ 2 \alpha_6 + \alpha_7 + \alpha_8$ & $A^{\rho^2}$  \\
\hline & \\
$\ \ \alpha_1 + 2 \alpha_2 + 3 \alpha_3 
+ 4 \alpha_4 + 3 \alpha_5 + 2 \alpha_6 +
2 \alpha_7 + \alpha_8$ & $X^{134568}$ \\
\hline & \\
$\ \ 2 \alpha_1 + 2 \alpha_2 + 3 \alpha_3 
+ 4 \alpha_4 + 3 \alpha_5 + 3 \alpha_6 +
2 \alpha_7 + \alpha_8$ & $Y^{124567}$ \\
\hline & \\
$\ \ 2 \alpha_1 + 2 \alpha_2 + 3 \alpha_3 +
4 \alpha_4 + 4 \alpha_5 + 3 \alpha_6 +
2 \alpha_7 + \alpha_8$ & $A^{\rho}$ \\
\hline & \\
$\ \ \alpha_1 + 2 \alpha_2 + 3 \alpha_3 +
5 \alpha_4 + 4 \alpha_5 + 3 \alpha_6 + 2 \alpha_7 +
\alpha_8$ & $A^{\rho^7}$ \\
\hline & \\
$\ \ \alpha_1 + 3 \alpha_2 + 3 \alpha_3 +
5 \alpha_4 + 4 \alpha_5 + 3 \alpha_6 + 
2 \alpha_7 + \alpha_8$ & $X^{134678}$ \\
\hline & \\
$\ \ 2 \alpha_1 + 3 \alpha_2 + 3 \alpha_3 +
5 \alpha_4 + 4 \alpha_5 + 3 \alpha_6 + 
2 \alpha_7 + \alpha_8$ & $Y^{123467}$ \\
\hline & \\ 
$\ \ 2 \alpha_1 + 2 \alpha_2 +
4 \alpha_3 + 5 \alpha_4 + 4 \alpha_5 + 
3 \alpha_6 + 2 \alpha_7 + \alpha_8$ & $A^{\rho^4}$ \\
\hline \end{tabular}
\end{center}
\begin{figure}[h]
\begin{center}
\caption{$E_8$ Correspondence}
\label{E5}
\end{center}
\end{figure}

\begin{center}
\begin{tabular}{|l|c|}
\hline & \\
$ -\alpha_6$ & $Y^{135678}$ \\
\hline & \\
$\ \ \alpha_6$ & $X^{123457}$ \\
\hline & \\
$\ \ \alpha_5 + \alpha_6 + \alpha_7$ & $Y^{245678}$ \\
\hline & \\
$\ \ \alpha_4 + \alpha_5 + \alpha_6 +
\alpha_7 + \alpha_8$ & $Y^{234568}$ \\
\hline & \\
$\ \ \alpha_2 + \alpha_3 + \alpha_4 +
\alpha_5 + \alpha_6 +
\alpha_7 + \alpha_8$ & $X^{123468}$ \\
\hline & \\
$\ \ \alpha_1 + \alpha_2 + \alpha_3 + 2 \alpha_4 +
\alpha_5 + \alpha_6 + \alpha_7$ & 
$Y^{124678}$ \\
\hline & \\
$\ \ \alpha_1 + \alpha_2 +  2 \alpha_3 + 2 \alpha_4 
+ 2 \alpha_5 + \alpha_6$ & 
$X^{123578}$ \\
\hline & \\
$\ \ \alpha_1 + \alpha_2 + 2 \alpha_3
+ 3 \alpha_4 + 2 \alpha_5 + \alpha_6$ & $Y^{123578}$ \\
\hline & \\
$\ \ \alpha_1 + 2 \alpha_2 + 2 \alpha_3 + 3 \alpha_4 + 2 \alpha_5 
+ \alpha_6 + \alpha_7$ & $X^{134567}$  \\
\hline & \\
$\ \ \alpha_1 + 2 \alpha_2 + 2 \alpha_3 
+ 3 \alpha_4 + 2 \alpha_5 + 2 \alpha_6 + 
\alpha_7 + \alpha_8$ & $Y^{134567}$ \\
\hline & \\
$\ \ \alpha_1 + \alpha_2 + 2 \alpha_3 
+ 3 \alpha_4 + 3 \alpha_5 + 2 \alpha_6 +
2 \alpha_7 + \alpha_8$ & $X^{234568}$ \\
\hline & \\
$\ \ \alpha_1 + \alpha_2 + 2 \alpha_3 +
3 \alpha_4 + 3 \alpha_5 + 3 \alpha_6 +
2 \alpha_7 + \alpha_8$ & $X^{245678}$ \\
\hline & \\
$\ \ \alpha_1 + 2 \alpha_2 + 2 \alpha_3 +
3 \alpha_4 + 3 \alpha_5 + 3 \alpha_6 + 2 \alpha_7 +
\alpha_8$ & $X^{124678}$ \\
\hline & \\
$\ \ \alpha_1 + 2 \alpha_2 + 2 \alpha_3 +
4 \alpha_4 + 3 \alpha_5 + 2 \alpha_6 + 
2 \alpha_7 + \alpha_8$ & $Y^{123468}$ \\
\hline & \\
$\ \ \alpha_1 + 2 \alpha_2 + 3 \alpha_3 +
4 \alpha_4 + 3 \alpha_5 + 2 \alpha_6 + 
\alpha_7 + \alpha_8$ & $X^{135678}$ \\
\hline & \\ 
$\ \ 2 \alpha_1 + 2 \alpha_2 +
3 \alpha_3 + 4 \alpha_4 + 3 \alpha_5 + 
2 \alpha_6 + \alpha_7$ & $Y^{123457}$ \\
\hline \end{tabular}
\end{center}
\begin{figure}[h]
\begin{center}
\caption{$E_8$ Correspondence}
\label{E6}
\end{center}
\end{figure}

\begin{center}
\begin{tabular}{|l|c|}
\hline & \\
$ -\alpha_7$ & $\Delta^{357}$ \\
\hline & \\
$\ \ \alpha_7$ & $Y^{145678}$ \\
\hline & \\
$\ \ \alpha_6 + \alpha_7 + \alpha_8$ & $X^{123456}$ \\
\hline & \\
$\ \ \alpha_5 + \alpha_6 + \alpha_7 +
\alpha_8$ & $\Delta^{268}$ \\
\hline & \\
$\ \ \alpha_4 + \alpha_5 + \alpha_6 +
\alpha_7$ & $\Delta^{248}$ \\
\hline & \\
$\ \ \alpha_2 + \alpha_3 + \alpha_4 + \alpha_5 +
\alpha_6 $ & $X^{123478}$ \\
\hline & \\
$\ \ \alpha_1 + \alpha_2 +  \alpha_3 + 2 \alpha_4 
+ \alpha_5 $ & $Y^{123678}$ \\
\hline & \\
$\ \ \alpha_1 + \alpha_2 + 2 \alpha_3
+ 2 \alpha_4 + \alpha_5$ & $\Delta^{135}$ \\
\hline & \\
$\ \ \alpha_1 + \alpha_2 + 2 \alpha_3 + 2 \alpha_4 + \alpha_5 
+ \alpha_6 $ & $\Delta^{157}$  \\
\hline & \\
$\ \ \alpha_1 + \alpha_2 + \alpha_3 
+ 2 \alpha_4 + 2 \alpha_5 + \alpha_6 + 
\alpha_7 + \alpha_8$ & $X^{234567}$ \\
\hline & \\
$\ \ \alpha_2 + \alpha_3 + 2 \alpha_4 
+ 2 \alpha_5 + 2 \alpha_6 + \alpha_7 + \alpha_8$ & $X^{345678}$ \\
\hline & \\
$\ \ \alpha_2 + \alpha_3 + 2 \alpha_4 +
2 \alpha_5 + 2 \alpha_6 + 2 \alpha_7 +
\alpha_8$ & $\Delta^{468}$ \\
\hline & \\
$\ \ \alpha_1 + \alpha_2 + \alpha_3 +
2 \alpha_4 + 2 \alpha_5 + 2 \alpha_6 + 2 \alpha_7 +
\alpha_8$ & $\Delta^{246}$ \\
\hline & \\
$\ \ \alpha_1 + \alpha_2 + 2 \alpha_3 +
2 \alpha_4 + 2 \alpha_5 + 2 \alpha_6 + 
\alpha_7 + \alpha_8$ & $X^{125678}$ \\
\hline & \\
$\ \ \alpha_1 + \alpha_2 + 2 \alpha_3 +
3 \alpha_4 + 2 \alpha_5 + \alpha_6 + 
\alpha_7$ & $Y^{123458}$ \\
\hline & \\ 
$\ \ \alpha_1 + 2 \alpha_2 +
2 \alpha_3 + 3 \alpha_4 + 2 \alpha_5 + 
\alpha_6$ & $\Delta^{137}$ \\
\hline \end{tabular}
\end{center}
\begin{figure}[h]
\begin{center}
\caption{$E_8$ Correspondence}
\label{E7}
\end{center}
\end{figure}

\begin{center}
\begin{tabular}{|l|c|}
\hline & \\
$ -\alpha_8$ & $\Delta^{457}$ \\
\hline & \\
$\ \ \alpha_8$ & $\Delta^{356}$ \\
\hline & \\
$\ \ \alpha_7 + \alpha_8$ & $\Delta^{168}$ \\
\hline & \\
$\ \ \alpha_6 + \alpha_7$ & $\Delta^{124}$ \\
\hline & \\
$\ \ \alpha_5 + \alpha_6$ & $\Delta^{278}$ \\
\hline & \\
$\ \ \alpha_4 + \alpha_5$ & $\Delta^{238}$ \\
\hline & \\
$\ \ \alpha_2 + \alpha_3 +  \alpha_4$ & $\Delta^{134}$ \\
\hline & \\
$\ \ \alpha_1 + \alpha_2 + \alpha_3 + \alpha_4$ & $\Delta^{167}$ \\
\hline & \\
$\ \ \alpha_1 + \alpha_3 + \alpha_4 + 
\alpha_5$ & $\Delta^{235}$  \\
\hline & \\
$\ \ \alpha_3 + \alpha_4 + \alpha_5 + \alpha_6$ & $\Delta^{578}$ \\
\hline & \\
$\ \ \alpha_2 + \alpha_4 + \alpha_5 + 
\alpha_6 + \alpha_7$ & $\Delta^{467}$ \\
\hline & \\
$\ \ \alpha_2 + \alpha_4 + \alpha_5 +
\alpha_6 + \alpha_7 + \alpha_8$ & $\Delta^{346}$ \\
\hline & \\
$\ \ \alpha_3 + \alpha_4 + \alpha_5 +
\alpha_6 + \alpha_7 + \alpha_8$ & $\Delta^{568}$ \\
\hline & \\
$\ \ \alpha_1 + \alpha_3 + \alpha_4 +
\alpha_5 + \alpha_6 + \alpha_7$ & $\Delta^{245}$ \\
\hline & \\
$\ \ \alpha_1 + \alpha_2 + \alpha_3 +
\alpha_4 + \alpha_5 + \alpha_6 $ & $\Delta^{127}$ \\
\hline & \\ 
$\ \ \alpha_2 + \alpha_3 +
2 \alpha_4 + \alpha_5$ & $\Delta^{138}$ \\
\hline \end{tabular}
\end{center}
\begin{figure}[h]
\begin{center}
\caption{$E_8$ Correspondence}
\label{E8}
\end{center}
\end{figure}
\bigskip

\section{Future Directions}

\bigskip
\noindent
One overriding direction
in the theory of cluster algebras,
posed as a conjecture by {\bf [Ref]},
is to show
that any cluster variable can be expressed as a
Laurent polynomial with {\bf positive} coefficients
in the cluster variables attached to any fixed
cluster. As indicated in section 2, this positivity
property has been verified for the homogeneous
coordinate ring of the Grassmannian $\Bbb{G}(2,n)$.

\bigskip
\noindent
Naturally its is expected that such positive Laurent
decompositions are valid for {\bf any} Grassmannian.
An intermediate step in the proof of this conjecture 
would be to show that
any cluster variable in $\Bbb{C} \Big[ \Bbb{G}(k,n) \Big]$
can be expressed as a positive Laurent polynomial
in the cluster variables arising from a Postnikov
diagram. Ideally the apparatus of Postnikov diagrams
might allow one to read off the positive coefficients
as some type of combinatorial statistic related
to the Postnikov diagram. Indeed this is true
for the quadrilateral Postnikov diagrams described in section 3
when $n= 2k$. As indicated in section 3, the quadrilateral
Postnikov diagram for $n=2k$ is related to a certain
double wiring arrangement studied by {\bf [Ref]}
in the context of {\it totally positive matrices}.
In this case the coefficients are positive integers
counting {\it vertex-disjoint paths} in some related
planar graph. 

\bigskip
\noindent
Another approach to this conjecture is tied to the
following observation. Let $J = \big\{j_1 
< \dots < j_k \big\}$ be a $k$-subset of $[1 \dots n]$
and let $\lambda_J = \big( \lambda_1, \dots, \lambda_k
\big)$ be the partition defined by $\lambda_r 
= j_r - r$. Whenever the short Pl\"ucker relation
 
\[ \Delta^{Iij} \Delta^{Ist} \ = \ \ \Delta^{Iis} \Delta^{Ijt}
\ + \ \Delta^{Iit} \Delta^{Ijs} \] 

\bigskip 
\noindent
is valid, then so is the following relation of Schur polynomials:

\[ s_{\lambda_{Iij}} s_{\lambda_{Ist}} \ = \ \
s_{\lambda_{Iis}} s_{\lambda_{Ijt}} \ + \ 
s_{\lambda_{Iit}} s_{\lambda_{Ijs}}. \]

\bigskip
\noindent
This coincidence indicates that perhaps the cluster algebra
$\Bbb{C}\Big[ \Bbb{G}(k,n) \Big]$ is isomorphic to
a subalgebra of the representation ring of $GL(k+1)$.
This model may reveal that the positivity conjecture
is a manifestation of certain Littlewood-Richardson
tensor decompositions rules.

\bigskip
\section{Acknowledgements}

\bigskip
\noindent
I am grateful for the discussions I had
with Alexander Postnikov, whose
theory of {\it $\pi_{k,n}$-diagrams}
proved so useful to this work.
I would like to acknowledge Sergey 
Fomin who provided the author with
the list of the $\tau$-chains for $E_8$.

\bigskip
\noindent
Many thanks to Pavel Etingof and David Vogan 
who generously hosted me at the mathematics dept.
of MIT where significant
advances for this publication transpired.

\bigskip
\noindent
I emphatically thank my advisor Andrei Zelevinsky
for an introduction to the combinatorics of
cluster algebras and to the border of current research.

\end{document}